\documentclass[12pt]{article}
\usepackage{amsmath,amsfonts,amssymb,amsthm}
\input amssym.def
\begin{document}
\def \Z{\Bbb Z}
\def \C{\Bbb C}
\def \R{\Bbb R}
\def \Q{\Bbb Q}
\def \N{\Bbb N}
\def \bR{\bf R}
\def \D{{\cal{D}}}
\def \E{{\cal{E}}}
\def \S{{\cal{S}}}
\def \R{{\cal{R}}}
\def \Y{{\cal{Y}}}
\def \wt{{\rm wt}}
\def \tr{{\rm tr}}
\def \span{{\rm span}}
\def \Res{{\rm Res}}
\def \End{{\rm End}\;}
\def \Ind {{\rm Ind}}
\def \Irr {{\rm Irr}}
\def \Aut{{\rm Aut}}
\def \Hom{{\rm Hom}}
\def \mod{{\rm mod}}
\def \ann{{\rm Ann}}
\def \ad{{\rm ad}}
\def \rank{{\rm rank}\;}

\def \<{\langle} 
\def \>{\rangle} 
\def \t{\tau }
\def \a{\alpha }
\def \e{\epsilon }
\def \l{\lambda }
\def \L{\Lambda }
\def \g{{\frak{g}}}
\def \h{{\frak{h}}}
\def \k{{\frak{k}}}
\def \sl{{\frak{sl}}}
\def \b{\beta }
\def \c{\chi}
\def \ch{\chi}
\def \cg{\chi_g}
\def \ag{\alpha_g}
\def \ah{\alpha_h}
\def \ph{\psi_h}
\def \be{\begin{equation}\label}
\def \ee{\end{equation}}
\def \bex{\begin{example}\label}
\def \eex{\end{example}}
\def \bl{\begin{lem}\label}
\def \el{\end{lem}}
\def \bt{\begin{thm}\label}
\def \et{\end{thm}}
\def \bp{\begin{prop}\label}
\def \ep{\end{prop}}
\def \br{\begin{rem}\label}
\def \er{\end{rem}}
\def \bc{\begin{coro}\label}
\def \ec{\end{coro}}
\def \bd{\begin{de}\label}
\def \ed{\end{de}}
\def \pf{{\bf Proof. }}

\newcommand{\n}{\:^{\times}_{\times}\:}
\newcommand{\nno}{\nonumber}
\newcommand{\nord}{\mbox{\scriptsize ${\circ\atop\circ}$}}
\newtheorem{thm}{Theorem}[section]
\newtheorem{prop}[thm]{Proposition}
\newtheorem{coro}[thm]{Corollary}
\newtheorem{conj}[thm]{Conjecture}
\newtheorem{example}[thm]{Example}
\newtheorem{lem}[thm]{Lemma}
\newtheorem{rem}[thm]{Remark}
\newtheorem{de}[thm]{Definition}
\newtheorem{hy}[thm]{Hypothesis}
\makeatletter
\@addtoreset{equation}{section}
\def\theequation{\thesection.\arabic{equation}}
\makeatother
\makeatletter

\begin{center}{\Large \bf Nonlocal vertex algebras generated by
formal vertex operators}
\end{center}

\begin{center}
{Haisheng Li\footnote{Partially supported by an NSA grant}\\
Department of Mathematical Sciences, Rutgers University, Camden, NJ 08102\\
and\\
Department of Mathematics, Harbin Normal University, Harbin, China}
\end{center}

\begin{center}
{To James Lepowsky and Robert Wilson
on Their 60th Birthdays}
\end{center}

\begin{abstract}
This is the first paper in a series to study vertex algebra-like
objects arising from infinite-dimensional quantum groups (quantum
affine algebras and Yangians).  In this paper we lay the foundation
for this study.  For any vector space $W$, we study what we call quasi
compatible subsets of $\Hom (W,W((x)))$ and we prove that any maximal
quasi compatible subspace has a natural nonlocal (namely
noncommutative) vertex algebra structure with $W$ as a natural faithful quasi
module in a certain sense and that any quasi compatible subset
generates a nonlocal vertex algebra with $W$ as a quasi module.  
In particular, taking $W$ to be a highest weight module
for a quantum affine algebra we obtain a nonlocal vertex
algebra with $W$ as a quasi module.  We also formulate and study a
notion of quantum vertex algebra and we give general constructions of
nonlocal vertex algebras, quantum vertex algebras and their modules.
\end{abstract}

\section{Introduction}
Vertex (operator) algebras, as a new fundamental class of  algebraic
structures, have been known to have deep connections with various fields 
in both mathematics and physics.  
Among the important examples of vertex operator algebras are
those associated to affine Kac-Moody Lie algebras (including
infinite-dimensional Heisenberg Lie algebras) and the Virasoro (Lie)
algebra, which underline
the algebraic study of the physical Wess-Zumino-Novikov-Witten model
and the minimal models in conformal field theory, respectively.  
A fundamental problem, which was first formulated in \cite{fj} (cf. \cite{efk}),
is to develop an appropriate theory
of quantum vertex (operator) algebras, so that quantum affine algebras
can be naturally associated to quantum vertex (operator) algebras in
the same way that affine Lie algebras are associated with vertex
algebras. 

Roughly speaking, there have been three notions of quantum (operator)
algebra in literature, which are E. Frenkel and Reshetikhin's notion
of deformed chiral algebra (see \cite{efr}), Etingof and Kazhdan's
notion of quantum vertex operator algebra (see \cite{ek}), and
Borcherds' notion of quantum vertex algebra (see \cite{b-qva}).  All
of these works produce interesting and important objects.  On the
other hand, as far as we understand, the previously mentioned problem
is still to be solved.  This is our main concern of this paper.

Starting with this paper in a series, we develop a new theory
eventually to give an explicit solution to this problem.  In this
paper we establish certain conceptual results and as a simple
application we show that quantum affine algebras (through highest
weight modules) in a certain natural way can be associated to certain nonlocal (namely
noncommutative) vertex algebras, 
namely weak axiomatic $G_{1}$-vertex algebras studied in
\cite{li-g1}, or field algebras studied in \cite{kacbook2} (cf. \cite{bk}).
(Nonlocal vertex algebras are analogues of (noncommutative) associative algebras, 
in contrast with that vertex algebras are analogues of commutative
associative algebras.)
This partially solves the problem and provides a machinery for a
complete solution.

Our key (new) idea is that first of all we should study what kind of
vertex algebra-like structures we possibly can get from the generating
functions of quantum affine algebras in the Drinfeld's realization
(see \cite{drinfeld}) on a highest weight module (see \cite{fj}).
This idea naturally came out of our earlier works \cite{li-local},
\cite{li-twisted} and \cite{li-gamma} where we established 
certain conceptual results and made natural
connections between (both untwisted and twisted) affine Lie algebras
and vertex algebras, and between certain quantum torus Lie algebras
and vertex algebras.

Let $W$ be any vector space. Set $\E(W)=\Hom (W,W((x)))$.
For $\alpha(x),\beta(x)\in \E(W)$, $n\in \Z$, define
\begin{eqnarray}\label{e-old-def}
\alpha(x)_{n}\beta(x)
=\Res_{x_{1}}\left( (x_{1}-x)^{n}\alpha(x_{1})\beta(x)
-(-x+x_{1})^{n}\beta(x)\alpha(x_{1})\right).
\end{eqnarray}
A subset $S$ of $\E(W)$ is said to be {\em local} if for any 
$\a(x),\b(x)\in S$,
$$(x_{1}-x_{2})^{k}\a(x_{1})\b(x_{2})=(x_{1}-x_{2})^{k}\b(x_{2})\a(x_{1})$$
for some nonnegative integer $k$.
It was proved (\cite{li-local}, cf. \cite{ll}) that
for any local subset $S$, there exists a unique smallest closed local subspace $\<S\>$ 
containing $S$ and the identity operator $1_{W}$ of $W$ and furthermore
$\<S\>$ is a vertex algebra with $W$ naturally as a module.
Take $W$ to be a highest weight module for an (untwisted) affine Lie algebra $\hat{\g}$
of level $\ell\in \C$.
For $a\in \g$, the associated generating function
$a(x)=\sum_{n\in \Z}(a\otimes t^{n})x^{-n-1}$ is naturally an element of $\E(W)$.
Furthermore, all the generating functions $a(x)$ on $W$
form a local subspace $S_{W}$ of $\E(W)$, where for any $a,b\in \g$,
$$(x_{1}-x_{2})^{2}[a(x_{1}),b(x_{2})]=0.$$
Then we have a vertex algebra $\<S_{W}\>$ (a subspace of $\E(W)$)
with $W$ as a module. It was proved that
the vertex algebra $\<S_{W}\>$ is also naturally a highest weight
$\hat{\g}$-module of level $\ell$ with $1_{W}$ as a highest weight vector
where for $a\in \g,\; n\in \Z$, $a\otimes t^{n}$ acts as $a(x)_{n}$.
If $W$ is taken to be what is called the Weyl vacuum module 
$V_{\hat{\g}}(\ell,0)$ of level $\ell$, 
$\<S_{V_{\hat{\g}}(\ell,0)}\>$ is isomorphic to $V_{\hat{\g}}(\ell,0)$
 as a $\hat{\g}$-module. Consequently,
$V_{\hat{\g}}(\ell,0)$ has a natural vertex algebra structure and 
any highest weight $\hat{\g}$-module of level $\ell$ is naturally a 
$V_{\hat{\g}}(\ell,0)$-module. This describes the connection between 
(untwisted) affine Lie algebras and vertex algebras.

Now, consider a quantum affine algebra $U_{q}(\hat{\g})$ (see
\cite{drinfeld}, \cite{fj}, cf.  \cite{jimbo}). One can show
that for any generating functions $a(x),b(x)$ in the Drinfeld's
realization (see \cite{drinfeld}) we have
\begin{eqnarray*}
f(x_{1}/x_{2})a(x_{1})b(x_{2})=g(x_{1}/x_{2})b(x_{2})a(x_{1})
\end{eqnarray*}
for some nonzero polynomials $f(x), g(x)$.
With this in mind, for any vector space $W$, 
we study what we call pseudo-local subsets of $\E(W)$, where
$U$ is said to be pseudo-local if for any $a(x),b(x)\in U$, there exist
$u_{i}(x),v_{i}(x)\in U$ and $0\ne p(x_{1},x_{2})\in
\C[x_{1}^{\pm 1},x_{2}^{\pm 1}],\; 
\psi_{i}(x_{1},x_{2})\in \C((x_{2}))((x_{1}))$ for
$i=1,\dots,r$ such that
\begin{eqnarray*}
p(x_{1},x_{2})a(x_{1})b(x_{2})
=\sum_{i=1}^{r}\psi_{i}(x_{1},x_{2})u_{i}(x_{2})v_{i}(x_{1}).
\end{eqnarray*}
In fact, in this paper we study more general subsets of $\E(W)$, which
we call quasi compatible subsets. A subset $U$ of $\E(W)$ is quasi compatible if
for any $a_{1}(x),\dots, a_{n}(x)\in U$, there exists a nonzero polynomial
$f(x,y)$ such that
$$\left(\prod_{1\le i<j\le n}f(x_{i},x_{j})\right)a_{1}(x_{1})\cdots a_{n}(x_{n})
\in \Hom (W,W((x_{1},\dots,x_{n}))).$$
(This notion and the notion of compatibility studied in \cite{li-g1}
were motivated by a notion of compatibility in \cite{b-Gva}.) 

Having chosen our study objects we now consider possible operations.
First,  note that the (old) operations defined in (\ref{e-old-def}) 
on the {\em whole space} $\E(W)$ are no longer appropriate to study
general quasi compatible subsets or pseudo local subsets.
For example, in the study of vertex superalgebras in \cite{li-local},
if $\alpha(x), \beta(x)\in \E(W)$ satisfy the relation
$$(x_{1}-x_{2})^{k}\a(x_{1})\b(x_{2})=-(x_{1}-x_{2})^{k}\b(x_{2})\a(x_{1}),$$
for $n\in \Z$, the definition of $\a(x)_{n}\b(x)$ should be modified
with the minus sign between the two main terms
on the right hand side being changed to plus sign.
In \cite{li-gamma}, new operations were introduced
where the new operations only use the product $\a(x_{1})\b(x_{2})$,
but not the product $\b(x_{2})\a(x_{1})$.
Also, they are ``partial'' operations defined for each ordered pair $(\a(x),\b(x))$
with a certain property which is satisfied by any pair from a 
quasi compatible subset of $\E(W)$.
Specifically, let $\a(x),\b(x)\in \E(W)$ be such that 
$$p(x_{1},x_{2})\a(x_{1})\b(x_{2})\in \Hom (W,W((x_{1},x_{2})))$$
for some nonzero polynomial $p(x_{1},x_{2})$.
The products $\a(x)_{n}\b(x)$ were defined in terms of the generating function
$Y_{\E}(\a(x),x_{0})\b(x)=\sum_{n\in \Z}\a(x)_{n}\b(x) x_{0}^{-n-1}$ by 
\begin{eqnarray}\label{edef-YE-int}
Y_{\E}(\a(x),x_{0})\b(x)
=\iota_{x,x_{0}} (1/p(x+x_{0},x))
\Res_{x_{1}}x_{1}^{-1}\delta\left(\frac{x+x_{0}}{x_{1}}\right)
\left(p(x_{1},x)\a(x_{1})\b(x)\right),
\end{eqnarray}  
where $\iota_{x,x_{0}}$ is the natural embedding of the field $\C(x,x_{0})$ 
of rational functions into $\C((x))((x_{0}))$.
Roughly speaking, $\a(x)_{n}\b(x)$ were defined by the operator product expansion
of $\a(x_{1})\b(x_{2})$.
These are the operations we use in this paper. 
(Note that the new definition agrees with the old one
when $\a(x), \b(x)$ are mutually local.)

We then prove that for any quasi compatible subset $U$
of $\E(W)$, there exists a (unique) smallest closed
quasi compatible subspace $\<U\>$ 
containing $U$ and $1_{W}$, such that $(\<U\>,Y_{\E},1_{W})$ is a nonlocal vertex
algebra,  and furthermore, $W$ is naturally a quasi module
(defined in Section 2) for the nonlocal vertex algebra with 
$Y_{W}(a(x),x_{0})=a(x_{0})$. 
We also show that if $U$ is pseudo-local, $\<U\>$ is also pseudo-local. Thus,
in addition that $W$ is a faithful quasi $\<U\>$-module,
$\{ Y_{W}(\alpha(x),x_{0})\;|\; \alpha(x)\in \<U\>\}$ is pseudo-local.
This naturally leads us to a notion of pseudo-local quasi module 
for a general nonlocal vertex algebra. In terms of this notion, we have
that for any vector space $W$ and for any pseudo-local subset $U$
of $\E(W)$, $W$ is a faithful pseudo-local quasi $\<U\>$-module. 
In particular, taking $W$ to be a highest weight $U_{q}(\hat{\g})$-module,
we have a nonlocal vertex algebra $V_{W}$ generated
by the generating functions of $U_{q}(\hat{\g})$ with $W$ 
as a faithful pseudo-local quasi module.

Furthermore, we prove that if $V$ is a general nonlocal vertex algebra 
and if $(W,Y_{W})$ is a pseudo-local quasi $V$-module 
satisfying a certain condition, then the pseudo-locality 
of the subspace $\{Y_{W}(v,x)\;|\; v\in V\}$ of $\E(W)$ 
uniquely defines a quantum Yang-Baxter operator 
$\S(x_{1},x_{2})$ on $V$ by 
$$\S(x,y)(v\otimes u)
=\sum_{i=1}^{r}v^{(i)}\otimes u^{(i)}\otimes f_{i}(x,y)
\in V\otimes V\otimes \C_{*}(x,y)$$
for $u,v\in V$, if
$$p(x_{1},x_{2})Y_{W}(u,x_{1})Y_{W}(v,x_{2})
=\sum_{i=1}^{r}p(x_{1},x_{2})\iota_{x_{2},x_{1}}(f_{i}(x_{2},x_{1}))
Y_{W}(v^{(i)},x_{2})Y_{W}(u^{(i)},x_{1})$$
for some $0\ne p(x,y)\in \C[x,y]$, 
where $\C_{*}(x,y)=\{ g/h\;|\; g\in \C[[x,y]],\; 0\ne h\in \C[x,y]\}$.
By definition a quantum Yang-Baxter operator 
$\S(x_{1},x_{2})$ satisfies the quantum Yang-Baxter equation
$$\S_{12}(x_{1},x_{2})\S_{13}(x_{1},x_{3})\S_{23}(x_{2},x_{3})
=\S_{23}(x_{2},x_{3})\S_{13}(x_{1},x_{3})\S_{12}(x_{1},x_{2}).$$
This naturally leads us to a notion of
$\S$-quasi $V$-module with respect to
a quantum Yang-Baxter operator $\S(x_{1},x_{2})$ on $V$.
(Recall that twisted modules for a vertex algebra $V$ are associated 
with finite order automorphisms of $V$.)

What we expect is that if $W$ is the basic highest weight module
constructed in \cite{fj} for the quantum affine algebra
$U_{q}(\hat{\g})$, $W$ is an $\S$-quasi module for the associated
nonlocal vertex algebra $V_{W}$ with respect to a (uniquely
determined) quantum Yang-Baxter operator $\S(x_{1},x_{2})$.  To
establish this we need to find an explicit characterization of on the
nonlocal vertex algebra $V_{W}$.  Recall that in the case with
(untwisted) affine Lie algebras $\hat{\g}$, the associated vertex
algebras are naturally highest weight $\hat{\g}$-modules so that they
can be easily determined.  In the case with $U_{q}(\hat{\g})$, on the
associated nonlocal vertex algebra $V_{W}$ the natural ``adjoint''
action does not give rise to a $U_{q}(\hat{\g})$-module structure and
the structure of $V_{W}$ is still to be determined.  To a certain
extent, this reflects the view point of \cite{efr} that the space of
states and the space of fields should be different in the quantum
case.  On the other hand, it is philosophically important to note that
this seemingly new phenomenon is rather similar to what happened in
the study of twisted modules and quasi modules for an ordinary vertex
algebra as we mention next.

For a general vertex (operator) algebra $V$, in addition to the notion
of $V$-module, we have the notion of $\sigma$-twisted $V$-module,
where $\sigma$ is an automorphism of $V$ of finite order.  Let
$\g$ be a finite-dimensional simple Lie algebra and let $\sigma$ be 
an automorphism of $\g$ of finite order $r$.  We have the
twisted affine Lie algebra $\hat{\g}[\sigma]=\prod_{i=0}^{r-1}
\g_{i}\otimes t^{-\frac{i}{r}}\C[t,t^{-1}]\oplus \C c$ (see
\cite{kacbook1}, \cite{flm}).  It was proved in \cite{li-twisted}
that for any highest weight $\hat{\g}[\sigma]$-module $W$ of level
$\ell$, the canonical generating functions of $\hat{\g}[\sigma]$ on
$W$, as elements of $\Hom (W,W((x^{1/r})))$, generate a vertex
algebra $V_{W}$ with $W$ as a $\sigma$-twisted module and
that $V_{W}$ is naturally a $\hat{\g}$-module of level $\ell$.  (There
does not exist a notion of twisted vertex algebra.)  That is, even
though $W$ is a module for the {\em twisted} affine Lie algebra
$\hat{\g}[\sigma]$, the generated vertex algebra $V_{W}$ is a module
for the {\em untwisted} affine Lie algebra $\hat{\g}$.

In \cite{li-gamma}, to associate vertex algebras to certain
infinite-dimensional Lie algebras we studied ``quasi local'' subsets
of $\E(W)$ for any given vector space $W$ and we proved that any quasi
local subset generates in a certain natural way a vertex algebra with
$W$ as a ``quasi module.''  As an application of this general result
we successfully associated vertex algebras in terms of quasi modules
to certain quantum torus Lie algebras.  For a suitably defined
restricted module $W$ for the quantum torus Lie algebras, the
canonical generating functions generate a vertex algebra $V_{W}$ with
$W$ as a quasi module, while the vertex algebra $V_{W}$ is a module
for a {\em new} Lie algebra, rather than a module for the quantum
torus Lie algebras.

These phenomena indicate that for a highest weight
$U_{q}(\hat{\g})$-module $W$, the nonlocal vertex algebra $V_{W}$
generated by the generating functions should be a module for some {\em
new} algebra.  It would be very interesting if this ``new'' algebra
turns out to be a new quantum algebra.  This will be investigated in a
sequel. Then we shall have a complete solution to the problem mentioned 
in the first paragraph.

While working on the structure of the nonlocal vertex algebra $\<U\>$ 
for a pseudo-local subset of $\E(W)$, we find that if $U$ is of a very special type, 
$\<U\>$ very much resembles
a quantum vertex operator algebra in the sense of \cite{ek}.
This leads us to notions of $\S$-locality and weak quantum 
vertex algebras.
A subset $U$ of $\E(W)$ is said to be $\S$-local if for any
$a(x),b(x)\in U$, there exist $k\in \N$,
$f_{1}(x),\dots,f_{r}(x)\in \C((x))$ and 
$u_{1}(x),v_{1}(x),\dots,u_{r}(x),v_{r}(x)\in U$ such that
$$(x_{1}-x_{2})^{k}a(x_{1})b(x_{2})=\sum_{i=1}^{r}(x_{1}-x_{2})^{k}
\iota_{x_{2},x_{1}}(f_{i}(x_{2}-x_{1}))u_{i}(x_{2})v_{i}(x_{1}).$$
A weak quantum vertex algebra is defined to be a nonlocal vertex algebra $V$
such that $\{ Y(v,x)\;|\; v\in V\}$ is an $\S$-local subspace of
$\E(V)$. We prove that for any $\S$-local subset $U$ of $\E(W)$,
$\<U\>$ is a weak quantum vertex algebra with $W$ as
a module instead of a quasi module.
Furthermore, motivated by the notion of quantum vertex operator
algebra in \cite{ek} we define
a quantum vertex algebra to be a weak 
quantum vertex algebra $V$ equipped with a
rational quantum Yang-Baxter
$\S: V\otimes V\rightarrow V\otimes V\otimes \C(x)$, which by
definition satisfies the rational quantum Yang-Baxter equation
$$\S_{12}(x)\S_{13}(x+z)\S_{23}(z)=\S_{23}(z)\S_{13}(x+z)\S_{12}(x),$$
satisfying the condition that
for $a,b\in V$ there exists $k\in \N$ such that 
$$(x_{1}-x_{2})^{k}Y(a,x_{1})Y(b,x_{2})=
\sum_{i=1}^{r}(x_{1}-x_{2})^{k}\iota_{x_{2},x_{1}}(\psi_{i}(x_{2}-x_{1}))
Y(u_{i},x_{2})Y(v_{i},x_{1}),$$
where $\S(x)(b\otimes a)=\sum u_{i}\otimes v_{i}\otimes \psi_{i}(x)$.
We prove that if a weak quantum vertex algebra
$V$ is nondegenerate in the sense of \cite{ek},
then the $\S$-locality uniquely
defines a unitary rational quantum Yang-Baxter operator $\S$ such that
$V$ equipped with $\S$ is a quantum vertex algebra.

The main results of this paper were reported in the conference
``Perspectives arising from vertex algebras,'' Osaka, Japan, November
7-12, 2004. We would like to thank Professors Nagatomo and Yamada for
their hospitality.

This paper is organized as follows: In Section 2, we study quasi
compatible subsets of $\E(W)$ and prove that any maximal quasi
compatible subspace of $\E(W)$ is naturally a nonlocal vertex algebra
with $W$ as a quasi module.  In Section 3, we study pseudo-local
subsets and $\S$-local subsets of $\E(W)$.  In Section 4, we study
pseudo-local quasi modules and $\S$-quasi modules for nonlocal vertex
algebras.  In Section 5, we study (weak) quantum vertex algebras.
In Section 6, we give some general construction results and examples.

\section{Constructing nonlocal vertex algebras from 
quasi compatible subsets of $\E(W)$}
In this section after recalling the notion of nonlocal vertex algebra
and defining a notion of quasi module,
we study what we call quasi compatible subsets of $\Hom (W,W((x)))$
for any vector space $W$. We prove that any maximal quasi compatible
subspace is naturally a nonlocal vertex algebra
with $W$ as a faithful quasi module and that any quasi compatible subset
generates a nonlocal vertex algebra with $W$ as a faithful quasi module.

First, throughout this paper,
$x,y, z,x_{0},x_{1},x_{2}, x_{3},\dots$ are mutually commuting independent
formal variables. Vector spaces are considered to be over the field $\C$ of
complex numbers. We use all the standard formal variable notations and
conventions. 
For example, $U[[x_{1},x_{2},\dots,x_{r}]]$ is the space of 
(infinite) nonnegative integral power series in $x_{1},\dots,x_{r}$ with coefficients in $U$
and $U((x_{1},x_{2},\dots,x_{r}))$ 
is the space of lower truncated (infinite) integral
power series in $x_{1},\dots,x_{r}$ with coefficients in $U$, i.e.,
\begin{eqnarray}
U((x_{1},\dots,x_{r}))=U[[x_{1},\dots,x_{r}]][x_{1}^{-1},\dots,x_{r}^{-1}].
\end{eqnarray}
Clearly,
\begin{eqnarray}
\C[x_{1},\dots,x_{r}]\subset \C[[x_{1},\dots,x_{r}]]
\subset \C((x_{1},\dots,x_{r}))
\subset \C((x_{1}))((x_{2}))\cdots ((x_{r})).
\end{eqnarray}
Note that $\C((x_{1}))$, $\C((x_{1}))((x_{2}))$,\dots,
$\C((x_{1}))((x_{2}))\cdots ((x_{r}))$ are fields.

Denote by $\C_{*}(x_{1},\dots,x_{r})$ the localization of 
$\C[[x_{1},\dots,x_{r}]]$ at
$\C[x_{1},\dots,x_{r}]^{\times}$ (the nonzero polynomials):
\begin{eqnarray}
\C_{*}(x_{1},\dots,x_{r})
=\left\{ \frac{f}{p}\;|\; f\in \C[[x_{1},\dots,x_{r}]],\; 
p\in \C[x_{1},\dots,x_{r}]^{\times}\right\}.
\end{eqnarray}
The field $\C(x_{1},\dots,x_{r})$ of rational functions, 
i.e., the fraction field of the polynomial ring $\C[x_{1},\dots,x_{r}]$,
is a subring of $\C_{*}(x_{1},\dots,x_{r})$. 
Since $\C[[x_{1},\dots,x_{r}]]$ is a subring of 
$\C((x_{1}))((x_{2}))\cdots ((x_{r}))$, there exists a unique
embedding of the ring $\C_{*}(x_{1},\dots,x_{r})$ into 
the field $\C((x_{1}))((x_{2}))\cdots ((x_{r}))$.
Denote this embedding by $\iota_{x_{1},\dots,x_{r}}$: 
\begin{eqnarray}
\iota_{x_{1},\dots,x_{r}}:
\C_{*}(x_{1},\dots,x_{r})\rightarrow \C((x_{1}))((x_{2}))\cdots
((x_{r})).
\end{eqnarray}
Similarly, for any permutation $\sigma$ on $\{1,\dots,r\}$,
there exists a unique embedding
\begin{eqnarray}
\iota_{x_{\sigma(1)},\dots,x_{\sigma(r)}}:
\C_{*}(x_{1},\dots,x_{r})\rightarrow \C((x_{\sigma(1)}))\cdots
((x_{\sigma(r)})).
\end{eqnarray}
Notice that partial differential operators $\frac{\partial}{\partial x_{i}}$ act on
both the domains and codomains and they
commute with the iota-maps. Furthermore, all the iota-maps
are $\C((x_{1},\dots,x_{r}))$-homomorphisms and they are identity
on $\C((x_{1},\dots,x_{r}))$.

These generalized iota maps, defined in \cite{li-gamma}, generalize
the usual iota maps introduced in \cite{flm}, \cite{fhl} and \cite{ll}. 
For example, for any nonzero complex number $\alpha$ and for any integer $n$, we have
\begin{eqnarray}
& &\iota_{x_{1},x_{2}}(x_{1}-\alpha x_{2})^{n}
=\sum_{i\ge 0}\binom{n}{i}(-\alpha)^{i}x_{1}^{n-i}x_{2}^{i}=(x_{1}-\alpha x_{2})^{n},\\
& &\iota_{x_{2},x_{1}}(x_{1}-\alpha x_{2})^{n}
=\sum_{i\ge 0}\binom{n}{i}(-\alpha)^{n-i}x_{2}^{n-i}x_{1}^{i}=(-\alpha x_{2}+x_{1})^{n},
\end{eqnarray}
where we are using the usual binomial series expansion conventions.

In formal calculus (\cite{flm}, \cite{fhl}, \cite{ll}), associativity and 
cancelation for products of formal series are subtle issues. 
Let $A$ be an associative algebra with identity 
(over $\C$) and let $U$ be an $A$-module. Let
$$F,G\in A[[x_{1}^{\pm 1},\dots,x_{n}^{\pm 1}]],\;\;
H\in U[[x_{1}^{\pm 1},\dots,x_{n}^{\pm 1}]].$$
The associativity $F(GH)=(FG)H$ does not hold in general 
and on the other hand, it does hold
if all the products $GH$, $FG$ and $FGH$ exist. 
As in associative algebra theory,
associativity sometimes leads to cancelation.
The following are some (simple) useful facts: 

\bl{lcancellation}
Let $U$ be any vector space. 
(a) Let $H(x_{1},x_{2}),K(x_{1},x_{2})\in
U((x_{1}))((x_{2}))$. If 
\begin{eqnarray}
f(x_{1},x_{2})H(x_{1},x_{2})=f(x_{1},x_{2})K(x_{1},x_{2})
\end{eqnarray}
for some $0\ne f(x_{1},x_{2})\in \C((x_{1}))((x_{2}))$,
then $H(x_{1},x_{2})=K(x_{1},x_{2})$.

(b) Let
$H_{1}\in U((x))((x_{1}))((x_{2}))$ and $H_{2}\in U((x))((x_{2}))((x_{1}))$.
If there exists a (nonzero) polynomial $g(x,x_{1},x_{2})$ with
$g(x,0,0)\ne 0$ such that
$$g(x,x_{1},x_{2})H_{1}=g(x,x_{1},x_{2})H_{2},$$
then $H_{1}=H_{2}$.
\el

\begin{proof} Part (a) was proved in \cite{li-gamma}. For (b), since $g(x,0,0)\ne 0$, we have
$$\iota_{x,x_{1},x_{2}}(1/g(x,x_{1},x_{2}))=\iota_{x,x_{2},x_{1}}(1/g(x,x_{1},x_{2}))
\in \C((x))[[x_{1},x_{2}]].$$
Multiplying both sides of the given identity 
by $\iota_{x,x_{1},x_{2}}(1/g(x,x_{1},x_{2}))$ we get $H_{1}=H_{2}$.
\end{proof}

We shall often use the following fact
which was proved in \cite{li-g1} (cf. \cite{fhl}, \cite{ll}):

\bl{lformaljacobiidentity}
Let $U$ be a vector space and let
\begin{eqnarray}
A(x_{1},x_{2})&\in& U((x_{1}))((x_{2})),\\
B(x_{1},x_{2})&\in& U((x_{2}))((x_{1})),\\
C(x_{0},x_{2})&\in& U((x_{2}))((x_{0})).
\end{eqnarray}
Then 
\begin{eqnarray}\label{eformaljacobiABC}
& &x_{0}^{-1}\delta\left(\frac{x_{1}-x_{2}}{x_{0}}\right)A(x_{1},x_{2})
-x_{0}^{-1}\delta\left(\frac{x_{2}-x_{1}}{-x_{0}}\right)B(x_{1},x_{2})
\nonumber\\
&=&x_{1}^{-1}\delta\left(\frac{x_{2}+x_{0}}{x_{1}}\right)C(x_{0},x_{2})
\end{eqnarray}
if and only if there exist nonnegative integers $k$ and $l$ such that
\begin{eqnarray}
(x_{1}-x_{2})^{k}A(x_{1},x_{2})
&=&(x_{1}-x_{2})^{k}B(x_{1},x_{2}),\label{eformalA=B}\\
(x_{0}+x_{2})^{l}A(x_{0}+x_{2},x_{2})
&=&(x_{0}+x_{2})^{l}C(x_{0},x_{2}).\label{eformalA=C}
\end{eqnarray}
\el

We now recall from \cite{li-g1} the notion of
nonlocal vertex algebra (cf. \cite{bk}):

\bd{dgva}
{\em A {\em nonlocal vertex algebra} is a vector space $V$ equipped
with a linear map 
\begin{eqnarray}
Y(\cdot,x): V &\rightarrow & \Hom (V,V((x)))\subset (\End V)[[x,x^{-1}]]\nonumber\\
v&\mapsto& Y(v,x)=\sum_{n\in \Z}v_{n}x^{-n-1}\;\;\; (v_{n}\in \End V),
\end{eqnarray}
or equivalently, a linear map from $V\otimes V$ to $V((x))$,
and equipped with a distinguished vector ${\bf 1}$ such that
for $v\in V$
\begin{eqnarray}
& &Y({\bf 1},x)v=v,\\
& &Y(v,x){\bf 1}\in V[[x]]\;\;\mbox{ and }\;\; \lim_{x\rightarrow
0}Y(v,x){\bf 1}=v
\end{eqnarray}
and such that for $u,v,w\in V$, there exists a nonnegative integer
$l$ such that
\begin{eqnarray}
(x_{0}+x_{2})^{l}Y(u,x_{0}+x_{2})Y(v,x_{2})w=
(x_{0}+x_{2})^{l}Y(Y(u,x_{0})v,x_{2})w.
\end{eqnarray}}
\ed

\br{rnonlocalfield-alg}
{\em Note that the notion of nonlocal vertex algebra is exactly
the same as the notion of weak axiomatic $G_{1}$-vertex algebra
studied in \cite{li-g1}. On the other hand, 
it follows from Lemma \ref{lcon1} that 
the notion of nonlocal vertex algebra is equivalent to the notion of field algebra
studied in \cite{bk} (cf. \cite{kacbook2}).}
\er

The following was proved in \cite{li-g1}:
 
\bl{lcon1}
Let $V$ be a nonlocal vertex algebra. Define ${\cal{D}}\in \End_{\C}V$ by
\begin{eqnarray}
{\cal{D}}(v)=\left({d\over dx}Y(v,x){\bf 1}\right)|_{x=0}
\left(=v_{-2}{\bf 1}\right)
\;\;\;\mbox{ for }v\in V.
\end{eqnarray}
Then
\begin{eqnarray}\label{edproperty}
[{\cal{D}},Y(v,x)]=Y({\cal{D}}(v),x)={d\over dx}Y(v,x)
\;\;\;\mbox{ for }v\in V.
\end{eqnarray}
Furthermore, for $v\in V$,
\begin{eqnarray}
& &e^{x{\cal{D}}}Y(v,x_{1})e^{-x{\cal{D}}}=Y(e^{xD}v,x_{1})=Y(v,x_{1}+x),
\label{econjugationformula1}\\
& &Y(v,x){\bf 1}=e^{x{\cal{D}}}v.\label{ecreationwithd}
\end{eqnarray}
\el

Let $S$ be a subset of a nonlocal vertex algebra $V$.
Denote by $\<S\>$ the {\em subalgebra of $V$ generated by $S$}, 
which by definition is the smallest subalgebra of
$V$, containing $S$. From \cite{li-g1} we have:

\bl{lsubalgebragenerated}
For any subset $S$ of $V$, the subalgebra $\<S\>$ generated by $S$ 
is linearly spanned by the vectors
\begin{eqnarray}\label{espannform}
u^{(1)}_{n_{1}}\cdots u^{(r)}_{n_{r}}{\bf 1}
\end{eqnarray}
for $r\ge 0,\; u^{(i)}\in S,\; n_{1},\dots,n_{r}\in \Z$.
\el

\bd{dmodule}
{\em A $V$-{\em module} is a vector space $W$ equipped with a linear map
\begin{eqnarray}
Y_{W}(\cdot,x): V &\rightarrow & \Hom (W,W((x)))\subset (\End W)[[x,x^{-1}]]\nonumber\\
v&\mapsto& Y_{W}(v,x)=\sum_{n\in \Z}v_{n}x^{-n-1}\;\;\; (v_{n}\in \End W)
\end{eqnarray}
such that all the following axioms hold: 
\begin{eqnarray}
Y_{W}({\bf 1},x)=1_{W}\;\;\;
\mbox{(where $1_{W}$ is the identity operator on $W$)},
\end{eqnarray}
and for any $u, v\in V,\; w\in W$, there exists $l\in \N$
such that
\begin{eqnarray}\label{emoduleweakassoc}
(x_{0}+x_{2})^{l}Y_{W}(u,x_{0}+x_{2})Y_{W}(v,x_{2})w
=(x_{0}+x_{2})^{l}Y_{W}(Y(u,x_{0})v,x_{2})w.
\end{eqnarray}}
\ed

We extend the notion of quasi module defined in \cite{li-gamma}
for a vertex algebra as follows:

\bd{dquasi-module}
{\em Let $V$ be a nonlocal vertex algebra. A {\em quasi $V$-module} is a vector space
$W$ equipped with a linear map $Y_{W}$ from $V$ to $\Hom (W,W((x)))$ such that
\begin{eqnarray}
Y_{W}({\bf 1},x)=1_{W}
\end{eqnarray}
and such that for $u,v\in V,\; w\in W$, there
exists a nonzero polynomial $f(x_{1},x_{2})$ such that
\begin{eqnarray}
f(x_{0}+x_{2},x_{2})Y_{W}(u,x_{0}+x_{2})Y_{W}(v,x_{2})w
=f(x_{0}+x_{2},x_{2})Y_{W}(Y(u,x_{0})v,x_{2})w.
\end{eqnarray}}
\ed

Slightly modifying the proof of Lemma \ref{lcon1} in \cite{li-g1} we have:

\bl{lquasi-dproperty}
Let $(W,Y_{W})$ be a quasi module for a nonlocal vertex algebra $V$.
Then 
\begin{eqnarray}\label{eDproperty-quasi-module}
Y_{W}(\D v,x)=\frac{d}{dx}Y_{W}(v,x)\;\;\;\mbox{ for }v\in V.
\end{eqnarray}
\el

{}From now on (for the rest of this section), 
let $W$ be a fixed general vector space over $\C$.
Set
\begin{eqnarray}
\E(W)=\Hom (W,W((x))).
\end{eqnarray}
The formal differential operator is a linear operator on $\E(W)$,
which we denote by $D$:
\begin{eqnarray}
D={d\over dx}.
\end{eqnarray}
Clearly, $\C((x))$ naturally acts on $\E(W)$ and so does
the Lie algebra $\C((x)D$.

We consider the following action of the multiplicative group $\C^{\times}$ of nonzero complex numbers 
on $\E(W)$:
\begin{eqnarray}
R_{\alpha}a(x)=a(\alpha x)\;\;\;\mbox{ for }\alpha\in \C^{\times},\;
a(x)\in \E(W).
\end{eqnarray}

The following notion generalizes the notions of compatibility in 
\cite{li-g1} and \cite{li-gamma}:

\bd{dcompatibility}
{\em An (ordered) sequence $(\psi^{(1)},\dots,\psi^{(r)})$ 
in $\E(W)$ is said to be {\em quasi compatible} 
if there exists a  nonzero polynomials $p(x,y)$ such that
\begin{eqnarray}
\left(\prod_{1\le i<j\le r}p(x_{i},x_{j})\right)
\psi^{(1)}(x_{1})\cdots \psi^{(r)}(x_{r})
\in \Hom (W,W((x_{1},\dots,x_{r}))).
\end{eqnarray}
A set or a subspace $S$ of $\E(W)$ is said to be 
{\em quasi compatible} if any finite sequence in $S$ 
is quasi compatible.}
\ed

The following partial operations on $\E(W)$ were defined in \cite{li-gamma}:

\bd{d-operation}
{\em Let $(a(x),b(x))$ be a quasi compatible (ordered) pair
in $\E(W)$. For $\alpha\in \C^{\times},\; n\in \Z$, we define
$a(x)_{(\alpha,n)}b(x)\in (\End W)[[x,x^{-1}]]$ in terms of generating
function
\begin{eqnarray}
Y_{\E}^{(\alpha)}(a(x),x_{0})b(x)=\sum_{n\in \Z}(a(x)_{(\alpha,n)}b(x)) x_{0}^{-n-1}
\in (\End W)[[x_{0}^{\pm 1},x^{\pm 1}]]
\end{eqnarray}
by
\begin{eqnarray}
& &Y_{\E}^{(\alpha)}(a(x),x_{0})b(x)\nonumber\\
&=&\Res_{x_{1}}\iota_{x,x_{0}}\left(p(x_{0}+\alpha x,x)^{-1}\right)
x_{1}^{-1}\delta\left(\frac{\alpha x+x_{0}}{x_{1}}\right)
\left(p(x_{1},x)a(x_{1})b(x)\right)\\
&=&\iota_{x,x_{0}}\left(p(x_{0}+\alpha x,x)^{-1}\right)
\left(p(x_{1},x)a(x_{1})b(x)\right)|_{x_{1}=\alpha x+x_{0}},\label{edef-3.15}
\end{eqnarray}
where $p(x_{1},x_{2})$ is any nonzero polynomial such that
$$p(x_{1},x_{2})a(x_{1})b(x_{2})\in \Hom (W,W((x_{1},x_{2}))).$$
We particularly set
\begin{eqnarray}
& &Y_{\E}(a(x),x_{0})b(x)=Y_{\E}^{(1)}(a(x),x_{0})b(x),\\
& &a(x)_{n}b(x)=a(x)_{(1,n)}b(x)\;\;\;\mbox{  for }n\in \Z.
\end{eqnarray}}
\ed

It was proved in \cite{li-gamma} that $Y_{\E}^{(\alpha)}(a(x),x_{0})b(x)$ 
is well defined,
that is, the expression exists and is independent of the choice of
the polynomial $p(x_{1},x_{2})$. Furthermore, the following 
Propositions \ref{pexistence}-\ref{pbasic-property2} were also  proved:

\bp{pexistence} 
Let $(a(x),b(x))$ be a quasi compatible (ordered) pair
in $\E(W)$. Then
\begin{eqnarray}\label{eew-closed}
a(x)_{(\alpha,n)}b(x)\in \E (W)\;(=\Hom (W,W((x))))
\;\;\;\mbox{ for }\alpha\in \C^{\times},\; n\in \Z.
\end{eqnarray}
Let $p(x_{1},x_{2})$ be any nonzero polynomial such that
\begin{eqnarray}\label{efabcondition}
p(x_{1},x_{2})a(x_{1})b(x_{2})\in \Hom(W,W((x_{1},x_{2})))
\end{eqnarray}
and let $k$ be an integer such that
$$x_{0}^{k}\iota_{x_{2},x_{0}}
\left(p(x_{0}+\alpha x_{2},x_{2})^{-1}\right)\in \C((x_{2}))[[x_{0}]].$$
Then
\begin{eqnarray}\label{e-truncation-ab}
a(x)_{(\alpha,n)}b(x)=0\;\;\;\mbox{ for }n\ge k.
\end{eqnarray}
Furthermore, for $w\in W$, let $l$ be a nonnegative integer depending on $w$ such that
$$x_{1}^{l}p(x_{1},x_{2})a(x_{1})b(x_{2})w\in W[[x_{1},x_{2}]][x_{2}^{-1}].$$
Then
\begin{eqnarray}\label{eweakassoc-def}
& &(x_{0}+\alpha x_{2})^{l}p(x_{0}+\alpha x_{2},x_{2})
(Y_{\E}^{(\alpha)}(a(x_{2}),x_{0})b(x_{2}))w\nonumber\\
&=&(x_{0}+\alpha x_{2})^{l}p(x_{0}+\alpha x_{2},x_{2})
a(x_{0}+\alpha x_{2})b(x_{2})w.
\end{eqnarray}
\ep

\bp{pvacuum}
For any $a(x)\in \E(W)$, the sequences $(1_{W}, a(x))$ and
$(a(x), 1_{W})$ are quasi compatible and for $\alpha\in \C^{\times}$ we have
\begin{eqnarray}
& &Y_{\E}^{(\alpha)}(1_{W},x_{0})a(x)=a(x),\label{e-vacuum-calculus}\\
& &Y_{\E}^{(\alpha)}(a(x),x_{0})1_{W}=a(\alpha x+x_{0})
=e^{\alpha^{-1}x_{0}{d\over dx}}a(\alpha x)
=e^{\alpha^{-1}x_{0}{d\over dx}}R_{\alpha}a(x).\label{e-creation-calculus}
\end{eqnarray}
In particular,
\begin{eqnarray}
& &a(x)_{(\alpha,-1)}1_{W}=R_{\alpha}a(x),\label{ea-alpha-1-R}\\
& &a(x)_{-2}1_{W}\;\left(=a(x)_{(1,-2)}1_{W}\right)
=\frac{d}{dx}a(x)=Da(x).\label{ea-2=da}
\end{eqnarray}
\ep

\bp{pbasic-property2}
Let $(a(x),b(x))$ be a quasi compatible (ordered) pair in $\E(W)$.
Then the ordered pairs $(a'(x),b(x)), (a(x),b'(x))$ and $(a(\alpha x), b(\beta x))$
for $\alpha,\beta\in \C^{\times}$ are quasi compatible.
Furthermore, we have
\begin{eqnarray}
& &Y_{\E}^{(\alpha)}(Da(x),x_{0})b(x)={\partial\over \partial x_{0}}
Y_{\E}^{(\alpha)}(a(x),x_{0})b(x),\label{e-Dderivative-calculus}\\
& & DY_{\E}^{(\alpha)}(a(x),x_{0})b(x)-Y_{\E}^{(\alpha)}(a(x),x_{0})D b(x)
=\alpha {\partial\over \partial x_{0}}
Y_{\E}^{(\alpha)}(a(x),x_{0})b(x),\label{e-Dbracket-calculus}\\
& &R_{\alpha}Y_{\E}^{(\beta)}(a(x),x_{0})b(x)
=Y_{\E}^{(\alpha\beta)}(a(x),x_{0})R_{\alpha}b(x),
\label{e-conjugation-calculus}\\
& &Y_{\E}^{(\alpha)}(R_{\beta}a(x),\beta^{-1}x_{0})b(x)
=Y_{\E}^{(\alpha\beta)}(a(x),x_{0})b(x),
\label{e-conjugation-calculus-2}
\end{eqnarray}
In particular,
\begin{eqnarray}\label{e-conjugation-calculus-3}
Y_{\E}^{(\alpha)}(a(x),x_{0})b(x)=R_{\alpha}Y_{\E}(a(x),x_{0})R_{\alpha^{-1}}b(x)
=Y_{\E}(R_{\alpha}a(x),\alpha^{-1}x_{0}).
\end{eqnarray}
\ep

We shall need the following technical result:

\bl{lproof-need}
Let $(a_{i}(x),b_{i}(x))$ $(i=1,\dots,n)$ be quasi compatible
ordered pairs in $\E(W)$. Suppose that 
\begin{eqnarray}\label{esum-lemma}
\sum_{i=1}^{n}g_{i}(z,x)a_{i}(z)b_{i}(x)\in \Hom (W,W((z,x)))
\end{eqnarray}
for some $g_{1}(z,x),\dots,g_{n}(z,x)\in \C[z,x]$. Then
\begin{eqnarray}
\sum_{i=1}^{n}g_{i}(x+x_{0},x)Y_{\E}(a_{i}(x),x_{0})b_{i}(x)
=\left(\sum_{i=1}^{n}g_{i}(z,x)a_{i}(z)b_{i}(x)\right)|_{z=x+x_{0}}.
\end{eqnarray}
\el

\begin{proof} Let $0\ne g(x,y)\in \C[x,y]$ be such that
$$g(z,x)a_{i}(z)b_{i}(x)\in \Hom (W,W((z,x)))\;\;\;\mbox{ for }i=1,\dots,n.$$
{}From Definition \ref{d-operation}, we have
$$g(x+x_{0},x)Y_{\E}(a_{i}(x),x_{0})b_{i}(x)
=\left(g(z,x)a_{i}(z)b_{i}(x)\right)\mid_{z=x+x_{0}}$$
for $i=1,\dots,n$. Then using (\ref{esum-lemma}) we have
\begin{eqnarray}
& &g(x+x_{0},x)\sum_{i=1}^{n}g_{i}(x+x_{0},x)Y_{\E}(a_{i}(x),x_{0})b_{i}(x)
\nonumber\\
&=&\sum_{i=1}^{n}g_{i}(x+x_{0},x)\left(g(z,x)a_{i}(z)b_{i}(x)\right)|_{z=x+x_{0}}
\nonumber\\
&=&\left(g(z,x)\sum_{i=1}^{n}g_{i}(z,x)a_{i}(z)b_{i}(x)\right)|_{z=x+x_{0}}
\nonumber\\
&=&g(x+x_{0},x)\left(\sum_{i=1}^{n}g_{i}(z,x)a_{i}(z)b_{i}(x)\right)|_{z=x+x_{0}}.
\end{eqnarray}
By (\ref{e-truncation-ab}) we have $\sum_{i=1}^{n}g_{i}(x+x_{0},x)
Y_{\E}(a_{i}(x),x_{0})b_{i}(x)\in (\Hom (W,W((x)))((x_{0}))$.
Now it follows immediately from the cancelation rule (Lemma \ref{lcancellation}).
\end{proof}

Let $U$ be a subspace of $\E(W)$ such that every ordered pair in $U$
is quasi compatible. We say $U$ is {\em closed} if
\begin{eqnarray}
a(x)_{n}b(x)\in U\;\;\;\mbox{ for }a(x),b(x)\in U,\; n\in \Z.
\end{eqnarray}
We are going to prove that any closed quasi compatible subspace containing $1_{W}$ of $\E(W)$ 
is a nonlocal vertex algebra and that
any maximal quasi compatible subspace of $\E(W)$ is closed and 
is a nonlocal vertex algebra. First we have:

\bl{lclosed}
Assume that $V$ is a subspace of $\E(W)$ such that any sequence
in $V$ of length $2$ or $3$ is quasi compatible and such that
$V$ is closed. Let $\psi(x),\phi(x),\theta(x)\in V$ and 
let $f(x,y)$ be a nonzero polynomial such that
\begin{eqnarray}
& &f(x,y)\phi(x)\theta(y)\in \Hom (W,W((x,y))),
\label{e4.31}\\
& &f(x,y)f(x,z)f(y,z)
\psi(x)\phi(y)\theta(z)\in \Hom (W,W((x,y,z))).\label{e4.32}
\end{eqnarray}
Then
\begin{eqnarray}
& &f(x+x_{1},x)f(x+x_{2},x)f(x+x_{1},x+x_{2})
Y_{\cal{E}}(\psi(x),x_{1})Y_{\cal{E}}(\phi(x),x_{2})\theta(x)\nonumber\\
&=&\left(f(y,x)f(z,x)f(y,z)
\psi(y)\phi (z)\theta(x)\right)|_{y=x+x_{1},z=x+x_{2}}.
\end{eqnarray}
\el

\begin{proof} With (\ref{e4.31}), from Definition \ref{d-operation} we have
\begin{eqnarray}\label{ehphi-theta}
f(x+x_{2},x)Y_{\cal{E}}(\phi (x),x_{2})\theta (x)
=\left(f(z,x)\phi(z)\theta(x)\right)|_{z=x+x_{2}},
\end{eqnarray}
which gives
\begin{eqnarray}
& &f(y,x)f(y,x+x_{2})f(x+x_{2},x)\psi(y)
Y_{\cal{E}}(\phi (x),x_{2})\theta(x)\nonumber\\
&=&\left(f(y,x)f(y,z) f(z,x)
\psi(y)\phi(z)\theta(x)\right)|_{z=x+x_{2}}.
\end{eqnarray}
{}From (\ref{e4.32}) we see that
the expression on the right-hand side lies in $(\Hom (W,W((y, x)))[[x_{2}]]$,
so does the expression on the left-hand side. 
That is,
\begin{eqnarray}
f(y,x)f(y,x+x_{2})f(x+x_{2},x)\psi(y)
Y_{\cal{E}}(\phi(x),x_{2})\theta(x)
\in (\Hom (W,W((y,x)))[[x_{2}]].
\end{eqnarray}
Multiplying by 
$\iota_{x,x_{2}}(f(x+x_{2},x)^{-1})$, which lies in 
$\C ((x))((x_{2}))\subset \C((y,x))((x_{2}))$,
we have
\begin{eqnarray}
f(y,x)f(y,x+x_{2})
\psi(y)Y_{\cal{E}}(\phi(x),x_{2})\theta(x)
\in (\Hom (W,W((y,x)))((x_{2})).
\end{eqnarray}
In view of Lemma \ref{lproof-need},
by considering the coefficient of each power of $x_{2}$, we have
\begin{eqnarray}
& &f(x+x_{1},x)f(x+x_{1},x+x_{2})
Y_{\cal{E}}(\psi(x),x_{1})Y_{\cal{E}}(\phi(x),x_{2})\theta(x)
\nonumber\\
&=&\left(f(y,x)f(y,x+x_{2})
\psi(y)(Y_{\cal{E}}(\phi(x),x_{2})\theta(x)\right)|_{y=x+x_{1}}.
\end{eqnarray}
Using this and (\ref{ehphi-theta}) we have
\begin{eqnarray*}
& &f(x+x_{2},x)f(x+x_{1},x)f(x+x_{1},x+x_{2})
Y_{\cal{E}}(\psi(x),x_{1})Y_{\cal{E}}(\phi(x),x_{2})\theta(x)\nonumber\\
&=&\left(f(x+x_{2},x)f(y,x)f(y,x+x_{2})
\psi(y)Y_{\cal{E}}(\phi (x),x_{2})\theta(x)\right)|_{y=x+x_{1}}\nonumber\\
&=&\left(f(z,x)f(y,x)f(y,z)
\psi(y)\phi (z)\theta(x)\right)|_{y=x+x_{1},z=x+x_{2}},
\end{eqnarray*}
completing the proof.
\end{proof}

Here we make a new notation for our convenience in the proof
of the following theorem.
For any formal series $a(x)=\sum_{n\in \Z}a_{n}x^{-n-1}$ 
(with coefficients $a_{n}$ in any vector space)
and for any $m\in \Z$, we set
\begin{eqnarray}
a(x)_{\ge m}=\sum_{n\ge m}a_{n}x^{-n-1}.
\end{eqnarray}
Then for any polynomial $q(x)$ we have
\begin{eqnarray}\label{esimplefactresidule}
\Res_{x}x^{m}q(x)a(x)=\Res_{x}x^{m}q(x)a(x)_{\ge m}.
\end{eqnarray}

Now we are in a position to prove our first key result:

\bt{tclosed}
Let $V$ be a subspace of $\E(W)$ such that any sequence
in $V$ of length $2$ or $3$ is quasi compatible and such that
$V$ contains $1_{W}$ and is closed.
Then $(V,Y_{\E},1_{W})$ carries the structure of a 
nonlocal vertex algebra with $W$ as a faithful quasi module 
where the vertex operator map $Y_{W}$ is given by
$Y_{W}(\alpha(x),x_{0})=\alpha(x_{0})$.
\et

\begin{proof} For the assertion on the nonlocal vertex algebra structure,
with Proposition \ref{pexistence},
it remains to prove the weak associativity, i.e.,
for $\psi,\phi,\theta\in V$, there exists a nonnegative 
integer $k$ such that
\begin{eqnarray}\label{eweakassocmainthem}
(x_{0}+x_{2})^{k}Y_{\cal{E}}(\psi,x_{0}+x_{2})Y_{\cal{E}}(\phi,x_{2})\theta
=(x_{0}+x_{2})^{k}Y_{\cal{E}}(Y_{\cal{E}}(\psi,x_{0})\phi,x_{2})\theta.
\end{eqnarray}
Let $f(x,y)$ be a nonzero polynomial such that
\begin{eqnarray*}
& &f(x,y)\psi(x)\phi(y)\in \Hom (W,W((x,y))),\\
& &f(x,y)\phi(x)\theta(y)\in \Hom (W,W((x,y))),\\
& &f(x,y)f(x,z)f(y,z)
\psi(x)\phi(y)\theta(z)\in \Hom (W,W((x,y,z))).
\end{eqnarray*}
By Lemma \ref{lclosed}, we have
\begin{eqnarray}\label{e5.76}
& &f(x+x_{2},x)f(x+x_{0}+x_{2},x)f(x+x_{0}+x_{2},x+x_{2})
Y_{\cal{E}}(\psi(x),x_{0}+x_{2})Y_{\cal{E}}(\phi(x),x_{2})\theta(x)
\nonumber\\
&=&\left(f(z,x)f(y,x)f(y,z)
\psi(y)\phi(z)\theta(x)\right)|_{y=x+x_{0}+x_{2},z=x+x_{2}}.
\end{eqnarray}

On the other hand, let $n\in \Z$ be {\em arbitrarily fixed}.
Since $\psi(x)_{m}\phi(x)=0$ for $m$ sufficiently large, there exists
a nonzero polynomial $p(x,y)$, depending on $n$, such that
\begin{eqnarray}\label{eges}
p(x+x_{2},x)
(Y_{\cal{E}}(\psi(x)_{m}\phi(x), x_{2})\theta(x)
=\left(p(z,x)(\psi(z)_{m}\phi(z))\theta(x)\right)|_{z=x+x_{2}}
\end{eqnarray}
for {\em all} $m\ge n$. 
With $f(x,y)\psi(x)\phi(y)\in \Hom (W,W((x,y)))$, from Definition \ref{d-operation} we have
\begin{eqnarray}\label{ethis}
f(x_{2}+x_{0},x_{2})(Y_{\cal{E}}(\psi(x_{2}),x_{0})\phi(x_{2}))\theta(x)
=\left(f(y,x_{2})\psi(y)\phi(x_{2})\theta(x)\right)|_{y=x_{2}+x_{0}}.
\end{eqnarray}
Using (\ref{esimplefactresidule}), (\ref{eges}) and (\ref{ethis}) we get
\begin{eqnarray}\label{e5.78}
& &\Res_{x_{0}}x_{0}^{n}f(x+x_{0}+x_{2},x)f(x+x_{0}+x_{2},x+x_{2})
p(x+x_{2},x)\nonumber\\
& &\ \ \ \ \cdot Y_{\cal{E}}(Y_{\cal{E}}(\psi(x),x_{0})\phi(x), x_{2})\theta(x)\nonumber\\
&=&\Res_{x_{0}}x_{0}^{n}f(x+x_{0}+x_{2},x)
f(x+x_{0}+x_{2},x+x_{2})p(x+x_{2},x)\nonumber\\
& &\ \ \ \ \cdot 
Y_{\cal{E}}(Y_{\cal{E}}(\psi(x),x_{0})_{\ge n}\phi(x), x_{2})\theta(x)\nonumber\\
&=&\Res_{x_{0}}x_{0}^{n}
f(x+x_{0}+x_{2},x)f(x+x_{0}+x_{2},x+x_{2})\nonumber\\
& &\ \ \ \ \cdot \left(p(z,x)
Y_{\cal{E}}(\psi(z),x_{0})_{\ge n}\phi(z))\theta(x)\right)|_{z=x+x_{2}}\nonumber\\
&=&\Res_{x_{0}}x_{0}^{n}
f(x+x_{0}+x_{2},x)f(x+x_{0}+x_{2},x+x_{2})
\nonumber\\
& &\ \ \ \ \cdot \left(p(z,x)Y_{\cal{E}}(\psi(z),x_{0})\phi(z))\theta(x)\right)|_{z=x+x_{2}}\nonumber\\
&=&\Res_{x_{0}}x_{0}^{n}\left(f(z+x_{0},x)f(z+x_{0},z)
p(z,x)(Y_{\cal{E}}(\psi(z),x_{0})\phi(z))\theta(x)\right)|_{z=x+x_{2}}\nonumber\\
&=&\Res_{x_{0}}x_{0}^{n}\left(f(y,x)f(y,z)
p(z,x)\psi(y)\phi(z)\theta(x)\right)|_{y=z+x_{0},z=x+x_{2}}
\nonumber\\
&=&\Res_{x_{0}}x_{0}^{n}\left(f(y,x)f(y,z)
p(z,x)\psi(y)\phi(z)\theta(x)\right)|_{y=x+x_{2}+x_{0},z=x+x_{2}}.
\end{eqnarray}
Combining (\ref{e5.78}) with (\ref{e5.76}) we get
\begin{eqnarray}\label{e5.79}
& &\Res_{x_{0}}x_{0}^{n}
f(x+x_{2},x)f(x+x_{0}+x_{2},x)f(x+x_{0}+x_{2},x+x_{2})
\nonumber\\
& &\ \ \ \ \cdot
p(x+x_{2},x)Y_{\cal{E}}(\psi(x),x_{0}+x_{2})
Y_{\cal{E}}(\phi(x),x_{2})\theta (x)\nonumber\\
&=&\Res_{x_{0}}x_{0}^{n}f(x_{2}+x,x)f(x+x_{0}+x_{2},x)f(x+x_{0}+x_{2},x+x_{2})
\nonumber\\
& &\ \ \ \ \cdot p(x+x_{2},x)
Y_{\cal{E}}(Y_{\cal{E}}(\psi(x),x_{0})\phi(x), x_{2})\theta(x).
\end{eqnarray}
Notice that both sides of (\ref{e5.79}) involve only finitely 
many negative powers of $x_{2}$.
Then multiplying both sides by 
$\iota_{x,x_{2}}p(x+x_{2},x)^{-1}f(x+x_{2},x)^{-1}$ we get
\begin{eqnarray*}
& &\Res_{x_{0}}x_{0}^{n}f(x+x_{0}+x_{2},x)
f(x+x_{0}+x_{2},x+x_{2})
Y_{\cal{E}}(\psi(x),x_{0}+x_{2})Y_{\cal{E}}(\phi(x),x_{2})\theta(x)\nonumber\\ 
&=&\Res_{x_{0}}x_{0}^{n}
f(x+x_{0}+x_{2},x)f(x+x_{0}+x_{2},x+x_{2})
Y_{\cal{E}}(Y_{\cal{E}}(\psi(x),x_{0})\phi(x), x_{2})\theta(x).\ \ \ \ 
\end{eqnarray*}
Since {\em $f(x,y)$ does not depend on $n$
and since $n$ is arbitrary}, we have
\begin{eqnarray}\label{enearfinal-new}
& &f(x+x_{0}+x_{2},x)f(x+x_{0}+x_{2},x+x_{2})
Y_{\cal{E}}(\psi(x),x_{0}+x_{2})Y_{\cal{E}}(\phi(x),x_{2})\theta(x)\nonumber\\
&=&f(x+x_{0}+x_{2},x)f(x+x_{0}+x_{2},x+x_{2})
Y_{\cal{E}}(Y_{\cal{E}}(\psi(x),x_{0})\phi(x), x_{2})\theta(x).
\end{eqnarray}
Write $f(x,y)=(x-y)^{k}g(x,y)$ for some $k\in \N,\; g(x,y)\in \C[x,y]$
with $g(x,x)\ne 0$.
Then
\begin{eqnarray*}
& &f(x+x_{0}+x_{2},x)=(x_{0}+x_{2})^{k}g(x+x_{0}+x_{2},x),\\
& &f(x+x_{0}+x_{2},x+x_{2})=x_{0}^{k}g(x+x_{0}+x_{2},x+x_{2}).
\end{eqnarray*}
In view of Lemma \ref{lcancellation} (b) we obtain
\begin{eqnarray*}
(x_{0}+x_{2})^{k}(Y_{\cal{E}}(\psi(x),x_{0}+x_{2})Y_{\cal{E}}(\phi(x),x_{2})
\theta(x)=(x_{0}+x_{2})^{k}
Y_{\cal{E}}(Y_{\cal{E}}(\psi(x),x_{0})\phi(x),x_{2})\theta(x),
\end{eqnarray*}
as desired.

For $a(x),b(x)\in V,\; w\in W$, in view of 
Proposition \ref{pexistence} there exist a nonzero polynomial $h(x,y)$ and a
nonnegative integer $l$ such that
$$(x_{0}+x_{2})^{l}h(x_{0}+x_{2},x_{2})(Y_{\cal{E}}(a(x),x_{0})b(x))(x_{2})w=
(x_{0}+x_{2})^{l}h(x_{0}+x_{2},x_{2})a(x_{0}+x_{2})b(x_{2})w.$$
That is,
\begin{eqnarray}
& &(x_{0}+x_{2})^{l}h(x_{0}+x_{2},x_{2})Y_{W}(Y_{\cal{E}}(a(x),x_{0})b(x),x_{2})w
\nonumber\\
&=&(x_{0}+x_{2})^{l}h(x_{0}+x_{2},x_{2})Y_{W}(a(x),x_{0}+x_{2})Y_{W}(b(x),x_{2})w.
\end{eqnarray}
Therefore $W$ is a quasi $V$-module with 
$Y_{W}(\alpha(x),x_{0})=\alpha(x_{0})$
for $\alpha(x)\in V$.
\end{proof}

Next we are going to prove that any quasi compatible subset $S$ of $\E(W)$
gives rise to a nonlocal vertex algebra.
To achieve this goal, we shall need the following key result:

\bp{pgeneratingcomplicatedone}
Let $\psi_{1}(x),\dots,\psi_{r}(x), a(x),b(x),\phi_{1}(x), \dots,\phi_{s}(x)
\in \E(W)$. Assume that 
the ordered sequences $(a(x), b(x))$ and 
$(\psi_{1}(x),\dots,\psi_{r}(x), a(x),b(x),\phi_{1}(x), \dots,\phi_{s}(x))$
 are quasi compatible.
Then for any $n\in \Z$, the ordered sequence 
$$(\psi_{1}(x),\dots,\psi_{r}(x),a(x)_{n}b(x),\phi_{1}(x),\dots,\phi_{s}(x))$$
 is quasi compatible.
\ep

\begin{proof} 
Let $f(x,y)$ be a nonzero polynomial such that
$$f(x_{1},x_{2})a(x_{1})b(x_{2})\in \Hom (W,W((x_{1},x_{2})))$$
and
\begin{eqnarray}\label{elong-exp}
& &\left(\prod_{1\le i<j\le r}f(y_{i},y_{j})\right)
\left(\prod_{1\le i\le r, 1\le j\le s}f(y_{i},z_{j})\right)
\left(\prod_{1\le i<j\le s}f(z_{i},z_{j})\right)\nonumber\\
& &\;\;\cdot f(x_{1},x_{2})
\left(\prod_{i=1}^{r}f(x_{1},y_{i})f(x_{2},y_{i})\right)
\left(\prod_{i=1}^{s}f(x_{1},z_{i})f(x_{2},z_{i})\right)\nonumber\\
& &\;\;\cdot \psi_{1}(y_{1})\cdots \psi_{r}(y_{r})
a(x_{1})b(x_{2})\phi_{1}(z_{1})\cdots \phi_{s}(z_{s})\nonumber\\
& &\in \Hom (W,W((y_{1},\dots, y_{r},x_{1},x_{2},z_{1},\dots,z_{s}))).
\end{eqnarray}
Set
$$P=\prod_{1\le i<j\le r}f(y_{i},y_{j}),\;\;\;\; 
Q=\prod_{1\le i<j\le s}f(z_{i},z_{j}),\;\;\;\;
R=\prod_{1\le i\le r,\; 1\le j\le s}f(y_{i},z_{j}).$$
Let $n\in \Z$ be {\em arbitrarily fixed}. 
There exists a nonnegative integer $k$ such that 
\begin{eqnarray}\label{etruncationpsiphi}
x_{0}^{k+n}\iota_{x_{2},x_{0}}\left(f(x_{0}+x_{2},x_{2})^{-1}\right)
\in \C ((x_{2}))[[x_{0}]].
\end{eqnarray}
Using (\ref{etruncationpsiphi}) and Definition \ref{d-operation}
we obtain
\begin{eqnarray}\label{ecompatibilitythreeproof}
& &\prod_{i=1}^{r}f(x_{2},y_{i})^{k}
\prod_{j=1}^{s}f(x_{2},z_{j})^{k}\nonumber\\
& &\ \ \ \ \cdot \psi_{1}(y_{1})\cdots\psi_{r}(y_{r})
(a(x)_{n}b(x))(x_{2})\phi_{1}(z_{1})\cdots\phi_{s}(z_{s})\nonumber\\
&=&\Res_{x_{0}}x_{0}^{n}\prod_{i=1}^{r}f(x_{2},y_{i})^{k}
\prod_{j=1}^{s}f(x_{2},z_{j})^{k}\nonumber\\
& &\ \ \ \ \cdot \psi_{1}(y_{1})\cdots \psi_{r}(y_{r})
(Y_{\cal{E}}(a,x_{0})b)(x_{2})\phi_{1}(z_{1})\cdots \phi_{s}(z_{s})
\nonumber\\
&=&\Res_{x_{1}}\Res_{x_{0}}x_{0}^{n}\prod_{i=1}^{r}f(x_{2},y_{i})^{k}
\prod_{j=1}^{s}f(x_{2},z_{j})^{k}\nonumber\\
& &\ \ \ \ \cdot \iota_{x_{2},x_{0}}(f(x_{2}+x_{0},x_{2})^{-1})
x_{1}^{-1}\delta\left(\frac{x_{2}+x_{0}}{x_{1}}\right)\nonumber\\
& &\ \ \ \ \cdot \left(f(x_{1},x_{2})\psi_{1}(y_{1})\cdots \psi_{r}(y_{r})
a(x_{1})b(x_{2})\phi_{1}(z_{1})\cdots
\phi_{s}(z_{s})\right)\nonumber\\
&=&\Res_{x_{1}}\Res_{x_{0}}x_{0}^{n}\prod_{i=1}^{r}f(x_{1}-x_{0},y_{i})^{k}
\prod_{j=1}^{s}f(x_{1}-x_{0},z_{j})^{k}\nonumber\\
& &\ \ \ \ \cdot \iota_{x_{2},x_{0}}(f(x_{2}+x_{0},x_{2})^{-1})
x_{1}^{-1}\delta\left(\frac{x_{2}+x_{0}}{x_{1}}\right)\nonumber\\
& &\ \ \ \ \cdot \left(f(x_{1},x_{2})\psi_{1}(y_{1})\cdots \psi_{r}(y_{r})
a(x_{1})b(x_{2})\phi_{1}(z_{1})\cdots
\phi_{s}(z_{s})\right)\nonumber\\
&=&\Res_{x_{1}}\Res_{x_{0}}x_{0}^{n}
e^{-x_{0}\frac{\partial}{\partial x_{1}}}
\left(\prod_{i=1}^{r}f(x_{1},y_{i})
\prod_{j=1}^{s}f(x_{1},z_{j})\right)^{k}\nonumber\\
& &\ \ \ \ \cdot \iota_{x_{2},x_{0}}(f(x_{2}+x_{0},x_{2})^{-1})
x_{1}^{-1}\delta\left(\frac{x_{2}+x_{0}}{x_{1}}\right)\nonumber\\
& &\ \ \ \ \cdot \left(f(x_{1},x_{2})\psi_{1}(y_{1})\cdots \psi_{r}(y_{r})
a(x_{1})b(x_{2})\phi_{1}(z_{1})\cdots
\phi_{s}(z_{s})\right)\nonumber\\
&=&\Res_{x_{1}}\Res_{x_{0}}\sum_{t=0}^{k-1}\frac{(-1)^{t}}{t!}x_{0}^{n+t}
\left(\frac{\partial}{\partial x_{1}}\right)^{t}
\left(\prod_{i=1}^{r}f(x_{1},y_{i})
\prod_{j=1}^{s}f(x_{1},z_{j})\right)^{k}\nonumber\\
& &\ \ \ \ \cdot \iota_{x_{2},x_{0}}(f(x_{2}+x_{0},x_{2})^{-1})
x_{1}^{-1}\delta\left(\frac{x_{2}+x_{0}}{x_{1}}\right)\nonumber\\
& &\ \ \ \ \cdot \left(f(x_{1},x_{2})\psi_{1}(y_{1})\cdots \psi_{r}(y_{r})
a(x_{1})b(x_{2})\phi_{1}(z_{1})\cdots
\phi_{s}(z_{s})\right).
\end{eqnarray}
Notice that for any polynomial $B$ and for $0\le t\le k-1$,
$\left(\frac{\partial}{\partial x_{1}}\right)^{t}B^{k}$
is a multiple of $B$.
Using (\ref{elong-exp}) we have
\begin{eqnarray*}
& &P Q R \prod_{i=1}^{r}f(x_{2},y_{i})
\prod_{j=1}^{s}f(x_{2},z_{j})
\sum_{t=0}^{k-1}\frac{(-1)^{t}}{t!}x_{0}^{n+t}
\left(\frac{\partial}{\partial x_{1}}\right)^{t}
\left(\prod_{i=1}^{r}f(x_{1},y_{i})
\prod_{j=1}^{s}f(x_{1},z_{j})\right)^{k}\nonumber\\
& &\ \ \ \ \cdot \iota_{x_{2},x_{0}}(f(x_{2}+x_{0},x_{2})^{-1})
x_{1}^{-1}\delta\left(\frac{x_{2}+x_{0}}{x_{1}}\right)\nonumber\\
& &\ \ \ \ \cdot \left(f(x_{1},x_{2})\psi_{1}(y_{1})\cdots \psi_{r}(y_{r})
a(x_{1})b(x_{2})\phi_{1}(z_{1})\cdots
\phi_{s}(z_{s})\right)\nonumber\\
&\in& \left(\Hom (W,W((y_{1},\dots,
y_{r},x_{2},z_{1},\dots,z_{s}))\right)
((x_{0}))[[x_{1}^{\pm 1}]].
\end{eqnarray*}
Then
\begin{eqnarray}
& &P Q R\prod_{i=1}^{r}f(x_{2},y_{i})^{k+1}
\prod_{j=1}^{s}f(x_{2},z_{j})^{k+1}
\psi_{1}(y_{1})\cdots\psi_{r}(y_{r})
(a(x)_{n}b(x))(x_{2})\phi_{1}(z_{1})\cdots\phi_{s}(z_{s})\nonumber\\
&\in& \Hom (W,W((y_{1},\dots, y_{r},x_{2},z_{1},\dots,z_{s}))).
\end{eqnarray}
This proves that the sequence 
$(\psi_{1}(x),\dots,\psi_{r}(x),a(x)_{n}b(x), \phi_{1}(x),
\dots,\phi_{s}(x))$ is quasi compatible.
\end{proof}

Now we have:

\bt{tmaximal}
Let $V$ be a maximal quasi compatible subspace of $\E(W)$.
Then $V$ contains $1_{W}$ and is closed and 
$(V,Y_{\cal{E}},1_{W})$ carries the structure of a 
nonlocal vertex algebra with $W$ as a faithful quasi module where the vertex
operator map $Y_{W}$ is given by $Y_{W}(\psi(x),x_{0})=\psi(x_{0})$.
Furthermore, for $a(x),b(x)\in V,\; \alpha\in \C^{\times},\; n\in \Z$ we have
\begin{eqnarray}
& &R_{\alpha}(a(x))\;(=a(\alpha x))\in V,
\label{emaximal-closed-C}\\
& &a(x)_{(\alpha,n)}b(x)\in V.\label{emaximal-closed-CZ}
\end{eqnarray}
\et

\begin{proof} 
Clearly the space spanned by $V$ and 
$1_{W}$ is still quasi compatible. With $V$ being
maximal we must have $1_{W}\in V$.
Now, let $a(x),b(x)\in V$ and $n\in \Z$. In view of 
Proposition \ref{pgeneratingcomplicatedone},
any (ordered) sequence in $V\cup \{ a(x)_{n}b(x)\}$
with one appearance of $a(x)_{n}b(x)$ is quasi compatible.
It follows from induction on the number of 
appearance of $a(x)_{n}b(x)$ and from
Proposition \ref{pgeneratingcomplicatedone} that
any (ordered) sequence in $V\cup \{ a(x)_{n}b(x)\}$
with any (finite) number of appearance of $a(x)_{n}b(x)$ is quasi compatible.
So the space spanned by $V$ and $a(x)_{n}b(x)$ is quasi compatible.
Again, with $V$ being maximal we must have $a(x)_{n}b(x)\in V$. 
This proves that $V$ is closed, and hence by
Theorem \ref{tclosed} $(V,Y_{\cal{E}},1_{W})$ carries the structure 
of a nonlocal vertex algebra with $W$ as a faithful quasi module.

Since $V$ is quasi compatible, it is clear that $\{ a(\alpha x)\;|\;
a(x)\in V,\; \alpha\in \C^{\times}\}$, containing $V$, is also quasi compatible. With
$V$ being maximal we must have $V=\{ a(\alpha x)\;|\; a(x)\in V,\;
\alpha\in \C^{\times}\}$.  This proves (\ref{emaximal-closed-C}). 
By Proposition \ref{pbasic-property2} we have
$$Y_{\E}^{(\alpha)}(a(x),x_{0})b(x)=R_{\alpha}Y_{\E}(a(x),x_{0})R_{\alpha^{-1}}b(x)$$
for $\alpha\in \C^{\times},\; a(x),b(x)\in V$.
Then (\ref{emaximal-closed-CZ}) follows immediately.
\end{proof}

Furthermore, we have (cf. \cite{li-g1}, \cite{li-gamma}):

\bt{tgeneratingthem}
Let $S$ be any quasi compatible subset of $\E(W)$ and 
let $\Gamma$ be any subgroup of $\C^{\times}$.
There exists a unique smallest closed quasi compatible subspace 
denoted by $\<S\>_{\Gamma}$ which contains $S\cup \{1_{W}\}$ and which is
closed under the actions of $R_{\alpha}$ for $\alpha\in \Gamma$.
The triple $(\<S\>_{\Gamma},Y_{\E},1_{W})$ carries the structure of 
a nonlocal vertex algebra with $W$ as a faithful quasi module 
with $Y_{W}(\psi(x),x_{0})=\psi(x_{0})$.
Furthermore, we have
\begin{eqnarray}\label{e-gamma}
\<S\>_{\Gamma}={\rm span }\{ a^{(1)}(x)_{n_{1}}\cdots a^{(r)}(x)_{n_{r}}1_{W}
\;|\; r\ge 0,\; a^{(i)}(x)\in R_{\Gamma}(S),\; n_{i}\in \Z\},
\end{eqnarray}
where $R_{\Gamma}(S)=\{ R_{\alpha}a(x)=a(\alpha x)\;|\; \alpha\in \Gamma,\; a(x)\in S\}$,
and for $a(x),b(x)\in \<S\>_{\Gamma},\; \alpha\in \Gamma,\; n\in \Z$, we have
\begin{eqnarray}
& &a(x)_{(\alpha,n)}b(x)\in \<S\>_{\Gamma},\label{eRgamma-closedness}\\
& &Y_{\E}^{(\alpha)}(a(x),x_{0})b(x)=R_{\alpha}Y_{\E}(a(x),x_{0})R_{\alpha^{-1}}b(x),
\label{econjugation1}\\
& &R_{\alpha}Y_{\E}(a(x),x_{0})R_{\alpha^{-1}}b(x)=Y_{\E}(R_{\alpha}a(x),\alpha^{-1}x_{0})b(x).
\label{econjugation}
\end{eqnarray}
\et

\begin{proof} In view of Zorn's lemma there exists a maximal quasi compatible subspace $V$ of 
$\E(W)$, containing $S$ and $1_{W}$. By Theorem \ref{tmaximal},
$V$ is closed, $(V,Y_{\cal{E}},1_{W})$ carries the structure of a 
nonlocal vertex algebra with $W$ as a faithful quasi module, and $V$ is closed 
under the actions of $R_{\alpha}$ for $\alpha\in \C^{\times}$.
Now the existence and uniqueness of $\<S\>_{\Gamma}$ is clear and
the second assertion is also proved.

Note that the last two identities have already been established 
in Proposition \ref{pbasic-property2}.
Then (\ref{eRgamma-closedness}) follows from (\ref{econjugation1}).
For (\ref{e-gamma}), it is clear that $\<S\>_{\Gamma}$ contains the space, say $T$,
on the right-hand side of (\ref{e-gamma}).
By Lemma \ref{lsubalgebragenerated}, 
$T$ is the nonlocal vertex subalgebra 
generated by $R_{\Gamma}(S)$, so it is closed.
It follows from (\ref{econjugation}) that 
$T$ is also closed under the actions of
$R_{\alpha}$ for $\alpha\in \Gamma$. With $\<S\>_{\Gamma}$ being minimal
we must have $\<S\>_{\Gamma}\subset T$. This proves $\<S\>_{\Gamma}=T$,
proving (\ref{e-gamma}).
\end{proof}

Just as with classical associative (or Lie) algebras, a quasi module
amounts to a ``representation'' for a nonlocal vertex algebra.

\bp{pextra-need}
Let $V$ be a nonlocal vertex algebra and let 
$W$ be a vector space equipped with a linear map $Y_{W}$ from $V$ to 
$\E(W)$ with $Y_{W}({\bf 1},x)=1_{W}$.
If $(W,Y_{W})$ is a quasi $V$-module, then
for $u,v\in V$, whenever $(Y_{W}(u,x),Y_{W}(v,x))$ is quasi compatible, we have
\begin{eqnarray}\label{ehom-property}
Y_{W}(Y(u,x_{0})v,x)=Y_{\E}(Y_{W}(u,x),x_{0})Y_{W}(v,x).
\end{eqnarray}
On the other hand, if for any $u,v\in V$, $(Y_{W}(u,x),Y_{W}(v,x))$ is quasi compatible
and (\ref{ehom-property}) holds, then $(W,Y_{W})$ is a quasi $V$-module.
\ep

\begin{proof} Assume that $(W,Y_{W})$ is a quasi $V$-module and 
let $u,v\in V$ be such that the ordered pair $(Y_{W}(u,x),Y_{W}(v,x))$ is quasi compatible.
Let $w\in W$. By Proposition \ref{pexistence}, there exists 
 a nonzero polynomial $p(x_{1},x_{2})$ such that
$$p(x_{0}+x,x)Y_{W}(u,x_{0}+x)Y_{W}(v,x)w
=p(x_{0}+x,x)\left(Y_{\E}(u(x),x_{0})v(x)\right)w$$
and by the definition of a quasi module, there exists 
a nonzero polynomial $q(x_{1},x_{2})$ such that
$$q(x_{0}+x,x)Y_{W}(u,x_{0}+x)Y_{W}(v,x)w
=q(x_{0}+x,x)Y_{W}(Y(u,x_{0})v,x)w.$$
Consequently,
\begin{eqnarray*}
& &p(x_{0}+x,x)q(x_{0}+x,x)Y_{W}(Y(u,x_{0})v,x)w\\
&=&p(x_{0}+x,x)q(x_{0}+x,x)(Y_{\E}(u(x),x_{0})v(x))w.
\end{eqnarray*}
Since $Y_{W}(Y(u,x_{0})v,x)w,\; (Y_{\E}(u(x),x_{0})v(x))w\in W((x))((x_{0}))$,
by cancelation we get
$$Y_{W}(Y(u,x_{0})v,x)w=(Y_{\E}(u(x),x_{0})v(x))w,$$
proving (\ref{ehom-property}).

Assume that for any $u,v\in V$ $(Y_{W}(u,x),Y_{W}(v,x))$ 
is quasi compatible and (\ref{ehom-property}) holds.
For $u,v\in V,\; w\in W$, from Proposition \ref{pexistence} 
there exists a nonzero polynomial $p(x,y)$ such that
$$p(x_{0}+x,x)\left(Y_{\E}(Y_{W}(u,x),x_{0})Y_{W}(v,x)\right)w
=p(x_{0}+x,x)Y_{W}(u,x_{0}+x)Y_{W}(v,x)w.$$ 
Using (\ref{ehom-property}) we get
$$p(x_{0}+x,x)Y_{W}(Y(u,x_{0})v,x)w=p(x_{0}+x,x)Y_{W}(u,x_{0}+x)Y_{W}(v,x)w.$$ 
This proves that $(W,Y_{W})$ is a quasi $V$-module.
\end{proof}

Considering a $V$-module instead of a quasi module we have:

\bp{pextra-need-2}
Let $V$ be a nonlocal vertex algebra and let $W$ be a vector space 
equipped with a linear map $Y_{W}$ from $V$ to $\E(W)$ with $Y_{W}({\bf 1},x)=1_{W}$. 
Then $(W,Y_{W})$ is a $V$-module
if and only if for $u,v\in V$, there exists a nonnegative integer $k$ such that
\begin{eqnarray}
& &(x_{1}-x_{2})^{k}Y_{W}(u,x_{1})Y_{W}(v,x_{2})\in \Hom(W,W((x_{1},x_{2}))),
\label{ecomp-module}\\
& &Y_{W}(Y(u,x_{0})v,x)=Y_{\E}(Y_{W}(u,x),x_{0})Y_{W}(v,x).\label{emodule-hom}
\end{eqnarray}
\ep

\begin{proof} Assume that $(W,Y_{W})$ is a $V$-module.
For $u,v\in W$, let $k$ be a nonnegative integer such that
$x^{k}Y(u,x)v\in V[[x]]$. For any $w\in W$, there exists a nonnegative integer $l$ such that
$$x_{0}^{k}(x_{0}+x_{2})^{l}Y_{W}(u,x_{0}+x_{2})Y_{W}(v,x_{2})w=
x_{0}^{k}(x_{0}+x_{2})^{l}Y_{W}(Y(u,x_{0})v,x_{2})w.$$
Noticing that the right hand side contains only nonnegative powers of $x_{0}$, we have 
$$x_{0}^{k}(x_{0}+x_{2})^{l}Y_{W}(u,x_{0}+x_{2})Y_{W}(v,x_{2})w
\in W[[x_{0},x_{2}]][x_{2}^{-1}].$$
Thus
$$(x_{1}-x_{2})^{k}x_{1}^{l}Y_{W}(u,x_{1})Y_{W}(v,x_{2})w
\in W[[x_{1},x_{2}]][x_{2}^{-1}],$$
which implies
$$(x_{1}-x_{2})^{k}Y_{W}(u,x_{1})Y_{W}(v,x_{2})w
\in W((x_{1},x_{2})).$$
As $k$ is independent of $w$, we have
$$(x_{1}-x_{2})^{k}Y_{W}(u,x_{1})Y_{W}(v,x_{2})\in \Hom (W,W((x_{1},x_{2}))).$$
By Proposition \ref{pextra-need} (\ref{emodule-hom}) holds. 

Conversely, for $u,v\in V$ and $w\in W$, let $k$ be a nonnegative integer such that
(\ref{ecomp-module}) holds. {}From Proposition \ref{pexistence} 
there exists a nonnegative integer $l$ such that
$$(x_{0}+x)^{l}x_{0}^{k}\left(Y_{\E}(Y_{W}(u,x),x_{0})Y_{W}(v,x)\right)w
=(x_{0}+x)^{l}x_{0}^{k}Y_{W}(u,x_{0}+x)Y_{W}(v,x)w.$$ 
As $Y_{\E}(Y_{W}(u,x),x_{0})Y_{W}(v,x)=Y_{W}(Y(u,x_{0})v,x)$ by assumption, 
it follows that
$$(x_{0}+x)^{l}Y_{W}(Y(u,x_{0})v,x)w=(x_{0}+x)^{l}Y_{W}(u,x_{0}+x)Y_{W}(v,x)w.$$ 
This proves that $(W,Y_{W})$ is a $V$-module.
\end{proof}

\br{rrep}
{\em The properties (\ref{ecomp-module}) and (\ref{emodule-hom})
amount to that $(Y_{W}(V),Y_{\E},1_{W})$ is a nonlocal vertex algebra and
$Y_{W}$ is a nonlocal vertex algebra homomorphism. }
\er

\section{Pseudo-local subsets of $\E(W)$}

In this section we study certain special quasi compatible subspaces,
which we call pseudo-local subsets, $\S$-local subsets and quasi
$\Gamma$-local subsets of $\E(W)$. The main result is that if $U$ is a
pseudo-local subset, or an $\S$-local subset, or a quasi
$\Gamma$-local subset of $\E(W)$, then the nonlocal vertex algebra
$\<U\>$ generated by $U$ is pseudo-local, or $\S$-local, or
quasi $\Gamma$-local, respectively.  The notions of quasi
$\Gamma$-locality and $\S$-locality are designed with quantum affine algebras
and Yangians as targets.

As in Section 2, we fix a vector space $W$ over $\C$ throughout this section.

\bd{dsemi-local}
{\em A subset $U$ of $\E(W)$ is said to be {\em pseudo-local}
if for any $a(x),b(x)\in U$, there exist 
$0\ne p(x,y)\in \C[x^{\pm 1},y^{\pm 1}]$, $q_{i}(x,y)\in \C_{*}(x,y)$ and 
$u_{i}(x),v_{i}(x)\in U$ 
for $i=1,\dots,r$ (finite) such that
\begin{eqnarray}\label{esemi-pab}
&&p(x_{1},x_{2})a(x_{1})b(x_{2})
=\sum_{i=1}^{r}\iota_{x_{2},x_{1}}(q_{i}(x_{1},x_{2}))
u_{i}(x_{2})v_{i}(x_{1}).
\end{eqnarray}}
\ed

We have:

\bl{lpracticalcase}
Every pseudo-local subset $U$ of $\E(W)$ is quasi compatible.
\el

\begin{proof} We must prove that any sequence in $U$ of finite length
is quasi compatible. We shall use induction on the length $n$ of sequences.
Let $(a(x),b(x))$ be an ordered pair (a sequence of length $2$)
in $U$. By assumption, there exist 
$0\ne p(x,y)\in \C[x^{\pm 1},y^{\pm 1}]$,
$a^{(i)}(x), b^{(i)}(x)\in U$ and
$\bar{q}_{i}(x,y)\in \C((y))((x))$ for $i=1,\dots,r$ such that 
\begin{eqnarray}\label{e4.3}
p(x_{1},x_{2})a(x_{1})b(x_{2})
=\sum_{i=1}^{r}\bar{q}_{i}(x_{1},x_{2})b^{(i)}(x_{2})a^{(i)}(x_{1}).
\end{eqnarray}
The expression on the left-hand side lies in 
$\Hom (W,W((x_{1}))((x_{2})))$
and the expression on the right-hand side lies 
in $\Hom (W,W((x_{2}))((x_{1})))$. This forces the expressions
on both sides to lie in $\Hom (W,W((x_{1},x_{2})))$.
Thus $(a(x),b(x))$ is quasi compatible, proving the case for $n=2$.

Assume that $n\ge 2$ and any sequence  in $S$ of length $n$ is quasi compatible.
Let $\psi^{(1)}(x),\dots,\psi^{(n+1)}(x)\in S$. From the inductive hypothesis,
there exists  $0\ne f(x,y)\in \C[x^{\pm 1},y^{\pm 1}]$ such that
\begin{eqnarray}\label{eproduct1}
\left(\prod_{2\le i<j\le n+1}f(x_{i},x_{j})\right)
\psi^{(2)}(x_{2})\cdots \psi^{(n+1)}(x_{n+1})
\in \Hom (W,W((x_{2},\dots,x_{n+2}))).
\end{eqnarray}
By assumption there exist 
$0\ne p(x,y)\in \C[x^{\pm 1},y^{\pm 1}]$,
$a^{(i)}(x), b^{(i)}(x)\in U$ and $\bar{q}_{i}(x,y)\in \C((y))((x))$ 
for $i=1,\dots,r$ such that
\begin{eqnarray}\label{eproduct2}
p(x_{1},x_{2})\psi^{(1)}(x_{1})\psi^{(2)}(x_{2})
=\sum_{i=1}^{r}\bar{q}_{i}(x_{1},x_{2})b^{(i)}(x_{2})a^{(i)}(x_{1}).
\end{eqnarray}
{}From the inductive hypothesis again, there exists
$0\ne g(x,y)\in \C[x^{\pm 1},y^{\pm 1}]$ such that
\begin{eqnarray}\label{eproduct3}
& &\left(
\prod_{1\le i<j\le n+1,\; i,j\ne 2}g(x_{i},x_{j})\right)
a^{(s)}(x_{1})\psi^{(3)}(x_{3})\cdots \psi^{(n+1)}(x_{n+1})
\nonumber\\
&\in& \Hom (W, W((x_{1},x_{3},x_{4},\dots,x_{n+1})))
\end{eqnarray}
for $s=1,\dots,r$. Using (\ref{eproduct2}) we have
\begin{eqnarray}\label{eproduct4}
& &\left(\prod_{2\le i<j\le n+1}f(x_{i},x_{j})
\prod_{1\le i<j\le n+1,\; i,j\ne 2}g(x_{i},x_{j})\right)
p(x_{1},x_{2})\psi^{(1)}(x_{1})\cdots \psi^{(n+1)}(x_{n+1})\nonumber\\
&=&\left(\prod_{2\le i<j\le n+1}f(x_{i},x_{j})
\prod_{1\le i<j\le n+1,\; i,j\ne 2}g(x_{i},x_{j})\right)\nonumber\\
& &\ \ \ \ \ \cdot \sum_{s=1}^{r}\bar{q}_{i}(x_{1},x_{2})
b^{(s)}(x_{2})a^{(s)}(x_{1})\psi^{(3)}(x_{3})\cdots
\psi^{(n+1)}(x_{n+1}).
\end{eqnarray}
{}From (\ref{eproduct1}), the expression on the left-hand side
of (\ref{eproduct4}) lies in
$$\Hom (W,W((x_{1}))((x_{2},x_{3},x_{4},\dots,x_{n+1}))),$$
and by (\ref{eproduct3}), the expression on the right-hand side 
of (\ref{eproduct4}) lies in
$$\Hom (W,W((x_{2}))((x_{1},x_{3},x_{4},\dots,x_{n+1}))).$$
This forces the
expressions on both sides to lie in
$\Hom (W,W((x_{1},x_{2},x_{3},x_{4},\dots,x_{n+1}))).$
In particular,  the
expression on the left-hand side to lie in
$\Hom (W,W((x_{1},x_{2},x_{3},x_{4},\dots,x_{n+1}))).$
This proves that the sequence $(\psi^{(1)}(x),\dots,\psi^{(n+1)}(x))$ 
is quasi compatible, completing the induction.
\end{proof}

\bd{dquasi-s-local}
{\em A subset $U$ of $\E(W)$ is said to be {\em $\S$-local} if 
for any $a(x),b(x)\in U$, there exist 
$\psi_{i}(x)\in \C((x)),\; u_{i}(x),v_{i}(x)\in U$
 for $i=1,\dots,r$ and a nonnegative integer $k$ such that
\begin{eqnarray}\label{esr-local-simple}
& &(x_{1}-x_{2})^{k}a(x_{1})b(x_{2})
=\sum_{i=1}^{r}(x_{1}-x_{2})^{k}
\iota_{x_{2},x_{1}}(\psi_{i}(x_{2}-x_{1})) u_{i}(x_{2})v_{i}(x_{1})
\nonumber\\
& &\hspace{3cm}\left(=\sum_{i=1}^{r}
\iota_{x_{2},x_{1}}((-1)^{k}(x_{2}-x_{1})^{k}
\psi_{i}(x_{2}-x_{1})) u_{i}(x_{2})v_{i}(x_{1})\right).
\end{eqnarray}}
\ed

Let $\Gamma$ be any subgroup of $\C^{\times}$. 
Set
\begin{eqnarray}
& &\C_{\Gamma}[x]=\< 1, x-\alpha\; |\; \alpha\in \Gamma\>\subset \C[x],\\
& &\C_{\Gamma}(x,y)=\C[[x,y]][x^{-1},y^{-1},(x-\alpha y)^{-1}\;|\;
\alpha\in \Gamma]\subset \C_{*}(x,y).
\end{eqnarray}
Then $\C_{\Gamma}(x,y)$ is a subalgebra and it is
closed under the partial differentiation operators $\partial/\partial
x$ and $\partial/\partial y$. 

\bd{dpseudo-gamma-local}
{\em Let $\Gamma$ be a subgroup of $\C^{\times}$.
A subset (or subspace) $U$ of $\E(W)$ is said to be 
{\em pseudo $\Gamma$-local} if 
for any $a(x),b(x)\in U$, there exist 
$0\ne p(x)\in \C_{\Gamma}[x]$ and 
$\psi_{i}(x,y)\in \C_{\Gamma}(x,y),\; u_{i}(x),v_{i}(x)\in U$ 
for $i=1,\dots,r$ such that
\begin{eqnarray}\label{efab=sum-gamma-pseudo}
p(x_{1}/x_{2})a(x_{1})b(x_{2})
&=&\sum_{i=1}^{r}p(x_{1}/x_{2})\iota_{x_{2},x_{1}}(\psi_{i}(x_{1},x_{2})) 
u_{i}(x_{2})v_{i}(x_{1})\nonumber\\
& &\left(=\sum_{i=1}^{r}
\iota_{x_{2},x_{1}}(p(x_{1}/x_{2})\psi_{i}(x_{1},x_{2})) 
u_{i}(x_{2})v_{i}(x_{1})\right).
\end{eqnarray}}
\ed

Clearly, $\S$-local subsets and pseudo $\Gamma$-local subsets
of $\E(W)$ are pseudo-local. Therefore they are quasi compatible
by Lemma \ref{lpracticalcase}.

We shall need the following two technical results:

\bl{lpseudocommutative-general-key}
Let $U$ be a quasi compatible subspace of $\E (W)$
and let $a(x),b(x), c(x)\in U$. Assume that
$0\ne p_{1}(x,y)\in \C[x^{\pm 1},y^{\pm 1}],\;q_{1i}(x,y)\in \C((y))((x))$
and $u_{i}(x),v_{i}(x)\in U$ for $i=1,\dots,r$
such that
\begin{eqnarray}\label{efab}
p_{1}(x_{1},x_{2})c(x_{1})a(x_{2})
=\sum_{i=1}^{r}q_{1i}(x_{1},x_{2})u_{i}(x_{2})v_{i}(x_{1}),
\end{eqnarray}
and assume that $0\ne p_{2}(x,y)\in \C[x^{\pm 1},y^{\pm 1}]$, 
$q_{2j}(x,y)\in \C((y))((x))$
and $\psi_{ij}(x),\phi_{ij}(x)\in U$ 
for $1\le i\le r,\; 1\le j\le s$ such that
\begin{eqnarray}\label{egca}
p_{2}(x_{1},x_{2})v_{i}(x_{1})b(x_{2})
=\sum_{j=1}^{s}q_{2j}(x_{1},x_{2})\psi_{ij}(x_{2})\phi_{ij}(x_{1})
\end{eqnarray}
for $i=1,\dots,r$. 
Then for any $n\in \Z$, there exists a nonnegative integer $k$ such that
\begin{eqnarray}
& &p_{1}(x_{3},x)^{k}p_{2}(x_{3},x)c(x_{3})(a(x)_{n}b(x))\nonumber\\
&=&\sum_{m=0}^{k-1}\frac{(-1)^{m}}{m!}
\sum_{i=1}^{r}\sum_{j=1}^{s}\sum_{t\ge 0}\frac{1}{t!}\left(\frac{\partial}{\partial x}\right)^{t}
\left(q_{1i}(x_{3},x)g_{m}(x_{3},x)\right)\nonumber\\
& &\ \ \ \ \cdot q_{2j}(x_{3},x)
\left(u_{i}(x)_{n+m+t}\psi_{ij}(x)\right)\phi_{ij}(x_{3}),
\end{eqnarray}
where $g_{m}(x,y)\in \C[x^{\pm 1},y^{\pm 1}]$ for $0\le m<k$ are defined by
\begin{eqnarray}\label{e4.9proof}
\left(\left(\frac{\partial}{\partial x_{1}}\right)^{m}
p_{1}(x_{3},x_{1})^{k}\right)=p_{1}(x_{3},x_{1})g_{m}(x_{3},x_{1}).
\end{eqnarray}
\el

\begin{proof} Combining (\ref{efab}) and (\ref{egca}) we have
\begin{eqnarray}\label{ecomb}
& &p_{1}(x_{3},x_{1})p_{2}(x_{3},x)c(x_{3})a(x_{1})b(x)
\nonumber\\
&=&\sum_{i=1}^{r}\sum_{j=1}^{s}q_{1i}(x_{3},x_{1})q_{2j}(x_{3},x)
u_{i}(x_{1})\psi_{ij}(x)\phi_{ij}(x_{3}).
\end{eqnarray}
Since $U$ is quasi compatible,
there exists a nonzero polynomial $f(x,y)$ such that
\begin{eqnarray}
& &f(x_{1},x_{2})a(x_{1})b(x_{2}),\;\;
f(x_{1},x_{2})u_{i}(x_{1})\psi_{ij}(x_{2})\in \Hom (W,W((x_{1},x_{2})))
\end{eqnarray}
for $i=1,\dots,r,\; j=1,\dots,s$. From Definition \ref{d-operation} we have
\begin{eqnarray}
& &Y_{\E}(a(x),x_{0})b(x)
=\Res_{x_{1}}
\iota_{x,x_{0}}(f(x+x_{0},x)^{-1})
x_{1}^{-1}\delta\left(\frac{x+x_{0}}{x_{1}}\right)
(f(x_{1},x)a(x_{1})b(x)),\label{eyab-proof}\\
& &Y_{\E}(u_{i}(x),x_{0})\psi_{ij}(x)
=\Res_{x_{1}}\iota_{x,x_{0}}(f(x+x_{0},x)^{-1})
x_{1}^{-1}\delta\left(\frac{x+x_{0}}{x_{1}}\right)
(f(x_{1},x)u_{i}(x_{1})\psi_{ij}(x)).\nonumber\\
& &\label{eyupsi-proof}
\end{eqnarray}

For any fixed $n\in \Z$, let $k$ be a positive integer depending on
 $n$ such that
\begin{eqnarray}\label{exk+n}
x_{0}^{k+n}\iota_{x,x_{0}}(f(x+x_{0},x)^{-1})\in \C((x))[[x_{0}]].
\end{eqnarray}
(Recall that $\iota_{x,x_{0}}(f(x+x_{0},x)^{-1})\in \C((x))((x_{0}))$.)
Using (\ref{eyab-proof}), 
(\ref{exk+n}), (\ref{e4.9proof}), (\ref{ecomb}) and (\ref{eyupsi-proof})
in a sequence we obtain
\begin{eqnarray}
& &p_{1}(x_{3},x)^{k}p_{2}(x_{3},x)c(x_{3})\left(a(x)_{n}b(x)\right)
\nonumber\\
&=&\Res_{x_{0}}x_{0}^{n}p_{1}(x_{3},x)^{k}p_{2}(x_{3},x)c(x_{3})
\left(Y_{\E}(a(x),x_{0})b(x)\right)\nonumber\\
&=&\Res_{x_{0}}\Res_{x_{1}}x_{0}^{n}p_{1}(x_{3},x)^{k}p_{2}(x_{3},x)c(x_{3})
\iota_{x,x_{0}}(f(x+x_{0},x)^{-1})\nonumber\\
& &\ \ \cdot x_{1}^{-1}\delta\left(\frac{x+x_{0}}{x_{1}}\right)
(f(x_{1},x)a(x_{1})b(x))
\nonumber\\
&=&\Res_{x_{0}}\Res_{x_{1}}x_{0}^{n}p_{1}(x_{3},x_{1}-x_{0})^{k}p_{2}(x_{3},x)c(x_{3})
\iota_{x,x_{0}}(f(x+x_{0},x)^{-1})\nonumber\\
& &\ \ \cdot x_{1}^{-1}\delta\left(\frac{x+x_{0}}{x_{1}}\right)
(f(x_{1},x)a(x_{1})b(x))
\nonumber\\
&=&\Res_{x_{0}}\Res_{x_{1}}x_{0}^{n}
\left(e^{-x_{0}\frac{\partial}{\partial x_{1}}}p_{1}(x_{3},x_{1})^{k}\right)p_{2}(x_{3},x)c(x_{3})
\iota_{x,x_{0}}(f(x+x_{0},x)^{-1})\nonumber\\
& &\ \ \cdot x_{1}^{-1}\delta\left(\frac{x+x_{0}}{x_{1}}\right)
(f(x_{1},x)a(x_{1})b(x))
\nonumber\\
&=&\Res_{x_{0}}\Res_{x_{1}}\sum_{m=0}^{k-1}\frac{(-1)^{m}}{m!}x_{0}^{n+m}
\left(\left(\frac{\partial}{\partial
x_{1}}\right)^{m}p_{1}(x_{3},x_{1})^{k}\right)p_{2}(x_{3},x)
\iota_{x,x_{0}}(f(x+x_{0},x)^{-1})\nonumber\\
& &\ \ \cdot x_{1}^{-1}\delta\left(\frac{x+x_{0}}{x_{1}}\right)
(f(x_{1},x)c(x_{3})a(x_{1})b(x))\nonumber\\
&=&\Res_{x_{0}}\Res_{x_{1}}\sum_{m=0}^{k-1}\frac{(-1)^{m}}{m!}x_{0}^{n+m}
p_{1}(x_{3},x_{1})g_{m}(x_{3},x_{1})p_{2}(x_{3},x)
\iota_{x,x_{0}}(f(x+x_{0},x)^{-1})\nonumber\\
& &\ \ \cdot x_{1}^{-1}\delta\left(\frac{x+x_{0}}{x_{1}}\right)
(f(x_{1},x)c(x_{3})a(x_{1})b(x))\nonumber\\
&=&\Res_{x_{0}}\Res_{x_{1}}\sum_{m=0}^{k-1}\frac{(-1)^{m}}{m!}x_{0}^{n+m}
\sum_{i=1}^{r}\sum_{j=1}^{s}\iota_{x,x_{0}}(f(x+x_{0},x)^{-1})
x_{1}^{-1}\delta\left(\frac{x+x_{0}}{x_{1}}\right)\nonumber\\
& &\ \ \ \cdot
\left(q_{1i}(x_{3},x_{1})g_{m}(x_{3},x_{1})q_{2j}(x_{3},x)
f(x_{1},x)u_{i}(x_{1})\psi_{ij}(x)\phi_{ij}(x_{3})\right)\nonumber\\
&=&\Res_{x_{0}}\Res_{x_{1}}\sum_{m=0}^{k-1}\frac{(-1)^{m}}{m!}x_{0}^{n+m}
\sum_{i=1}^{r}\sum_{j=1}^{s}\iota_{x,x_{0}}(f(x+x_{0},x)^{-1})
x_{1}^{-1}\delta\left(\frac{x+x_{0}}{x_{1}}\right)\nonumber\\
& &\ \ \ \cdot
q_{1i}(x_{3},x+x_{0})g_{m}(x_{3},x+x_{0})q_{2j}(x_{3},x)
\left(f(x_{1},x)u_{i}(x_{1})\psi_{ij}(x)\phi_{ij}(x_{3})\right)\nonumber\\
&=&\Res_{x_{0}}\sum_{m=0}^{k-1}\frac{(-1)^{m}}{m!}x_{0}^{n+m}
\sum_{i=1}^{r}\sum_{j=1}^{s}
q_{1i}(x_{3},x+x_{0})g_{m}(x_{3},x+x_{0})\nonumber\\
& &\ \ \ \ \cdot q_{2j}(x_{3},x)
\left(Y_{\E}(u_{i}(x),x_{0})\psi_{ij}(x)\right)\phi_{ij}(x_{3})\nonumber\\
&=&\sum_{m=0}^{k-1}\frac{(-1)^{m}}{m!}
\sum_{i=1}^{r}\sum_{j=1}^{s}\sum_{t\ge 0}\frac{1}{t!}\left(\frac{\partial}{\partial x}\right)^{t}
\left(q_{1i}(x_{3},x)g_{m}(x_{3},x)\right)\nonumber\\
& &\ \ \ \ \cdot q_{2j}(x_{3},x)
\left(u_{i}(x)_{n+m+t}\psi_{ij}(x)\right)\phi_{ij}(x_{3}).
\end{eqnarray}
This completes the proof.
\end{proof}

Clearly, the following parallel assertion also holds:

\bl{lpseudocommutative-general-key-right}
Let $U$ be a quasi compatible subspace of $\E (W)$
and let $a(x),b(x), c(x)\in U$. Assume that
$0\ne p_{1}(x,y)\in \C[x^{\pm 1},y^{\pm 1}],\;q_{1i}(x,y)\in \C((y))((x))$
and $u_{i}(x), v_{i}(x)\in U$ for $i=1,\dots,r$
such that
\begin{eqnarray}\label{efab-right}
p_{1}(x_{1},x_{2})b(x_{1})c(x_{2})
=\sum_{i=1}^{r}q_{1i}(x_{1},x_{2})u_{i}(x_{2})v_{i}(x_{1}),
\end{eqnarray}
and assume that $0\ne p_{2}(x,y)\in \C[x^{\pm 1},y^{\pm 1}]$, 
$q_{2j}(x,y)\in \C((y))((x))$
and $\psi_{ij}(x),\phi_{ij}(x)\in U$ 
for $1\le i\le r,\; 1\le j\le s$ such that
\begin{eqnarray}\label{eright-gca}
p_{2}(x_{1},x_{2})a(x_{1})u_{i}(x_{2})
=\sum_{j=1}^{s}q_{2j}(x_{1},x_{2})\psi_{ij}(x_{2})\phi_{ij}(x_{1})
\end{eqnarray}
for $i=1,\dots,r$. 
Then for any $n\in \Z$, there exists a nonnegative integer $k$ such that
\begin{eqnarray}
& &p_{1}(x,x_{3})p_{2}(x,x_{3})^{k}(a(x)_{n}b(x))c(x_{3})\nonumber\\
&=&\sum_{m=0}^{k-1}\frac{(-1)^{m}}{m!}
\sum_{i=1}^{r}\sum_{j=1}^{s}\sum_{t\ge 0}\frac{1}{t!}\left(\frac{\partial}{\partial x}\right)^{t}
\left(q_{1i}(x,x_{3})g_{m}(x,x_{3})\right)\nonumber\\
& &\ \ \ \ \cdot q_{2j}(x,x_{3})
\psi_{ij}(x_{3})\left(\phi_{ij}(x)_{n+m+t}v_{i}(x)\right),
\end{eqnarray}
where $g_{m}(x,y)\in \C[x^{\pm 1},y^{\pm 1}]$ for $0\le m<k$ are defined by
\begin{eqnarray}\label{eright-commu}
\left(\left(\frac{\partial}{\partial x_{1}}\right)^{m}
p_{2}(x_{1},x_{3})^{k}\right)=p_{2}(x_{1},x_{3})g_{m}(x_{1},x_{3}).
\end{eqnarray}
\el

\bd{dLC(U)}
{\em Let $S$ and $U$ be subsets of $\E(W)$. 
Define $PC_{U}(S)$ to consist of 
$a(x)\in U$ satisfying the condition that
for any $b(x)\in S$, there exists 
$0\ne p(x,y)\in \C[x^{\pm 1},y^{\pm 1}]$ such that
\begin{eqnarray}\label{eLCU}
p(x_{1},x_{2})a(x_{1})b(x_{2})
=\sum_{i=1}^{r}\iota_{x_{2},x_{1}}(\psi_{i}(x_{1},x_{2}))
c_{i}(x_{2})d_{i}(x_{1})
\end{eqnarray}
for some $\psi_{i}(x,y)\in
\C_{*}(x,y),\;c_{i}(x)\in S, \;d_{i}(x)\in U$ for $i=1,\dots,r$
and such that \begin{eqnarray}\label{eRCU}
p(x_{1},x_{2})b(x_{1})a(x_{2})
=\sum_{j=1}^{s}\iota_{x_{2},x_{1}}(\phi_{j}(x_{1},x_{2}))
v_{j}(x_{2})u_{j}(x_{1})
\end{eqnarray}
for some $\phi_{j}(x,y)\in \C_{*}(x,y)$ and 
$u_{j}(x)\in S, \;v_{j}(x)\in U$ for $j=1,\dots,s$.}
\ed

We also define $PC_{U}^{\S}(S)$ by requiring that
$p(x,y)\in \C[(x-y)],\;\psi_{i}(x,y),\phi_{j}(x,y)\in \C((x-y))$.
Furthermore, for any subgroup $\Gamma$ of $\C^{\times}$,
we define $PC_{U}^{\Gamma}(S)$
by requiring that $p(x,y)\in \C_{\Gamma}[x/y],\;
\psi_{i}(x,y),\phi_{j}(x,y) \in \C_{\Gamma}(x,y)$.

Clearly, if $U$ is a subspace of $\E(W)$, $PC_{U}(S)$,
$PC_{U}^{\S}(S)$, and $PC_{U}^{\Gamma}(S)$ are subspaces of $U$. 
{}From Lemmas \ref{lpseudocommutative-general-key} and 
\ref{lpseudocommutative-general-key-right} we immediately have:

\bc{cpseudocommutative-key}
Let $U$ be a closed quasi compatible subspace of $\E (W)$ and
let $S$ be a subset of $U$. 
Then $PC_{U}(S)$ and
$PC_{U}^{\S}(S)$ are closed subspaces of $U$. 
For any subgroup $\Gamma$ of $\C^{\times}$, the subspace
$PC_{U}^{\Gamma}(S)$ is also closed.
\ec

Recall from Lemma \ref{lpracticalcase} that 
any pseudo-local subset $S$ of $\E (W)$
is quasi compatible. By Theorem \ref{tgeneratingthem}, 
we have the smallest closed quasi
compatible subspace $\<U\>$ of $\E(W)$, containing $S\cup \{1_{W}\}$.

\bp{pgeneral-locality-key} 
If $S$ is a pseudo-local subset of $\E(W)$, then the nonlocal vertex algebra 
$\<S\>$ generated by $S$ is also pseudo-local.  
If $S$ is an $\S$-local subset of $\E(W)$, then $\<S\>$ is $\S$-local.  
\ep

\begin{proof} 
Taking $U=\<S\>$ in Corollary \ref{cpseudocommutative-key}, 
we see that
$PC_{U}(S)$, containing $S\cup \{1_{W}\}$, is a closed
subspace of $U$. Thus, $U=\<S\>\subset PC_{U}(S)$.
This implies that $S\subset PC_{U}(U)$.
Using Corollary \ref{cpseudocommutative-key} again 
we have $U=\<S\>\subset PC_{U}(U)$. Therefore,
$U=\<S\>= PC_{U}(U)$.
That is, $U=\<S\>$ is pseudo-local. 

For the second assertion, using
$PC_{U}^{\S}(S)$, instead of $PC_{U}(S)$ we see that
$\<S\>$ is an $\S$-local subspace of $\E(W)$.
\end{proof}

Let $\Gamma$ be a subgroup of $\C^{\times}$
and let $S$ be a pseudo $\Gamma$-local subset of $\E(W)$.
Let $\Gamma_{0}$ be a subgroup of $\Gamma$.
By Theorem \ref{tgeneratingthem}, we have the smallest closed quasi
compatible subspace $\<S\>_{\Gamma_{0}}$ of $\E(W)$,  
which contains $S\cup \{1_{W}\}$ 
and which is closed under the actions of
$R_{\alpha}$ for $\alpha\in \Gamma_{0}$. {}From the proof
of Proposition \ref{pgeneral-locality-key} we immediately have:

\bp{pgamma-locality-key}
Let $\Gamma$ be a subgroup of $\C^{\times}$
and let $S$ be a pseudo $\Gamma$-local subset of $\E (W)$.
Then for any subgroup $\Gamma_{0}$ of $\Gamma$,
$\<S\>_{\Gamma_{0}}$ is also pseudo $\Gamma$-local.
\ep

\br{rquasi-gamma-local}
{\em Let $\Gamma$ be a subgroup of $\C^{\times}$. Set
\begin{eqnarray*}
\C_{\Gamma}(x)=\C[[x]][x^{-1},(x-\alpha)^{-1}\;|\; \alpha\in \Gamma]
\subset \C(x).
\end{eqnarray*}
Let $U$ be a subset of $\E(W)$ satisfying the condition that
for any $a(x),b(x)\in U$, there exist 
$0\ne p(x)\in \C_{\Gamma}[x]$ and 
$\psi_{i}(x)\in \C_{\Gamma}(x),\; u_{i}(x),v_{i}(x)\in U$ 
for $i=1,\dots,r$ such that
\begin{eqnarray}\label{efab=sum-trigonometric}
p(x_{1}/x_{2})a(x_{1})b(x_{2})
=\sum_{i=1}^{r}\iota_{x_{2},x_{1}}(\psi_{i}(x_{1}/x_{2})) 
u_{i}(x_{2})v_{i}(x_{1}).
\end{eqnarray}
Clearly, $U$ is pseudo $\Gamma$-local.
By Proposition \ref{pgamma-locality-key},
the nonlocal vertex algebra $\<U\>_{\Gamma}$ 
is pseudo $\Gamma$-local. But, in general, $\<U\>_{\Gamma}$ does not satisfy 
the same (stronger pseudo $\Gamma$-locality) property
that $U$ satisfies. 
This is due to the fact that 
$\C_{\Gamma}[x/y]$ is not closed under the partial
differential operators with respect to $x$ and $y$ 
(cf. Lemmas \ref{lpseudocommutative-general-key} and 
\ref{lpseudocommutative-general-key-right}).}
\er

In Section 2, we have discussed the actions of $R_{\alpha}$ for $\alpha\in \C^{\times}$
on $\E(W)$ and its maximal quasi compatible subspaces. In fact,
there are also other natural actions as we discuss next.

Let $t$ be a formal variable. Then $\E(W)$ is naturally
a vector space over $\C((t))$ with
\begin{eqnarray}
f(t)a(x)=f(x)a(x)\;\;\;\mbox{ for }f(t)\in \C((t)),\; a(x)\in \E(W).
\end{eqnarray}

\bp{pnew-facts}
If $(a(x),b(x))$ is a quasi compatible pair
in $\E(W)$, then for any $f(t),g(t)\in \C((t))$, $(f(t)a(x),g(t)b(x))$ 
is a quasi compatible pair and
\begin{eqnarray}\label{et-linear-property}
Y_{\E}(f(t)a(x),x_{0})g(t)b(x)=f(t+x_{0})g(t)Y_{\E}(a(x),x_{0})b(x).
\end{eqnarray}
The $\C((t))$-linear span of any quasi compatible subset of $\E(W)$ is
quasi compatible and the $\C((t))$-linear span of any closed quasi compatible
$\C$-subspace of $\E(W)$ is closed.
Any maximal quasi compatible $\C$-subspace of $\E (W)$
is a $\C((t))$-subspace. 
\ep

\begin{proof} Let $a(x),b(x)\in \E(W)$ and $0\ne p(x,y)\in \C[x,y]$ be such that
$p(x_{1},x_{2})a(x_{1})b(x_{2})\in \Hom (W,W((x_{1},x_{2})))$. 
For any $f(t),g(t)\in \C((t))$, we have
$$p(x_{1},x_{2})(f(t)a(x_{1}))(g(t)b(x_{2}))=p(x_{1},x_{2})f(x_{1})g(x_{2})a(x_{1})b(x_{2})
\in \Hom (W,W((x_{1},x_{2}))).$$ 
By definition we have
\begin{eqnarray}
& &Y_{\E}(f(t)a(x),x_{0})(g(t)b(x))\nonumber\\
&=&\iota_{x,x_{0}}\left(p(x_{0}+x,x)^{-1}\right)
\left(p(x_{1},x)f(x_{1})g(x)a(x_{1})b(x)\right)|_{x_{1}=x+x_{0}}\nonumber\\
&=&f(x+x_{0})g(x)\iota_{x,x_{0}}\left(p(x_{0}+x,x)^{-1}\right)
\left(p(x_{1},x)a(x_{1})b(x)\right)|_{x_{1}=x+x_{0}}\nonumber\\
&=&f(t+x_{0})g(t)Y_{\E}(a(x),x_{0})b(x).
\end{eqnarray}
This proves the first assertion.

If $(\psi_{1}(x),\dots,\psi_{r}(x))$ is a
quasi compatible sequence in $\E(W)$, it is clear that
for any $f_{1}(t),\dots, f_{r}(t)\in \C((t))$, the sequence
$(f_{1}(x)\psi_{1}(x),\dots,f_{r}(x)\psi_{r}(x))$
is quasi compatible, 
that is, the sequence $(f_{1}(t)\psi_{1}(x),\dots,f_{r}(t)\psi_{r}(x))$
is quasi compatible. It follows that the $\C((t))$-linear span of 
any quasi compatible subset of $\E(W)$ is quasi compatible. This furthermore implies that
any maximal quasi compatible $\C$-subspace of $\E (W)$ is a $\C((t))$-subspace.
{}From (\ref{et-linear-property}), the $\C((t))$-linear span of any closed quasi compatible
$\C$-subspace of $\E(W)$ is closed.
\end{proof}

\bp{pconvert}
Let $V$ be a closed quasi compatible subspace of $\E(W)$.
Suppose that the following relation holds
\begin{eqnarray}
p(x_{1},x_{2})a(x_{1})b(x_{2})
=\sum_{i=1}^{r}\iota_{x_{2},x_{1}}(q_{i}(x_{1},x_{2}))
u_{i}(x_{2})v_{i}(x_{1}),
\end{eqnarray}
where $a(x),b(x),u_{i}(x),v_{i}(x)\in V$ and 
$p(x,y)\in \C[x,y],\; q_{i}(x,y)\in \C_{*}(x,y)$.
Then
\begin{eqnarray}\label{epqsum}
& &p(t+x_{1},t+x_{2})Y_{\E}(a(x),x_{1})Y_{\E}(b(x),x_{2})
\nonumber\\
&=&\sum_{i=1}^{r}\iota_{t,x_{2},x_{1}}(q_{i}(x_{1}+t,x_{2}+t))
Y_{\E}(u_{i}(x),x_{2})Y_{\E}(v_{i}(x),x_{1}).
\end{eqnarray}
Furthermore, if $q_{i}(x,y)=p(x,y)f_{i}(y-x)$ 
where $f_{i}(x)\in \C((x))$ for $1\le i\le r$, we have
\begin{eqnarray}\label{epfsum}
& &(x_{1}-x_{2})^{k}Y_{\E}(a(x),x_{1})Y_{\E}(b(x),x_{2})
\nonumber\\
&=&\sum_{i=1}^{r}(x_{1}-x_{2})^{k}\iota_{x_{2},x_{1}}(f_{i}(x_{2}-x_{1}))
Y_{\E}(u_{i}(x),x_{2})Y_{\E}(v_{i}(x),x_{1}),
\end{eqnarray}
where $k$ is the multiplicity of zero $x_{1}=x_{2}$ for $p(x_{1},x_{2})$.
\ep

\begin{proof} Let $\theta(x)\in V,\; w\in W$ be arbitrarily fixed vectors.
Let $f(x,y)$ be a nonzero polynomial such that
\begin{eqnarray*}
& &f(x,y)b(x)\theta(y)\in \Hom (W,W((x,y)),\\
& &f(x,y)f(x,z)f(y,z)a(x)b(y)\theta(z)\in \Hom (W,W((x,y,z)).
\end{eqnarray*}
Furthermore, let $l$ be a nonnegative integer such that
\begin{eqnarray*}
x^{l}y^{l}f(x,y)f(x,z)f(y,z)a(x)b(y)\theta(z)w\in W[[x,y,z]][z^{-1}].
\end{eqnarray*}
Using Lemma \ref{lclosed}, we have
\begin{eqnarray*}
& &(x_{1}+x)^{l}(x_{2}+x)^{l}f(x_{1}+x,x_{2}+x)f(x_{1}+x,x)f(x_{2}+x,x)\\
& &\ \ \ \ 
\cdot \left(Y_{\E}(a(x),x_{1})Y_{\E}(b(x),x_{2})\theta (x)\right)w\\
&=&\left(y^{l}z^{l}f(y,z)f(y,x)f(z,x)
a(y)b(z)\theta (x)w\right)|_{y=x+x_{1},z=x+x_{2}}\\
&=&\left(y^{l}z^{l}f(y,z)f(y,x)f(z,x)
a(y)b(z)\theta (x)w\right)|_{y=x_{1}+x,z=x_{2}+x}\\
&=&(x_{1}+x)^{l}(x_{2}+x)^{l}f(x_{1}+x,x_{2}+x)f(x_{1}+x,x)f(x_{2}+x,x)
a(x_{1}+x)b(x_{2}+x)\theta (x)w.
\end{eqnarray*}
Write $f(x,y)=(x-y)^{s}F(x,y)$ where $s\ge 0,\; F(x,y)\in \C[x,y]$
with $F(x,x)\ne 0$. Using cancelation we get
\begin{eqnarray*}
& &(x_{1}+x)^{l}(x_{2}+x)^{l}F(x_{1}+x,x_{2}+x)F(x_{1}+x,x)F(x_{2}+x,x)
\nonumber\\
& &\ \ \ \ 
\cdot \left(Y_{\E}(a(x),x_{1})Y_{\E}(b(x),x_{2})\theta (x)\right)w
\nonumber\\
&=&(x_{1}+x)^{l}(x_{2}+x)^{l}F(x_{1}+x,x_{2}+x)F(x_{1}+x,x)F(x_{2}+x,x)
a(x_{1}+x)b(x_{2}+x)\theta (x)w.
\end{eqnarray*}
It is clear that replacing $f(x,y)$ with a multiply of
$f(x,y)$ if necessary, we can also obtain
\begin{eqnarray*}
& &(x_{1}+x)^{l}(x_{2}+x)^{l}F(x_{1}+x,x_{2}+x)F(x_{1}+x,x)F(x_{2}+x,x)
\nonumber\\
& &\ \ \ \ \cdot
\left(Y_{\E}(u_{i}(x),x_{2})Y_{\E}(v_{i}(x),x_{1})\theta (x)\right)w
\nonumber\\
&=&(x_{1}+x)^{l}(x_{2}+x)^{l}F(x_{1}+x,x_{2}+x)F(x_{1}+x,x)F(x_{2}+x,x)
u_{i}(x_{2}+x)v_{i}(x_{1}+x)\theta (x)w
\end{eqnarray*}
for $i=1,\dots,r$.
Therefore, we obtain
\begin{eqnarray}
& &F(x_{1}+x,x_{2}+x)(x_{1}+x)^{l}F(x_{1}+x,x)(x_{2}+x)^{l}F(x_{2}+x,x)\nonumber\\
& &\ \ \ \ \cdot p(x_{1}+x,x_{2}+x)
\left(Y_{\E}(a(x),x_{1})Y_{\E}(b(x),x_{2})\theta (x)\right)w\nonumber\\
&=&F(x_{1}+x,x_{2}+x)(x_{1}+x)^{l}F(x_{1}+x,x)(x_{2}+x)^{l}F(x_{2}+x,x)\nonumber\\
& &\ \ \ \ \cdot \sum_{i=1}^{r}\iota_{x,x_{2},x_{1}}(q_{i}(x+x_{1},x+x_{2}))\left(Y_{\E}(u_{i}(x),x_{2})
Y_{\E}(v_{i}(x),x_{1})\theta (x)\right)w.
\end{eqnarray}
Invoking the cancelation rule (Lemma \ref{lcancellation}) we get
\begin{eqnarray}
& &p(x_{1}+x,x_{2}+x)
\left(Y_{\E}(a(x),x_{1})Y_{\E}(b(x),x_{2})\theta (x)\right)w\nonumber\\
&=&\sum_{i=1}^{r}\iota_{x,x_{2},x_{1}}(q_{i}(x+x_{1},x+x_{2}))
\left(Y_{\E}(u_{i}(x),x_{2})
Y_{\E}(v_{i}(x),x_{1})\theta (x)\right)w,
\end{eqnarray}
proving (\ref{epqsum}).

Now assume that
\begin{eqnarray}
p(x_{1},x_{2})a(x_{1})b(x_{2})
=\sum_{i=1}^{r}p(x_{1},x_{2})
\iota_{x_{2},x_{1}}(f_{i}(x_{2}-x_{1}))u_{i}(x_{2})v_{i}(x_{1}),
\end{eqnarray}
where $f_{i}(x)\in \C((x))$ for $i=1,\dots,r$.
Then
\begin{eqnarray}
& &p(x_{1}+x,x_{2}+x)
\left(Y_{\E}(a(x),x_{1})Y_{\E}(b(x),x_{2})\theta (x)\right)w\nonumber\\
&=&\sum_{i=1}^{r}p(x_{1}+x,x_{2}+x)
\iota_{x_{2},x_{1}}(f_{i}(x_{2}-x_{1}))\left(Y_{\E}(u_{i}(x),x_{2})
Y_{\E}(v_{i}(x),x_{1})\theta (x)\right)w.
\end{eqnarray}
We may assume $p(x,y)\ne 0$.
Write $p(x,y)=(x-y)^{k}p_{1}(x,y)$, where $p_{1}(x,y)\in \C[x,y]$ 
with $p_{1}(x,x)\ne 0$.
Using the cancelation law (Lemma \ref{lcancellation}) we get
\begin{eqnarray}
& &(x_{1}-x_{2})^{k}Y_{\E}(a(x),x_{1})Y_{\E}(b(x),x_{2})\nonumber\\
&=&\sum_{i=1}^{r}(x_{1}-x_{2})^{k}\iota_{x_{2},x_{1}}(f_{i}(x_{2}-x_{1}))
Y_{\E}(u_{i}(x),x_{2})Y_{\E}(v_{i}(x),x_{1}).
\end{eqnarray}
This completes the proof.
\end{proof}

\section{Pseudo-local quasi modules for nonlocal vertex algebras}

In this section we formulate and study notions of 
pseudo-local quasi module and $\S$-quasi module for a nonlocal vertex algebra
with respect to a quantum Yang-Baxter operator $\S$.
We also associate nonlocal vertex algebras to quantum affine algebras and their
restricted modules.

Recall that 
$$\C_{*}(x_{1},x_{2})=\{ p(x_{1},x_{2})/q(x_{1},x_{2})\;|\;
p(x_{1},x_{2})\in \C[[x_{1},x_{2}]],\; 
0\ne q(x_{1},x_{2})\in \C[x_{1},x_{2}]\}.$$
We shall use the following subspaces of $\C_{*}(x_{1},x_{2})$:
\begin{eqnarray*}
\C_{*}(x_{1}/x_{2})&=&\{ f(x_{1}/x_{2})\;|\; f(x)\in \C[[x]][x^{-1}]\}
\subset \C_{*}(x_{1},x_{2}),\\
\C_{*}(x_{1}-x_{2})&=&\{ f(x_{1}-x_{2})\;|\; f(x)\in \C[[x]][x^{-1}]\}
\subset \C_{*}(x_{1},x_{2}).
\end{eqnarray*}

\bd{dpseudo-local-quasi-module}
{\em Let $V$ be a nonlocal vertex algebra.
A {\em pseudo-local quasi $V$-module}
is a quasi $V$-module $(W,Y_{W})$ satisfying the condition that
for any $u,v\in V$, there exist 
$f_{i}(x_{1},x_{2})\in \C_{*}(x_{1},x_{2})$ 
and $u^{(i)}, v^{(i)}\in V$ for $i=1,\dots,n$ such that
\begin{eqnarray}\label{epseudo-local-quasi-module}
p(x_{1},x_{2})Y_{W}(u,x_{1})Y_{W}(v,x_{2})
=\sum_{i=1}^{n}p(x_{1},x_{2})\iota_{x_{2},x_{1}}(f_{i}(x_{1},x_{2}))
Y_{W}(v^{(i)},x_{2})Y_{W}(u^{(i)},x_{1})
\end{eqnarray}
for some nonzero polynomial $p(x_{1},x_{2})$. That is,
a pseudo-local quasi $V$-module
is a quasi $V$-module $(W,Y_{W})$ such that
$\{ Y_{W}(v,x)\;|\; v\in V\}$ is a pseudo-local subspace of $\E(W)$.}
\ed

For a subgroup $\Gamma$ of $\C^{\times}$, we define a {\em pseudo $\Gamma$-local quasi module}
for $V$ to be a quasi $V$-module $(W,Y_{W})$ such that
$\{ Y_{W}(v,x)\;|\; v\in V\}$ is a pseudo $\Gamma$-local subspace of $\E(W)$.

In terms of these notions, from Propositions \ref{pgeneral-locality-key} 
and \ref{pgamma-locality-key}, we immediately have:

\bt{tpseudo-locality-key}
Let $W$ be a vector space and let $S$ be any pseudo-local subset of $\E (W)$.
Then $W$ is a faithful pseudo-local quasi module 
for the nonlocal vertex algebra $\<S\>$ (generated by $S$)
with $Y_{W}(\alpha(x),x_{0})=\alpha(x_{0})$
for $\alpha(x)\in \<S\>$. Furthermore, if $S$ is pseudo $\Gamma$-local for some
subgroup $\Gamma$ of $\C^{\times}$, then $W$ is a pseudo $\Gamma$-local quasi module. 
\et

\br{rYBoperator}
{\em We here recall the definitions of quantum Yang-Baxter operators
(cf. \cite{cp}, \cite{efk}).
Let $U$ be a vector space. A {\em unitary
quantum Yang-Baxter operator} (with two spectral parameters) 
on $U$ is a linear map
$\S (x_{1},x_{2}): U\otimes U\rightarrow U\otimes U\otimes
\C_{*}(x_{1},x_{2})$ such that
\begin{eqnarray}
& &\hspace{2cm} \S_{21}(x_{2},x_{1})\S(x_{1},x_{2})=1,\\
& &\S_{12}(x_{1},x_{2})\S_{13}(x_{1},x_{3})\S_{23}(x_{2},x_{3})
=\S_{23}(x_{2},x_{3})\S_{13}(x_{1},x_{3})\S_{12}(x_{1},x_{2}).
\label{eYBE-def0}
\end{eqnarray}
A {\em unitary trigonometric quantum Yang-Baxter operator} on $U$ 
is to a linear map
$\S: U\otimes U\rightarrow U\otimes U\otimes \C((x))$ such that
\begin{eqnarray}
& &\hspace{2cm} \S_{21}(x)\S(x^{-1})=1,\\
& &\S_{12}(x)\S_{13}(xz)\S_{23}(z)
=\S_{23}(z)\S_{13}(zx)\S_{12}(x)\label{eYBE-def-trig}
\end{eqnarray}
and a {\em unitary rational quantum Yang-Baxter operator} on $U$ 
is to a linear map
$\S: U\otimes U\rightarrow U\otimes U\otimes \C((x))$ such that
\begin{eqnarray}
& &\hspace{2cm} \S_{21}(x)\S(-x)=1,\\
& &\S_{12}(x)\S_{13}(x+z)\S_{23}(z)
=\S_{23}(z)\S_{13}(z+x)\S_{12}(x).\label{eYBE-def-rational}
\end{eqnarray}
It is clear that for any unitary rational quantum 
Yang-Baxter operator $\S(x)$, $\S(x_{1}-x_{2})$ is a unitary quantum 
Yang-Baxter operator with two spectral parameters.
We here use the standard notations in literature.
In particular, we have
$$\S_{21}(x_{1},x_{2})=\sigma_{12}\S_{12}(x_{1},x_{2})\sigma_{12},$$
where $\sigma$ is the flip map
$U\otimes U\rightarrow U\otimes U, u\otimes v\mapsto v\otimes u$.}
\er

\br{rbraid-relation}
{\em For any linear map $S(x_{1},x_{2}): U\otimes U\rightarrow U\otimes
U\otimes \C_{*}(x_{1},x_{2})$, define a linear map $R(x_{1},x_{2})$ of the
same type by $R(x_{1},x_{2})=S(x_{1},x_{2})\sigma$, i.e., 
$R(x_{1},x_{2})(u\otimes v)=S(x_{1},x_{2})(v\otimes u)$. Then
the relation (\ref{eYBE-def0}) is equivalent to the following
braided relation:
\begin{eqnarray}
R_{12}(x_{1},x_{2})R_{23}(x_{1},x_{3})R_{12}(x_{2},x_{3})
=R_{23}(x_{2},x_{3})R_{12}(x_{1},x_{3})R_{23}(x_{1},x_{2}).
\end{eqnarray}}
\er

For any linear map 
$\S: U\otimes U\rightarrow U\otimes U\otimes \C_{*}(x_{1},x_{2})$,
we set
\begin{eqnarray}
\tilde{\S}(x_{1},x_{2})
=(1\otimes 1\otimes \iota_{x_{1},x_{2}})\S(x_{1},x_{2}),
\end{eqnarray}
which is a linear map from $U\otimes U$ to 
$U\otimes U\otimes \C((x_{1}))((x_{2}))$. Note that
$\tilde{\S}(x_{2},x_{1})$ is a linear map from $U\otimes U$ to 
$U\otimes U\otimes \C((x_{2}))((x_{1}))$.

\bd{dSquasimodule}
{\em Let $V$ be a nonlocal vertex algebra and let $\S (x_{1},x_{2})$ 
be a unitary quantum Yang-Baxter operator on $V$.
An {\em $\S$-quasi $V$-module} is a quasi $V$-module $(W,Y_{W})$ satisfying the condition that
for $u,v\in V$, there exists $0\ne p(x_{1},x_{2})\in
\C[x_{1},x_{2}]$
such that
\begin{eqnarray}\label{eSquasimodule}
p(x_{1},x_{2})Y_{W}(u,x_{1})Y_{W}(v,x_{2})
=p(x_{1},x_{2})\sum_{i=1}^{r}\iota_{x_{2},x_{1}}(f_{i}(x_{2},x_{1}))
Y_{W}(v^{(i)},x_{2})Y_{W}(u^{(i)},x_{1}),
\end{eqnarray}
where $\S(x_{2},x_{1})(v\otimes u)=\sum_{i=1}^{r}v^{(i)}\otimes u^{(i)}\otimes f_{i}(x_{2},x_{1})$.}
\ed

\br{rfact-sim-q}
{\em Let $W$ be a vector space. For
$$F(x_{1},\dots,x_{n}),G(x_{1},\dots,x_{n})
\in (\End W)[[x_{1}^{\pm 1},\dots,x_{n}^{\pm 1}]],$$
we write $F\sim G$ if 
there exists a nonzero polynomial $p(x_{1},\dots,x_{n})$ such that
\begin{eqnarray}
p(x_{1},\dots,x_{n})F(x_{1},\dots,x_{n})
=p(x_{1},\dots,x_{n})G(x_{1},\dots,x_{n}).
\end{eqnarray}
Clearly, the relation ``$\sim$'' is an equivalent relation
on $(\End W)[[x_{1}^{\pm 1},\dots,x_{n}^{\pm 1}]]$. Furthermore, if
$$F(x_{1},\dots,x_{n}),G(x_{1},\dots,x_{n})
\in \Hom (W,W((x_{1}))\cdots ((x_{n}))),$$
then $F\sim G$ if and only if $F=G$.
Furthermore, all the partial derivative operators $\partial/\partial x_{i}$ 
preserve the equivalence relation and
the left multiplication by a Laurent polynomial also preserves the
equivalence relation.}
\er

\bd{dpinw}
{\em Let $V$ be a nonlocal vertex algebra and let $(W,Y_{W})$ be a quasi $V$-module. 
For any positive integer $n$ we define a linear map 
$$\pi^{W}_{n}: V^{\otimes n}\otimes \C((x_{1}))\cdots ((x_{n}))
\rightarrow \Hom (W,W((x_{1}))\cdots ((x_{n})))$$
by}
$$\pi_{n}^{W}(v_{1}\otimes \cdots \otimes v_{n}\otimes f)(w)
=f(x_{1},\dots,x_{n})Y_{W}(v_{1},x_{1})\cdots Y_{W}(v_{n},x_{n})w.$$
\ed

The following theorem was inspired by a result of \cite{ek}
(Proposition 1.11.):

\bt{tfirst}
Let $V$ be a nonlocal vertex algebra
and let $(W,Y_{W})$ be a pseudo-local quasi $V$-module.
Assume that the linear map $\pi_{3}^{W}$ is injective.
Then the pseudo-locality (\ref{epseudo-local-quasi-module}) 
uniquely defines a unitary quantum Yang-Baxter operator 
$\S(x_{1},x_{2}):
V\otimes V\rightarrow V\otimes V\otimes \C_{*}(x_{1},x_{2})$ 
with $\S(x_{2},x_{1})(v\otimes u)=\sum_{i=1}^{r}v^{(i)}\otimes
u^{(i)}\otimes f_{i}(x_{2},x_{1})$
such that
$W$ is an $\S$-quasi $V$-module and
$$\hspace{1cm} [1\otimes
\D,\S(x_{1},x_{2})]=-\frac{\partial}{\partial x_{1}}\S(x_{1},x_{2}).$$
Furthermore, if $\{Y_{W}(v,x)\;|\; v\in V\}$ is $\S$-local,
then $\S(x_{1},x_{2})$ is given 
by a unitary rational quantum Yang-Baxter operator.
\et

\begin{proof} With $Y_{W}({\bf 1},x)=1_{W}$, 
it is clear that $\pi_{2}^{W}$ is injective.
For $u,v\in V$, by assumption
there exists $\Psi(x_{1},x_{2})\in V\otimes V\otimes
\C_{*}(x_{1},x_{2})$ such that
\begin{eqnarray}\label{eslocality-prop}
Y_{W}(u,x_{1})Y_{W}(v,x_{2})
\sim Y_{W}(x_{2})(1\otimes Y_{W}(x_{1}))\iota_{x_{2},x_{1}}\Psi(x_{2},x_{1}).
\end{eqnarray}
Suppose that we also have 
\begin{eqnarray}
Y_{W}(u,x_{1})Y_{W}(v,x_{2})
\sim Y_{W}(x_{2})(1\otimes Y_{W}(x_{1}))\iota_{x_{2},x_{1}}\Phi(x_{2},x_{1})
\end{eqnarray}
for $\Phi(x_{1},x_{2})\in V\otimes V\otimes \C_{*}(x_{1},x_{2})$.
Then
\begin{eqnarray*}
Y_{W}(x_{2})(1\otimes Y_{W}(x_{1}))\iota_{x_{2},x_{1}}\Psi(x_{2},x_{1})
\sim Y_{W}(x_{2})(1\otimes Y_{W}(x_{1}))\iota_{x_{2},x_{1}}\Phi (x_{2},x_{1}).
\end{eqnarray*}
In view of Remark \ref{rfact-sim-q}, we have
\begin{eqnarray*}
Y_{W}(x_{2})(1\otimes Y_{W}(x_{1}))\iota_{x_{2},x_{1}}\Psi (x_{2},x_{1})
=Y_{W}(x_{2})(1\otimes Y_{W}(x_{1}))\iota_{x_{2},x_{1}}\Phi(x_{2},x_{1}).
\end{eqnarray*}
As $\pi_{2}^{W}$ is injective,  we have $\Psi(x_{2},x_{1})=\Phi(x_{2},x_{1})$.
Thus, we have a well defined linear map 
$\S(x_{1},x_{2}): V\otimes V\rightarrow V\otimes
V\otimes \C_{*}(x_{1},x_{2})$ such that for $u,v\in V$,
\begin{eqnarray}
Y_{W}(u,x_{1})Y_{W}(v,x_{2})\sim
Y_{W}(x_{2})(1\otimes Y_{W}(x_{1}))\tilde{\S}(x_{2},x_{1})(v\otimes u).
\end{eqnarray}

Let $u,v\in V,\; w\in W$ and let $p(x_{1},x_{2})$ be a nonzero
polynomial such that
$$p(x_{1},x_{2})\S(x_{1},x_{2})(v\otimes u)
=p(x_{1},x_{2})\tilde{\S}(x_{1},x_{2})(v\otimes u)
\in V\otimes V\otimes \C[[x_{1},x_{2}]]$$
and such that
$$p(x_{1},x_{2})Y_{W}(x_{2})(1\otimes Y_{W}(x_{1}))(u\otimes v\otimes w)
=p(x_{1},x_{2})Y_{W}(x_{1})(1\otimes Y_{W}(x_{2}))
(\tilde{\S}(x_{1},x_{2})(v\otimes u)\otimes w). $$
Let $q(x_{1},x_{2})$ be another nonzero polynomial such that
\begin{eqnarray*}
& &q(x_{1},x_{2})Y_{W}(x_{1})(1\otimes Y_{W}(x_{2}))
(p(x_{1},x_{2})\S(x_{1},x_{2})(v\otimes u)\otimes w)\nonumber\\
&=& q(x_{1},x_{2})Y_{W}(x_{2})(1\otimes Y_{W}(x_{1}))
\tilde{\S}^{12}(x_{2},x_{1})\sigma_{12}
\left(p(x_{1},x_{2})\S(x_{1},x_{2})(v\otimes u)\otimes w\right).
\end{eqnarray*}
(Note that $\tilde{\S}(x_{2},x_{1})\sigma
\tilde{\S}(x_{1},x_{2})(v\otimes u)$ 
does not exist in general.) Then
\begin{eqnarray*}
& &p(x_{1},x_{2})q(x_{1},x_{2})Y_{W}(x_{2})(1\otimes Y_{W}(x_{1}))(u\otimes v\otimes w)
\nonumber\\
&=&p(x_{1},x_{2})q(x_{1},x_{2})Y_{W}(x_{1})(1\otimes Y_{W}(x_{2}))
(\tilde{\S}(x_{1},x_{2})(v\otimes u)\otimes w)
\nonumber\\
&=&q(x_{1},x_{2})Y_{W}(x_{1})(1\otimes Y_{W}(x_{2}))
\left(p(x_{1},x_{2})\S(x_{1},x_{2})(v\otimes u)\otimes w)\right)
\nonumber\\
&=&q(x_{1},x_{2})Y_{W}(x_{2})(1\otimes Y_{W}(x_{1}))
\tilde{\S}^{12}(x_{2},x_{1})\sigma_{12}
\left(p(x_{1},x_{2})\S(x_{1},x_{2})(v\otimes u)\otimes w\right).
\end{eqnarray*}
As $\pi_{2}^{W}$ is injective, we get
\begin{eqnarray*}
q(x_{1},x_{2})p(x_{1},x_{2})(u\otimes v)
=q(x_{1},x_{2})\tilde{\S}(x_{2},x_{1})\sigma
\left(p(x_{1},x_{2})\S(x_{1},x_{2})(v\otimes
u)\right),
\end{eqnarray*}
{}from which we get
\begin{eqnarray}
\S_{21}(x_{2},x_{1})\S(x_{1},x_{2})(v\otimes u)=(v\otimes u).
\end{eqnarray}
This proves that $\S_{21}(x_{2},x_{1})\S(x_{1},x_{2})=1$.

As in Remark \ref{rbraid-relation}, 
set $R(x_{1},x_{2})=\S(x_{1},x_{2})\sigma$.
Let $a,b,c,d\in V$. Using the pseudo locality we write the product 
$Y_{W}(a,x_{1})Y_{W}(b,x_{2})Y_{W}(c,x_{3})$
in terms of the reversely ordered product 
$Y_{W}(c,x_{3})Y_{W}(b,x_{2})Y_{W}(a,x_{1})$ in two different ways
as follows:
\begin{eqnarray}
& &Y_{W}(x_{1})(1\otimes Y_{W}(x_{2})(1\otimes 1\otimes Y_{W}(x_{3}))
(a\otimes b\otimes c\otimes d)\nonumber\\
&\sim& Y_{W}(x_{1})(1\otimes Y_{W}(x_{3}))(1\otimes 1\otimes
Y_{W}(x_{2})) \tilde{R}_{23}(x_{3},x_{2})
(a\otimes b\otimes c\otimes d)\nonumber\\
&\sim& Y_{W}(x_{3})(1\otimes Y_{W}(x_{1}))(1\otimes 1\otimes Y_{W}(x_{2}))
\tilde{R}_{12}(x_{3},x_{1})\tilde{R}_{23}(x_{3},x_{2})
(a\otimes b\otimes c\otimes d)\nonumber\\
&\sim& Y_{W}(x_{3})(1\otimes Y_{W}(x_{2}))(1\otimes 1\otimes Y_{W}(x_{1}))
\tilde{R}_{23}(x_{2},x_{1})\nonumber\\
& &\ \ \cdot \tilde{R}_{12}(x_{3},x_{1})
\tilde{R}_{23}(x_{3},x_{2})
(a\otimes b\otimes c\otimes d)
\end{eqnarray}
and
\begin{eqnarray}
& &Y_{W}(x_{1})(1\otimes Y_{W}(x_{2})(1\otimes 1\otimes Y_{W}(x_{3}))
(a\otimes b\otimes c\otimes d)\nonumber\\
&\sim& Y_{W}(x_{2})(1\otimes Y_{W}(x_{1}))(1\otimes 1\otimes Y_{W}(x_{3}))
\tilde{R}_{12}(x_{2},x_{1})
(a\otimes b\otimes c\otimes d)\nonumber\\
&\sim& Y_{W}(x_{2})(1\otimes Y_{W}(x_{3}))(1\otimes 1\otimes Y_{W}(x_{1}))
\tilde{R}_{23}(x_{3},x_{1})\tilde{R}_{12}(x_{2},x_{1})
(a\otimes b\otimes c\otimes d)\nonumber\\
&\sim& Y_{W}(x_{3})(1\otimes Y_{W}(x_{2}))(1\otimes 1\otimes Y_{W}(x_{1}))
\tilde{R}_{12}(x_{3},x_{2})\nonumber\\
& &\ \ \ \ \cdot \tilde{R}_{23}(x_{3},x_{1})
\tilde{R}_{12}(x_{2},x_{1})
(a\otimes b\otimes c\otimes d).
\end{eqnarray}
Thus
\begin{eqnarray}
 & &Y_{W}(x_{3})(1\otimes Y_{W}(x_{2}))(1\otimes 1\otimes Y_{W}(x_{1}))
\tilde{R}_{23}(x_{2},x_{1})\nonumber\\
& &\ \ \ \ \cdot \tilde{R}_{12}(x_{3},x_{1})\tilde{R}_{23}(x_{3},x_{2})
(a\otimes b\otimes c\otimes d)\nonumber\\
&\sim&Y_{W}(x_{3})(1\otimes Y_{W}(x_{2}))(1\otimes 1\otimes Y_{W}(x_{1}))
\tilde{R}_{12}(x_{3},x_{2})\nonumber\\
& &\ \ \ \ \cdot \tilde{R}_{23}(x_{3},x_{1})\tilde{R}_{12}(x_{2},x_{1})
(a\otimes b\otimes c\otimes d).
\end{eqnarray}
In view of Remark \ref{rfact-sim-q} we have
\begin{eqnarray}
 & &Y_{W}(x_{3})(1\otimes Y_{W}(x_{2}))(1\otimes 1\otimes Y_{W}(x_{1}))
\tilde{R}_{23}(x_{2},x_{1})\nonumber\\
& &\ \ \ \ \cdot \tilde{R}_{12}(x_{3},x_{1})\tilde{R}_{23}(x_{3},x_{2})
(a\otimes b\otimes c\otimes d)\nonumber\\
&=&Y_{W}(x_{3})(1\otimes Y_{W}(x_{2}))(1\otimes 1\otimes Y_{W}(x_{1}))
\tilde{R}_{12}(x_{3},x_{2})\nonumber\\
& &\ \ \ \ \cdot \tilde{R}_{23}(x_{3},x_{1})\tilde{R}_{12}(x_{2},x_{1})
(a\otimes b\otimes c\otimes d).
\end{eqnarray}
Since $\pi_{3}^{W}$ is injective by assumption, we have
\begin{eqnarray}
\tilde{R}_{23}(x_{2},x_{1})\tilde{R}_{12}(x_{3},x_{1})
\tilde{R}_{23}(x_{3},x_{2})
=\tilde{R}_{12}(x_{3},x_{2})\tilde{R}_{23}(x_{3},x_{1})
\tilde{R}_{12}(x_{2},x_{1}).
\end{eqnarray}
Thus
\begin{eqnarray}\label{eSTrelation}
R_{23}(x_{2},x_{1})R_{12}(x_{3},x_{1})R_{23}(x_{3},x_{2})
=R_{12}(x_{3},x_{2})R_{23}(x_{3},x_{1})R_{12}(x_{2},x_{1}).
\end{eqnarray}
{}From Remark \ref{rbraid-relation} we have
\begin{eqnarray}
\S_{23}(x_{2},x_{1})\S_{13}(x_{3},x_{1})\S_{12}(x_{3},x_{2})
=\S_{12}(x_{3},x_{2})\S_{13}(x_{3},x_{1})\S_{23}(x_{2},x_{1}).
\end{eqnarray}
This proves that $\S(x_{1},x_{2})$ is a unitary quantum Yang-Baxter operator.

For $u,v\in V$, by definition we have
\begin{eqnarray}
Y_{W}(x_{1})(1\otimes Y_{W}(x_{2}))(\D u\otimes v\otimes w)
\sim Y_{W}(x_{2})(1\otimes Y_{W}(x_{1}))(\tilde{\S}(x_{2},x_{1})
(v\otimes \D u)\otimes w).
\end{eqnarray}
On the other hand, using Lemma \ref{lquasi-dproperty} we have
\begin{eqnarray}
& &Y_{W}(x_{1})(1\otimes Y_{W}(x_{2}))(\D u\otimes v\otimes w)\nonumber\\
&=&\frac{\partial}{\partial x_{1}}Y_{W}(x_{1})(1\otimes Y_{W}(x_{2}))
(u\otimes v\otimes w)\nonumber\\
&\sim& \frac{\partial}{\partial x_{1}}
\left( Y_{W}(x_{2})(1\otimes Y_{W}(x_{1}))(\tilde{\S}(x_{2},x_{1})
(v\otimes u)\otimes w)\right)\nonumber\\
&=&\left(\frac{\partial}{\partial x_{1}} 
Y_{W}(x_{2})(1\otimes Y_{W}(x_{1}))\right)
(\tilde{\S}(x_{2},x_{1})(v\otimes u)\otimes w)\nonumber\\
& &+Y(x_{2})(1\otimes Y(x_{1}))
\frac{\partial}{\partial x_{1}}(\tilde{\S}(x_{2},x_{1})(v\otimes u)\otimes w)
\nonumber\\
&=&Y_{W}(x_{2})(1\otimes Y_{W}(x_{1}))(1\otimes \D\otimes 1)
(\tilde{\S}(x_{2},x_{1})(v\otimes u)\otimes w)\nonumber\\
& &+Y_{W}(x_{2})(1\otimes Y_{W}(x_{1}))
\frac{\partial}{\partial x_{1}}(\tilde{\S}(x_{2},x_{1})(v\otimes u)\otimes w).
\end{eqnarray}
Thus
\begin{eqnarray}
& &Y_{W}(x_{2})(1\otimes Y_{W}(x_{1})(\tilde{\S}(x_{2},x_{1})(v\otimes
\D u)\otimes w)\nonumber\\
&\sim&Y_{W}(x_{2})(1\otimes Y_{W}(x_{1}))(1\otimes \D \otimes 1)
(\tilde{\S}(x_{2},x_{1})(v\otimes u)\otimes w)\nonumber\\
& &+Y_{W}(x_{2})(1\otimes Y_{W}(x_{1}))
\frac{\partial}{\partial x_{1}}(\tilde{\S}(x_{2},x_{1})(v\otimes u)\otimes w).
\end{eqnarray}
In view of Remark \ref{rfact-sim-q} the above equivalent quantities
on both sides are actually equal. Again, as $\pi_{2}^{W}$ is
injective, we have
\begin{eqnarray*}
\tilde{\S}(x_{2},x_{1})(v\otimes \D u)
=(1\otimes \D)\tilde{\S}(x_{2},x_{1})(v\otimes u)+
\frac{\partial}{\partial x_{1}}\tilde{\S}(x_{2},x_{1})(v\otimes u).
\end{eqnarray*}
This proves
\begin{eqnarray}
[1\otimes \D, \S(x_{2},x_{1})]
=-\frac{\partial}{\partial x_{1}}\S(x_{2},x_{1}).
\end{eqnarray}
The last assertion is clear. Now the proof is complete.
\end{proof}

Next we associate nonlocal vertex algebras to quantum affine algebras. 
First, we follow \cite{drinfeld} and \cite{fj} to present the quantum affine algebras.
Let $A=(A_{ij})_{i,j=1}^{l}$ be a Cartan matrix 
of type $A_{l}\; (l\ge 1), D_{l} \;(l\ge 4)$, or $E_{l}$ $(l=6,7,8)$.
Let $q$ be a nonzero complex number.
For $1\le i,j\le l$, set
\begin{eqnarray}
f_{ij}(x)=(q^{A_{ij}}x-1)/(x-q^{A_{ij}})\in \C(x).
\end{eqnarray}
Then we set
\begin{eqnarray}
g_{ij}(x)^{\pm 1}=\iota_{x,0} f_{ij}(x)^{\pm 1}\in \C[[x]],
\end{eqnarray}
where $\iota_{x,0} f_{ij}(x)^{\pm 1}$ are the formal Taylor series expansions of
$f_{ij}(x)^{\pm 1}$ at $0$.
The quantum affine algebra $U_{q}(\hat{\g})$ is (isomorphic to) the associative
algebra with unit $1$ and the generators
\begin{eqnarray}
\{ X_{ik}^{\pm},\; \phi_{im},\; \psi_{in}, \gamma^{1/2},\gamma^{-1/2}
\;|\; i=1,\dots,l; \;k\in \Z,\; m\in -\Z_{+},\;n\in \Z_{+}\},
\end{eqnarray}
where $\gamma^{\pm 1/2}$ are central,
satisfying the relations below written in terms of the following
generating functions in a formal variable $z$:
\begin{eqnarray}
X_{i}^{\pm}(z)=\sum_{k\in \Z}X_{ik}^{\pm}z^{-k},\ \ \ \ 
\phi_{i}(z)=\sum_{m\in -\Z_{+}}\phi_{im}z^{-m},\ \ \ \ 
\psi_{i}(z)=\sum_{n\in \Z_{+}}\psi_{in}z^{-n}.
\end{eqnarray}
The relations are
\begin{eqnarray}
& &\gamma^{1/2}\gamma^{-1/2}=\gamma^{-1/2}\gamma^{1/2}=1,\\
& &\phi_{i0}\psi_{i0}=\psi_{i0}\phi_{i0}=1,\\
& &[\phi_{i}(z),\phi_{j}(w)]=0,\ \ \ \ 
[\psi_{i}(z),\psi_{j}(w)]=0,\\
& &\phi_{i}(z)\psi_{j}(w)\phi_{i}(z)^{-1}\psi_{j}(w)^{-1}
=g_{ij}(zw^{-1}\gamma^{-1})/g_{ij}(zw^{-1}\gamma),\\
& &\phi_{i}(z)X^{\pm}_{j}(w)\phi_{i}(z)^{-1}
=g_{ij}(zw^{-1}\gamma^{\mp 1/2})^{\pm 1}X^{\pm}_{j}(w),\\
& &\psi_{i}(z)X^{\pm}_{j}(w)\psi_{i}(z)^{-1}
=g_{ij}(z^{-1}w\gamma^{\mp 1/2})^{\mp 1}X^{\pm}_{j}(w),\\
& &(z-q^{\pm 4A_{ij}}w)X^{\pm}_{i}(z)X^{\pm}_{j}(w)
=(q^{\pm 4A_{ij}}z-w)X^{\pm}_{j}(w)X^{\pm}_{i}(z),\\
& &[X^{+}_{i}(z),X^{-}_{j}(w)]=\delta_{ij}
\left(\delta(zw^{-1}\gamma^{-1})\psi_{i}(w\gamma^{1/2})
-\delta(zw^{-1}\gamma)\phi_{i}(z\gamma^{1/2})\right)/(q-q^{-1})
\ \ \ \ 
\end{eqnarray}
and there is one more set of relations of Serre type.

{}From the above relations we get
\begin{eqnarray}
& &\psi_{j}(w)\phi_{i}(z)
=g_{ij}(zw^{-1}\gamma)/g_{ij}(zw^{-1}\gamma^{-1}) \phi_{i}(z)\psi_{j}(w),\\
& &\psi_{i}(z)X^{\pm}_{j}(w)
=g_{ij}(z^{-1}w\gamma^{\mp 1/2})^{\mp 1}X^{\pm}_{j}(w)\psi_{i}(z),\\
& &X^{\pm}_{j}(w)\psi_{i}(z)=g_{ij}(z^{-1}w\gamma^{\mp 1/2})^{\pm 1}\psi_{i}(z)X^{\pm}_{j}(w),\\
& &\psi_{i}(z)X^{\pm}_{j}(w)
=g_{ij}(z^{-1}w\gamma^{\mp 1/2})^{\mp 1}X^{\pm}_{j}(w)\psi_{i}(z),\\
& &X^{\pm}_{j}(w)\psi_{i}(z)=g_{ij}(z^{-1}w\gamma^{\mp 1/2})^{\pm 1}\psi_{i}(z)X^{\pm}_{j}(w),\\
& &(z-w\gamma)(z\gamma-w)[X^{+}_{i}(z),X^{-}_{j}(w)]=0.
\end{eqnarray}

Let $W$ be a {\em restricted} $U_{q}(\hat{\g})$-module in the sense that
$X_{i}^{\pm}(x), \phi_{i}(x),\psi_{i}(x)$ acting on $W$ are
elements of $\E(W)$. This amounts to that for any $w\in W$,
$X_{ik}w=0$ and $\psi_{ik}w=0$ for $i=1,\dots,l$ and for $k$ sufficiently large.
Furthermore, assume that $W$ is of {\em level $\ell$} in the sense that
$\gamma^{\pm 1/2}$ act on $W$ as scalars $q^{\pm \ell/4}$.
(Rigorously speaking, one needs to choose a branch of $\log q$.)
Set
\begin{eqnarray}
\Gamma (q,\ell)=\{ q^{m+n\ell/4}\;|\; m,n\in \Z\}\subset \C^{\times}.
\end{eqnarray}
We see that $S_{W}=\{ \phi_{i}(x), \psi_{i}(x), X^{\pm}(x)
\;|\; i=1,\dots,l\}$ is a pseudo $\Gamma(q,\ell)$-local subset of $\E(W)$.
In view of Theorem \ref{tgeneratingthem}, there exists a unique smallest closed 
subspace $V_{W}$ of $\E(W)$ that contains $S_{W}\cup \{ 1_{W}\}$
and is stable under the actions of $R_{\alpha}$ for $\alpha\in \Gamma(q,\ell)$.
Furthermore, $V_{W}$ is a nonlocal vertex algebra generated by $R_{\Gamma(q,\ell)}(S_{W})$
with $W$ as a quasi module. By Theorem \ref{tpseudo-locality-key},
$W$ is a pseudo  $\Gamma (q,\ell)$-local quasi module for $V_{W}$. To summarize we have:

\bp{pqaffine-algebra}
Let $q$ be a nonzero complex number and let $W$ be any restricted module of level 
$\ell\in \C$ for the quantum
affine algebra $U_{q}(\hat{\g})$. Set
$$S_{W}=\{ 1_{W}, \phi_{i}(x), \psi_{i}(x), X^{\pm}(x)
\;|\; i=1,\dots,l\}.$$ 
Then $S_{W}$ is a pseudo $\Gamma(q,\ell)$-local subset of $\E(W)$
and there exists a unique smallest closed 
pseudo $\Gamma (q,\ell)$-local subspace $V_{W}$ of $\E(W)$
that contains $S_{W}\cup \{ 1_{W}\}$ and is stable under the actions of
$R_{\alpha}$ for $\alpha\in \Gamma(q,\ell)$.
Furthermore, $V_{W}$ is a nonlocal vertex algebra
with $W$ as a faithful pseudo $\Gamma(q,\ell)$-local quasi module.
\ep

We conjecture that if $W$ is the level-one highest weight $U_{q}(\hat{\g})$-module 
constructed by Frenkel-Jing in \cite{fj}, then the adjoint $V_{W}$-module is irreducible and
$W$ is an $\S$-quasi module with respect to a uniquely determined
quantum Yang-Baxter operator $\S(x_{1},x_{2})$ on $V_{W}$.

\section{Quantum vertex algebras and their modules}

In this section we formulate and study notions of 
weak quantum vertex algebra and quantum vertex algebra.
For our need we also prove certain basic results
for general nonlocal vertex algebras and modules.
Our notion of quantum vertex algebra was partly motivated 
by Etingof-Kazhdan's notion of quantum vertex operator algebra in \cite{ek}
and we also use some of their important ideas.

First, using the same argument as in \cite{ll} (cf. \cite{fhl}) we have:

\bl{lquasi-property}
Let $V$ be a nonlocal vertex algebra and 
let $(W,Y_{W})$ be a quasi $V$-module. 
Then for any $u,v\in V,\; w\in
W,\; w^{*}\in W^{*}$, the formal series
$$\<w^{*},Y_{W}(u,x_{1})Y_{W}(v,x_{2})w\>,$$
an element of $\C((x_{1}))((x_{2}))$, lies in the range of
$\iota_{x_{1},x_{2}}$, and the formal series
$$\<w^{*},Y_{W}(Y(u,x_{0})v,x_{2})w\>,$$
an element of $\C((x_{2}))((x_{0}))$, lies in the range of
$\iota_{x_{2},x_{0}}$.
Furthermore,
\begin{eqnarray}
\left(\iota_{x_{1},x_{2}}^{-1}\<w^{*},Y_{W}(u,x_{1})Y_{W}(v,x_{2})w\>\right)|_{x_{1}=x_{0}+x_{2}}
=\iota_{x_{2},x_{0}}^{-1}\<w^{*},Y_{W}(Y(u,x_{0})v,x_{2})w\>.
\end{eqnarray}
If $(W,Y_{W})$ is a $V$-module, then
\begin{eqnarray}
\iota_{x_{1},x_{2}}^{-1}\<w^{*},Y_{W}(u,x_{1})Y_{W}(v,x_{2})w\>
\in \C[[x_{1},x_{2}]][x_{1}^{-1},x_{2}^{-1},(x_{1}-x_{2})^{-1}].
\end{eqnarray}
\el

We also have (\cite{ll}, Theorem 3.6.3):

\bp{pskew-slocality}
Let $V$ be a nonlocal vertex algebra, let $(W,Y_{W})$ be a $V$-module,
and let $$u,v,u^{(1)},v^{(1)},\dots,u^{(r)},v^{(r)}\in V,\; 
f_{1}(x),\dots,f_{r}(x)\in \C((x)).$$
 If
\begin{eqnarray}
Y(u,x)v
=f_{1}(-x)e^{x\D}Y(v^{(1)},-x)u^{(1)}+\cdots 
+f_{r}(-x)e^{x\D}Y(v^{(r)},-x)u^{(r)}
\label{e-Sskew-relation}
\end{eqnarray}
(the $\S$-skew symmetry), 
then 
\begin{eqnarray}\label{e-Slocality-relation}
(x_{1}-x_{2})^{k}Y_{W}(u,x_{1})Y_{W}(v,x_{2})
=(x_{1}-x_{2})^{k}\sum_{i=1}^{r}
f_{i}(x_{2}-x_{1})Y_{W}(v^{(i)},x_{2})Y_{W}(u^{(i)},x_{1})
\end{eqnarray}
for some nonnegative integer $k$.
If $(W,Y_{W})$ is faithful, the converse is also true.
\ep

\begin{proof} 
For the first assertion,
following the proof of Theorem 3.6.3 of \cite{ll},
using Lemma \ref{lquasi-property}, the $\S$-skew symmetry (\ref{e-Sskew-relation})
and the $\D$-derivative property (Lemma \ref{lquasi-dproperty}), we get
\begin{eqnarray}\label{e-sjacobi-module-1}
& &x_{0}^{-1}\delta\left(\frac{x_{1}-x_{2}}{x_{0}}\right)Y_{W}(u,x_{1})Y_{W}(v,x_{2})\nonumber\\
& &\ \ \ \ -x_{0}^{-1}\delta\left(\frac{x_{2}-x_{1}}{-x_{0}}\right)
\sum_{i=1}^{r}f_{i}(x_{2}-x_{1})Y_{W}(v^{(i)},x_{2})Y_{W}(u^{(i)},x_{1})\nonumber\\
&=&x_{2}^{-1}\delta\left(\frac{x_{1}-x_{0}}{x_{2}}\right)Y_{W}(Y(u,x_{0})v,x_{2}),
\end{eqnarray}
{}from which we get (\ref{e-Slocality-relation}).

Assume that $(W,Y_{W})$ is faithful and (\ref{e-Slocality-relation}) holds. 
In view of Lemma \ref{lformaljacobiidentity},
(\ref{e-sjacobi-module-1}) holds and we have
\begin{eqnarray*}
& &-x_{0}^{-1}\delta\left(\frac{x_{2}-x_{1}}{-x_{0}}\right)
\sum_{i=1}^{r}f_{i}(x_{2}-x_{1})Y_{W}(v^{(i)},x_{2})Y_{W}(u^{(i)},x_{1})\nonumber\\
& &\ \ \ \ -(-x_{0})^{-1}\delta\left(\frac{x_{1}-x_{2}}{x_{0}}\right)
Y_{W}(u,x_{1})Y_{W}(v,x_{2})\nonumber\\
&=&x_{1}^{-1}\delta\left(\frac{x_{2}+x_{0}}{x_{1}}\right)
\sum_{i=1}^{r}f_{i}(-x_{0})Y_{W}(Y(v^{(i)},-x_{0})u^{(i)},x_{1})\nonumber\\
&=&x_{1}^{-1}\delta\left(\frac{x_{2}+x_{0}}{x_{1}}\right)
\sum_{i=1}^{r}Y_{W}(f_{i}(-x_{0})Y(v^{(i)},-x_{0})u^{(i)},x_{2}+x_{0})\nonumber\\
&=&x_{1}^{-1}\delta\left(\frac{x_{2}+x_{0}}{x_{1}}\right)
\sum_{i=1}^{r}Y_{W}(f_{i}(-x_{0})e^{x_{0}\D}Y(v^{(i)},-x_{0})u^{(i)},x_{2}),
\end{eqnarray*}
where we are also using the $\D$-derivative property 
(Lemma \ref{lquasi-dproperty}).
Consequently we have
$$Y_{W}(Y(u,x_{0})v,x_{2})
=\sum_{i=1}^{r}Y_{W}(f_{i}(-x_{0})e^{x_{0}\D}
Y(v^{(i)},-x_{0})u^{(i)},x_{2}).$$
As $Y_{W}$ is injective we obtain (\ref{e-Sskew-relation}).
\end{proof}

With the adjoint module $(V,Y)$ being faithful, 
as immediate consequences of Proposition \ref{pskew-slocality} we have:

\bc{cSlocality-Sskew-algebra}
Let $V$ be a nonlocal vertex algebra and let
$$u,v,u^{(1)},v^{(1)},\dots,u^{(r)},v^{(r)}\in V,\; 
f_{1}(x),\dots,f_{r}(x)\in \C((x)).$$
Then
\begin{eqnarray}
Y(u,x)v=f_{1}(-x)e^{x\D}Y(v^{(1)},-x)u^{(1)}+\cdots 
+f_{r}(-x)e^{x\D}Y(v^{(r)},-x)u^{(r)}
\end{eqnarray}
if and only if
\begin{eqnarray}
(x_{1}-x_{2})^{k}Y(u,x_{1})Y(v,x_{2})
=(x_{1}-x_{2})^{k}\sum_{i=1}^{r}
f_{i}(x_{2}-x_{1})Y(v^{(i)},x_{2})Y(u^{(i)},x_{1})
\end{eqnarray}
for some nonnegative integer $k$.
\ec

\bc{cSlocality-algebra-module}
Let $V$ be a nonlocal vertex algebra, let $(W,Y_{W})$ be a $V$-module and let
$$u,v,u^{(1)},v^{(1)},\dots,u^{(r)},v^{(r)}\in V,\; 
f_{1}(x),\dots,f_{r}(x)\in \C((x)).$$
If
\begin{eqnarray}
(x_{1}-x_{2})^{k}Y(u,x_{1})Y(v,x_{2})
=(x_{1}-x_{2})^{k}\sum_{i=1}^{r}
f_{i}(x_{2}-x_{1})Y(v^{(i)},x_{2})Y(u^{(i)},x_{1})
\end{eqnarray}
for some nonnegative integer $k$, then 
\begin{eqnarray}\label{eSlocality-module-cor}
(x_{1}-x_{2})^{l}Y_{W}(u,x_{1})Y_{W}(v,x_{2})
=(x_{1}-x_{2})^{l}\sum_{i=1}^{r}
f_{i}(x_{2}-x_{1})Y_{W}(v^{(i)},x_{2})Y_{W}(u^{(i)},x_{1})
\end{eqnarray}
for some nonnegative integer $l$. If $(W,Y_{W})$ is faithful, the converse is also true.
\ec

Now we are ready to introduce the notions of weak quantum vertex algebra and 
quantum vertex algebra and present the main results.

\bd{dweak-qva}
{\em A {\em weak quantum vertex algebra} (over $\C$)
is a vector space $V$ equipped with a linear map
$$Y: V\rightarrow \E(V)=\Hom (V,V((x)));\;\; v\mapsto Y(v,x)$$
and a distinguished vector ${\bf 1}\in V$ such that
\begin{eqnarray}
& &Y({\bf 1},x)=1,\\
& &Y(v,x){\bf 1}\in V[[x]]\;\;\mbox{ and }
\; \lim_{x\rightarrow 0}Y(v,x){\bf 1}=v\;\;\;\mbox{ for }v\in V,
\end{eqnarray}
and such that for $u,v\in V$,
there exists 
$$\sum_{i=1}^{r}v^{(i)}\otimes u^{(i)}\otimes f_{i}(x)\in
V\otimes V\otimes \C((x))$$
such that
\begin{eqnarray}\label{e-Rjacobi}
& &x_{0}^{-1}\delta\left(\frac{x_{1}-x_{2}}{x_{0}}\right)
Y(u,x_{1})Y(v,x_{2})\nonumber\\
& &\ \ \ \ -x_{0}^{-1}\delta\left(\frac{x_{2}-x_{1}}{-x_{0}}\right)
\sum_{i=1}^{r} f_{i}(-x_{0})Y(v^{(i)},x_{2})Y(u^{(i)},x_{1})\nonumber\\
&=&x_{2}^{-1}\delta\left(\frac{x_{1}-x_{0}}{x_{2}}\right)
Y(Y(u,x_{0})v,x_{2}).
\end{eqnarray}}
\ed

Clearly, the notion of quantum vertex algebra generalizes the notions of 
vertex algebra and vertex superalgebra.

\br{rcomments}
{\em In view of Lemma \ref{lformaljacobiidentity}, we see that
a weak quantum vertex algebra amounts to a
nonlocal vertex algebra $V$ satisfying the following
{\em $\S$-locality}: For any $u,v\in V$, there exist
$f_{i}(x)\in \C((x)),\; u^{(i)},v^{(i)}\in V$, $i=1,\dots,r$ such that
\begin{eqnarray}\label{eslocality}
(x_{1}-x_{2})^{k}Y(u,x_{1})Y(v,x_{2})
=(x_{1}-x_{2})^{k}\sum_{i=1}^{r}f_{i}(x_{2}-x_{1})
Y(v^{(i)},x_{2})Y(u^{(i)},x_{1})
\end{eqnarray}
for some nonnegative integer $k$.}
\er

For a weak quantum vertex algebra $V$, we define a {\em $V$-module}
 to be a module for $V$ viewed as a nonlocal vertex algebra.

\bl{lwqva-module-def}
Let $V$ be a weak quantum vertex algebra and let $(W,Y_{W})$ be 
a $V$-module, namely, a module for $V$ viewed as a nonlocal vertex algebra.
For $u,v\in V$, 
$f_{i}(x)\in \C((x)),\; u^{(i)},v^{(i)}\in V$, $i=1,\dots,r$, 
if (\ref{e-Rjacobi}) holds, then we have
\begin{eqnarray}\label{emodule-wqva-Sjacobi}
& &x_{0}^{-1}\delta\left(\frac{x_{1}-x_{2}}{x_{0}}\right)Y_{W}(u,x_{1})Y_{W}(v,x_{2})\nonumber\\
& &\ \ \ \ -x_{0}^{-1}\delta\left(\frac{x_{2}-x_{1}}{-x_{0}}\right)
\sum_{i=1}^{r}f_{i}(-x_{0})Y_{W}(v^{(i)},x_{2})Y_{W}(u^{(i)},x_{1})\nonumber\\
&=&x_{2}^{-1}\delta\left(\frac{x_{1}-x_{0}}{x_{2}}\right)Y_{W}(Y(u,x_{0})v,x_{2}).
\end{eqnarray}
\el

\begin{proof} In view of Lemma \ref{lformaljacobiidentity},
for a nonlocal vertex algebra, (\ref{e-Rjacobi}) is equivalent to 
(\ref{eslocality}) for some nonnegative integer $k$
while for a module (\ref{emodule-wqva-Sjacobi}) is equivalent to (\ref{eSlocality-module-cor})
for  some nonnegative integer $l$.
Then it follows immediately from 
Corollary \ref{cSlocality-algebra-module}.
\end{proof}

The following is one of the main results of this section:

\bt{tSlocality-key}
Let $W$ be any vector space and let $S$ be
any $\S$-local subset of $\E(W)$. Then the nonlocal vertex algebra
$\<S\>$ generated by $S$ is a weak quantum vertex algebra
and $W$ is a faithful $\<S\>$-module with $Y_{W}(\alpha(x),x_{0})=\alpha(x_{0})$
for $\alpha(x)\in \<S\>$. 
\et

\begin{proof} {}From Theorem \ref{tgeneratingthem} (with $\Gamma=\{1\}$),
$\<S\>$ is a nonlocal vertex algebra with $W$ as a faithful quasi module.
By Proposition \ref{pgeneral-locality-key}, $\<S\>$ is $\S$-local.
For any $a(x),b(x)\in \<S\>$, there exists a
nonnegative integer $k$ such that
$$(x_{1}-x_{2})^{k}a(x_{1})b(x_{2})\in \Hom (W,W((x_{1},x_{2}))).$$
By Proposition \ref{pexistence}, for any $w\in W$,
there exists a nonnegative integer $l$ such that
\begin{eqnarray*}
(x_{0}+x)^{l}x_{0}^{k}(Y_{\E}(a(x),x_{0})b(x))w
&=& (x_{0}+x)^{l}x_{0}^{k}a(x_{0}+x)b(x)w\\
&=&(x_{0}+x)^{l}x_{0}^{k}Y_{W}(a(x),x_{0}+x)Y_{W}(b(x),x)w.
\end{eqnarray*}
Noticing that $Y_{W}(Y_{\E}(a(x),x_{0})b(x),x_{2})w=(Y_{\E}(a(x_{2}),x_{0})b(x_{2}))w$,
we get
\begin{eqnarray*}
& &(x_{0}+x_{2})^{l}x_{0}^{k}Y_{W}(Y_{\E}(a(x),x_{0})b(x),x_{2})w\\
&=&(x_{0}+x_{2})^{l}x_{0}^{k}Y_{W}(a(x),x_{0}+x_{2})Y_{W}(b(x),x_{2})w.
\end{eqnarray*}
This proves that $W$ is actually an $\<S\>$-module.
It follows from Corollary \ref{cSlocality-algebra-module}
and Remark \ref{rcomments} that $\<S\>$ is a weak quantum vertex algebra.
\end{proof}

\bp{prep-wqva}
Let $V$ be a weak quantum vertex algebra and let
$W$ be a vector space equipped with a linear map $Y_{W}$ from $V$ to
$\Hom (W,W((x)))\;(=\E(W))$ with $Y_{W}({\bf 1},x)=1_{W}$.
Then $(W,Y_{W})$ carries the structure of a $V$-module if and only if
$\{ Y_{W}(v,x)\;|\; v\in V\}\subset \E(W)$ is $\S$-local and closed and
$Y_{W}$ is a nonlocal vertex algebra homomorphism.
\ep

\begin{proof} If $(W,Y_{W})$ carries the structure of a $V$-module,
{}from Corollary \ref{cSlocality-algebra-module}, $\{Y_{W}(v,x)\;|\;
v\in V\}$ is $\S$-local. By Proposition \ref{pextra-need-2}, $Y_{W}(V)$
is closed and $Y_{W}$ is a homomorphism of weak quantum vertex
algebras.  Conversely, under the hypothesis that $Y_{W}(V)$ is
$\S$-local and closed, by Theorem \ref{tSlocality-key}, $Y_{W}(V)$ is
a weak quantum vertex algebra with $W$ naturally as a module.
As $Y_{W}$ is a homomorphism, $W$ is a $V$-module.
\end{proof}

\bd{d-qva-rational}
{\em A {\em quantum vertex algebra} (over $\C$) 
is a weak quantum vertex algebra $V$ 
equipped with a unitary rational quantum Yang-Baxter operator $\S(x)$ on $V$
such that for $u,v\in V$, (\ref{e-Rjacobi}) holds with
$$\sum_{i=1}^{r}v^{(i)}\otimes u^{(i)}\otimes f_{i}(x)= \S(x)(v\otimes u).$$}
\ed

Immediately from Theorem \ref{tfirst} with $W=V$ we have:
 
\bt{twqva-qva}
Let $V$ be a weak quantum vertex algebra. 
Assume that the  linear map $\pi_{3}^{V}$ is injective.
Then there exists a unitary rational quantum Yang-Baxter operator $\S$ on $V$
such that 
$(Y,\S)$ is a quantum vertex algebra and such $\S$ is uniquely determined by
the $\S$-locality.
\et

The following notion (for a vertex algebra) is due to 
Etingof and Kazhdan (see \cite{ek}):

\bd{dnondeg}
{\em A nonlocal vertex algebra $V$ is said to be {\em nondegenerate}
if for any positive integer $n$, the linear map
\begin{eqnarray}
Z_{n}&=&Y(x_{1})(1\otimes Y(x_{2}))\cdots 
(1^{\otimes (n-1)}\otimes Y(x_{n}))(1^{\otimes n}\otimes {\bf 1}):
\nonumber\\
& &V^{\otimes n}\otimes \C ((x_{1}))\cdots ((x_{n}))
\rightarrow V((x_{1}))\cdots ((x_{n}))
\end{eqnarray}
defined by
\begin{eqnarray}
Z_{n}(v_{1}\otimes \cdots \otimes v_{n}\otimes f)=
f(x_{1},\dots,x_{n}) Y(v_{1},x_{1})\cdots Y(v_{n},x_{n}){\bf 1}
\end{eqnarray}
for $v_{1},\dots,v_{n}\in V,\; f\in \C((x_{1}))\cdots ((x_{n}))$
is injective.}
\ed

If $Z_{n}$ is injective for some positive integer $n$, 
it is clear that $\pi_{n}^{V}$ is injective. Then we immediately have:

\bc{cnondeg-wqva}
Every nondegenerate weak quantum vertex algebra is a quantum vertex algebra
where the quantum Yang-Baxter operator $\S(x)$ is uniquely determined.
\ec

\section{General constructions of nonlocal vertex algebras and 
their quasi modules}

In this section we give certain general construction theorems for 
nonlocal vertex algebras and their quasi modules. 

First we formulate the following analogue
of a result of \cite{li-form}:

\bl{lvacuum-like-vector}
Let $V$ be a nonlocal vertex algebra, let $W$ be a quasi $V$-module
and let $e\in W$ be a vacuum-like vector in the sense that
\begin{eqnarray}
Y(v,x)e\in W[[x]]\;\;\;\mbox{ for }v\in V.
\end{eqnarray}
Then the linear map $\phi: V\rightarrow W;\; v\mapsto v_{-1}e$
is a $V$-homomorphism uniquely determined by the condition 
${\bf 1}\mapsto e$. 
\el

\begin{proof} For any $u,v\in V$, there exists $0\ne p(x,y)\in
\C[x,y]$ such that
$$p(x_{0}+x_{2},x_{2})Y_{W}(Y(u,x_{0})v,x_{2})e
=p(x_{0}+x_{2},x_{2})Y_{W}(u,x_{0}+x_{2})Y_{W}(v,x_{2})e.$$
As $Y_{W}(Y(u,x_{0})v,x_{2})e\in W[[x_{2}]]((x_{0}))$ and
$Y_{W}(u,x_{0}+x_{2})Y_{W}(v,x_{2})e\in W((x_{0}))[[x_{2}]]$,
we have $Y_{W}(Y(u,x_{0})v,x_{2})e
=Y_{W}(u,x_{0}+x_{2})Y_{W}(v,x_{2})e$. Then it follows immediately.
\end{proof}

Furthermore we have:

\bp{pL(-1)abstract}
Let $V$ be a nonlocal vertex algebra and 
let $(W,Y_{W})$ be a quasi $V$-module.
Suppose that $L(-1)\in \End_{\C}W$ and $e\in W$ be such that 
$L(-1)e=0$ and 
\begin{eqnarray}
[L(-1),Y_{W}(u,x)]=\frac{d}{dx}Y_{W}(u,x)
\end{eqnarray}
for $u\in U$, where $U$ is a generating subset of $V$. 
Then the linear map $\phi: V\rightarrow W;\; v\mapsto v_{-1}e$
is a $V$-homomorphism uniquely determined by the condition 
${\bf 1}\mapsto e$.
\ep

\begin{proof} In view of Lemma \ref{lvacuum-like-vector}
it suffices to prove that $e$ is a vacuum-like vector.
If we can prove
$[L(-1),Y_{W}(v,x)]=\frac{d}{dx}Y_{W}(v,x)$ for all $v\in V$,
the same argument of Proposition 3.3
of \cite{li-form} shows that $e$ is a vacuum-like vector.
Set
$$T=\{ v\in V\;|\;[L(-1),Y_{W}(v,x)]=\frac{d}{dx}Y_{W}(v,x)\}.$$ 
We must prove $T=V$. By assumption, $U\subset T$, and clearly
${\bf 1}\in T$. As $U$ generates $V$, it suffices to prove that
$T$ is a subalgebra.
Let $u,v\in T,\; w\in W$. There exists a nonzero polynomial $p(x,y)$ 
such that
\begin{eqnarray*}
& &p(x_{0}+x_{2},x_{2})Y_{W}(Y(u,x_{0})v,x_{2})w
=p(x_{0}+x_{2},x_{2})Y_{W}(u,x_{0}+x_{2})Y_{W}(v,x_{2})w,\\
& &p(x_{0}+x_{2},x_{2})Y_{W}(Y(u,x_{0})v,x_{2})L(-1)w
=p(x_{0}+x_{2},x_{2})Y_{W}(u,x_{0}+x_{2})Y_{W}(v,x_{2})L(-1)w.
\end{eqnarray*}
Then
\begin{eqnarray*}
& &p(x_{0}+x_{2},x_{2})^{2}[L(-1),Y_{W}(Y(u,x_{0})v,x_{2})]w\\
&=&p(x_{0}+x_{2},x_{2})^{2}[L(-1),Y_{W}(u,x_{0}+x_{2})Y_{W}(v,x_{2})]w\\
&=&p(x_{0}+x_{2},x_{2})^{2}\frac{\partial}{\partial x_{2}}
Y_{W}(u,x_{0}+x_{2})Y_{W}(v,x_{2})w\\
&=&\frac{\partial}{\partial x_{2}}\left(p(x_{0}+x_{2},x_{2})^{2}
Y_{W}(u,x_{0}+x_{2})Y_{W}(v,x_{2})w\right)\\
& &-2p(x_{0}+x_{2},x_{2})
\left(\frac{\partial}{\partial x_{2}}p(x_{0}+x_{2},x_{2})\right)
Y_{W}(u,x_{0}+x_{2})Y_{W}(v,x_{2})w\\
&=&\frac{\partial}{\partial x_{2}}\left(p(x_{0}+x_{2},x_{2})^{2}
Y_{W}(Y(u,x_{0})v,x_{2})w\right)\\
& &-2p(x_{0}+x_{2},x_{2})
\left(\frac{\partial}{\partial x_{2}}p(x_{0}+x_{2},x_{2})\right)
Y_{W}(Y(u,x_{0})v,x_{2})w\\
&=&p(x_{0}+x_{2},x_{2})^{2}\frac{\partial}{\partial x_{2}}
Y_{W}(Y(u,x_{0})v,x_{2})w.
\end{eqnarray*}
Consequently
$$[L(-1),Y_{W}(Y(u,x_{0})v,x_{2})]w=
\frac{\partial}{\partial x_{2}}Y_{W}(Y(u,x_{0})v,x_{2})w.$$
It follows that $Y(u,x_{0})v\in T((x_{0}))$, proving that
$K$ is a subalgebra of $V$. 
\end{proof}

The following is an analogue of a theorem of [FKRW] and [MP]:

\bt{tqva-construction}
Let $V$ be a vector space, ${\bf 1}$ a vector of $V$,
$L(-1)$ a linear operator on $V$, $U$ a subset of $V$, and 
\begin{eqnarray}
Y_{0}:U\rightarrow \Hom (V,V((x))), \; 
u\mapsto Y_{0}(u,x)=u(x)=\sum_{n\in \Z}u_{n}x^{-n-1},
\end{eqnarray}
a map. Assume that all the following conditions hold:
\begin{eqnarray}
& &L(-1){\bf 1}=0,\\
& &[L(-1),u(x)]=\frac{d}{dx}u(x),\label{eassumption-2}\\
& &u(x){\bf 1}\in V[[x]]\;\;\mbox{ and  }\;\;
\lim_{x\rightarrow 0}u(x){\bf 1}=u\;\;\;\mbox{ for }u\in U,
\label{ecreat-assump}
\end{eqnarray}
the subset $Y_{0}(U)=\{ u(x)\;|\;u\in U\}$ 
of $\E(V)$ is quasi compatible, and
\begin{eqnarray}\label{eassumption-span}
V={\rm span}\{u^{1}_{n_{1}}\cdots u^{r}_{n_{r}}{\bf 1}\;|\;
r\in \N,\; u^{i}\in U,\; n_{i}\in \Z\}.
\end{eqnarray}
In addition we assume that there exists a linear map 
$\psi$ from $V$ into $\<Y_{0}(U)\>$
such that $\psi ({\bf 1})=1_{W}$ and
\begin{eqnarray}\label{ehom-assump}
\psi (u_{n}v)=u(x)_{n}\psi (v)\;\;\;\mbox{ for }u\in U,\; v\in V,\;
n\in \Z.
\end{eqnarray}
Then $Y_{0}$ can be extended uniquely to a linear map
$Y: V\rightarrow \Hom(V,V((x)))$ such that
$(V,Y,{\bf 1})$ carries the structure of a nonlocal vertex algebra.
In the assumption, if $Y_{0}(U)$ is $\S$-local, then $V$ is a
weak quantum vertex algebra.
\et

\begin{proof} With the spanning property (\ref{eassumption-span}), 
the uniqueness is clear.
Now we must prove the existence. 
With $Y_{0}(U)$ being quasi compatible, by Theorem \ref{tgeneratingthem}, 
we have a nonlocal vertex algebra $H=\<Y_{0}(U)\>$, which is generated by $Y_{0}(U)$ inside $\E(V)$, 
and $V$ is a faithful quasi
$H$-module with $Y_{V}(\alpha(x),x_{0})=\alpha(x_{0})$.
{}From (\ref{eassumption-2}), for $u\in U$  we have
$$[L(-1),Y_{V}(u(x),x_{0})]=[L(-1),u(x_{0})]=\frac{d}{dx_{0}}u(x_{0})
=\frac{\partial}{\partial x_{0}}Y_{V}(u(x),x_{0})$$
and $L(-1){\bf 1}=0$. By Proposition \ref{pL(-1)abstract},
we have an $H$-homomorphism  $\phi$ from $H$ onto $V$ where
$$\phi (\alpha(x))
=\lim_{x_{0}\rightarrow 0}Y_{V}(\alpha(x),x_{0}){\bf 1}
=\lim_{x_{0}\rightarrow 0}\alpha(x_{0}){\bf 1}\;\;\;
\mbox{ for }\alpha(x)\in H. $$
For $u\in U$ we have
$$\phi(1_{W})={\bf 1},\ \ 
\phi(u(x))=\lim_{x_{0}\rightarrow 0}u(x_{0}){\bf 1}=u.$$
Using (\ref{ecreat-assump}) and (\ref{ehom-assump}) we have
$$\psi(u)=\psi(u_{-1}{\bf 1})=u(x)_{-1}1_{W}=u(x)\;\;\;\mbox{ for }u\in U.$$
It follows that $\phi$ and $\psi$ are inverse each other.
Then the nonlocal vertex algebra structure on $H$ can be transported 
onto $V$ so that $V$ is a nonlocal vertex algebra.
Using the fact that $\phi$ is an $H$-homomorphism we have
$$Y(u,x_{0})=\phi Y_{\E}(u(x),x_{0})\phi^{-1}
=Y_{V}(u(x),x_{0})=u(x_{0})$$
for $u\in U$. Thus, $Y$ indeed extends $Y_{0}$.

If $Y_{0}(U)$ is $\S$-local, by Proposition \ref{pgeneral-locality-key}, 
$H=\<Y_{0}(U)\>$ is a weak quantum vertex algebra. Thus $V$ is a
weak quantum vertex algebra.
\end{proof}

\br{rwhy}
{\em In contrast with the theorem in [FKRW] and [MP] (cf. \cite{ll}), 
here we have one more assumption (the last one with $\psi$). 
We mention that this assumption is necessary in this situation.
For example, let $A$ be any nonassociative algebra
with identity. Set ${\bf 1}=1$ and $L(-1)=0$. 
Define a linear map $Y: A\rightarrow \End A$ by $Y(a,x)b=ab$.
Then all the assumptions in Theorem \ref{tqva-construction} hold 
except the last one. 
We see that $(A,Y,1)$ carries the structure of a nonlocal vertex algebra 
if and only if $A$ is associative.}
\er

We also have the following analogue of a result of \cite{ll}:

\bt{twva-module}
Let $V$ be a nonlocal vertex algebra,
$U$ a generating subspace
and $W$ a vector space equipped with a linear map $Y_{W}^{0}$ from $U$
to $\Hom (W,W((x)))$. 
Assume that $\{Y_{W}^{0}(u,x)\;|\; u\in U\}\subset \E(W)$ is quasi compatible.
Then $Y_{W}^{0}$ can be extended (uniquely) to a linear map
$Y_{W}$ from $V$ to $\Hom (W,W((x)))$ such that
$(W,Y_{W})$ carries the structure of a quasi $V$-module if and only if
there exists a linear map $f$ from $V$ to $\<Y_{W}^{0}(U)\>$ such that 
$f({\bf 1})=1_{W}$ and
\begin{eqnarray}\label{esemi-hom}
f(u_{n}v)=u(x)_{n}f(v)\;\;\;\mbox{ for }u\in U,\; v\in V.
\end{eqnarray}
\et

\begin{proof} Assume the existence of $f$. 
Note that since $Y_{W}^{0}(U)$ is a quasi compatible subset of $\E(W)$, 
$\<Y_{W}^{0}(U)\>$ is a nonlocal vertex algebra
with $W$ as a quasi module (by Theorem \ref{tgeneratingthem}). 
As $U$ generates $V$,
it follows from (\ref{esemi-hom}) (cf. [LL]) that $f$ is a 
nonlocal vertex algebra homomorphism. 
Then $W$ is naturally a quasi module 
for $V$ through the map $f$. Clearly, the map $Y_{W}$ extends
$Y_{W}^{0}$.

Conversely, assume that $Y_{W}^{0}$ can be extended to 
a quasi module structure $Y_{W}$ for $V$. Set $f=Y_{W}$, a linear map
{}from $V$ to $\E(W)$. If we can prove 
$$Y_{W}(v,x)\in \<Y_{W}^{0}(U)\>\;\;\;\mbox{ for all }v\in V,$$
then by Proposition \ref{pextra-need}, 
$f$ is a linear map from $V$ to $\<Y_{W}^{0}(U)\>$
with the desired properties.
Set
$$T=\{ v\in V\;|\; Y_{W}(v,x)\in \<Y_{W}^{0}(U)\>\}.$$
Now it suffices to prove $V=T$.
Clearly, $T$ contains $U$ and ${\bf 1}$. For any $a,b\in T$,
since $(Y_{W}(a,x),Y_{W}(b,x))$ is quasi compatible 
(because $\<Y_{W}^{0}(U)\>$ is a quasi compatible set), 
by Proposition \ref{pextra-need} we have
$$Y_{W}(Y(a,x_{0})b,x)=Y_{\E}(Y_{W}(a,x),x_{0})Y_{W}(b,x)
\in \<Y_{W}^{0}(U)\>((x_{0})).$$
It follows that $T$ is closed. 
As $U$ generates $V$, we have $V=T$, completing the proof.
\end{proof}

In practical applications, the following result, which follows 
immediately {}from Theorem \ref{tclosed} and the following
Proposition \ref{pmodule-algebra-relations},
shall often be used as a companion of 
Theorems \ref{tqva-construction} and \ref{twva-module}: 

\bp{pmodule-algebra-relations-case}
Let $W$ be a vector space, let $V$ be a closed quasi compatible 
subspace of $\E(W)$, containing $1_{W}$, and let
$$u(x),v(x),u^{(1)}(x),v^{(1)}(x),\dots,u^{(r)}(x),v^{(r)}(x),c^{0}(x),\dots,c^{s}(x) \in V,\; 
f_{1}(x),\dots,f_{r}(x)\in \C((x)).$$
Suppose that
\begin{eqnarray}\label{e-cross-bracket-W}
& &(x_{1}-x_{2})^{n}u(x_{1})v(x_{2})
-(-x_{2}+x_{1})^{n}\sum_{i=1}^{r}f_{i}(x_{2}-x_{1})v^{(i)}(x_{2})u^{(i)}(x_{1})
\nonumber\\
&=&
\sum_{j=0}^{s}c^{j}(x_{2})\frac{1}{j!}\left(\frac{\partial}{\partial x_{2}}\right)^{j}
x_{2}^{-1}\delta\left(\frac{x_{1}}{x_{2}}\right)
\end{eqnarray}
for some integer $n$, then 
\begin{eqnarray}\label{e-cross-bracket-algebra-case}
& &(x_{1}-x_{2})^{n}Y_{\E}(u(x),x_{1})Y_{\E}(v(x),x_{2})\nonumber\\
& &\ \ \ \ \ \ 
-(-x_{2}+x_{1})^{n}\sum_{i=1}^{r}f_{i}(x_{2}-x_{1})Y_{\E}(v^{(i)}(x),x_{2})Y_{\E}(u^{(i)}(x),x_{1})
\nonumber\\
&=&\sum_{j=0}^{s}Y_{\E}(c^{j}(x),x_{2})\frac{1}{j!}\left(\frac{\partial}{\partial x_{2}}\right)^{j}
x_{2}^{-1}\delta\left(\frac{x_{1}}{x_{2}}\right).
\end{eqnarray}
\ep

The following is a generalization of
Corollary \ref{cSlocality-algebra-module}:

\bp{pmodule-algebra-relations}
Let $V$ be a nonlocal vertex algebra, let $(W,Y_{W})$ be a $V$-module and let
$$n\in \Z,\; u,v,u^{(1)},v^{(1)},\dots,u^{(r)},v^{(r)},c^{0},\dots,c^{s} \in V,\; 
f_{1}(x),\dots,f_{r}(x)\in \C((x)).$$
If
\begin{eqnarray}\label{e-cross-bracket}
& &(x_{1}-x_{2})^{n}Y(u,x_{1})Y(v,x_{2})
-(-x_{2}+x_{1})^{n}\sum_{i=1}^{r}f_{i}(x_{2}-x_{1})Y(v^{(i)},x_{2})Y(u^{(i)},x_{1})
\nonumber\\
&=&
\sum_{j=0}^{s}Y(c^{j},x_{2})\frac{1}{j!}\left(\frac{\partial}{\partial x_{2}}\right)^{j}
x_{2}^{-1}\delta\left(\frac{x_{1}}{x_{2}}\right),
\end{eqnarray}
then 
\begin{eqnarray}\label{e-cross-bracket-module}
& &(x_{1}-x_{2})^{n}Y_{W}(u,x_{1})Y_{W}(v,x_{2})
-(-x_{2}+x_{1})^{n}\sum_{i=1}^{r}f_{i}(x_{2}-x_{1})Y_{W}(v^{(i)},x_{2})Y_{W}(u^{(i)},x_{1})
\nonumber\\
&=&\sum_{j=0}^{s}Y_{W}(c^{j},x_{2})\frac{1}{j!}\left(\frac{\partial}{\partial x_{2}}\right)^{j}
x_{2}^{-1}\delta\left(\frac{x_{1}}{x_{2}}\right).
\end{eqnarray}
If  $(W,Y_{W})$ is faithful, the converse is also true.
\ep

\begin{proof} Assume (\ref{e-cross-bracket}). 
Then for any nonnegative integer $k$ with $k\ge n,s+1$ we have
\begin{eqnarray}\label{e-Slocality-property-proof}
(x_{1}-x_{2})^{k}Y(u,x_{1})Y(v,x_{2})
=(x_{1}-x_{2})^{k}\sum_{i=1}^{r}f_{i}(x_{2}-x_{1})
Y(v^{(i)},x_{2})Y(u^{(i)},x_{1}).
\end{eqnarray}
In view of Lemma \ref{lformaljacobiidentity} we have
\begin{eqnarray}\label{etemp-jacobi}
& &x_{0}^{-1}\delta\left(\frac{x_{1}-x_{2}}{x_{0}}\right)Y(u,x_{1})Y(v,x_{2})\nonumber\\
& &\ \ \ \ -x_{0}^{-1}\delta\left(\frac{x_{2}-x_{1}}{-x_{0}}\right)
\sum_{i=1}^{r}f_{i}(x_{2}-x_{1})Y(v^{(i)},x_{2})Y(u^{(i)},x_{1})\nonumber\\
&=&x_{2}^{-1}\delta\left(\frac{x_{1}-x_{0}}{x_{2}}\right)Y(Y(u,x_{0})v,x_{2}),
\end{eqnarray}
which gives
\begin{eqnarray}\label{etemp-cross}
& &(x_{1}-x_{2})^{n}Y(u,x_{1})Y(v,x_{2})
-(-x_{2}+x_{1})^{n}\sum_{i=1}^{r}f_{i}(x_{2}-x_{1})Y(v^{(i)},x_{2})Y(u^{(i)},x_{1})
\nonumber\\
&=&
\sum_{j\ge 0}Y(u_{n+j}v,x_{2})\frac{1}{j!}\left(\frac{\partial}{\partial x_{2}}\right)^{j}
x_{2}^{-1}\delta\left(\frac{x_{1}}{x_{2}}\right).
\end{eqnarray}
Combining this with (\ref{e-cross-bracket})
we get $c^{j}=u_{j+n}v$ for $j=0,\dots,s$ and $u_{n+j}v=0$ for $j>s$,
as $(V,Y)$ being faithful.

Also, with (\ref{e-Slocality-property-proof}), by Corollary \ref{cSlocality-algebra-module}
we have
\begin{eqnarray}\label{e-Slocality-property-proof-module}
(x_{1}-x_{2})^{l}Y_{W}(u,x_{1})Y_{W}(v,x_{2})
=(x_{1}-x_{2})^{l}\sum_{i=1}^{r}f_{i}(x_{2}-x_{1})
Y_{W}(v^{(i)},x_{2})Y_{W}(u^{(i)},x_{1})
\end{eqnarray}
for some nonnegative integer $l$. In view of Lemma \ref{lformaljacobiidentity} we have
\begin{eqnarray}\label{eSjacobi-identity-module-proof}
& &x_{0}^{-1}\delta\left(\frac{x_{1}-x_{2}}{x_{0}}\right)Y_{W}(u,x_{1})Y_{W}(v,x_{2})\nonumber\\
& &\ \ \ \ -x_{0}^{-1}\delta\left(\frac{x_{2}-x_{1}}{-x_{0}}\right)
\sum_{i=1}^{r}f_{i}(x_{2}-x_{1})Y_{W}(v^{(i)},x_{2})Y_{W}(u^{(i)},x_{1})\nonumber\\
&=&x_{2}^{-1}\delta\left(\frac{x_{1}-x_{0}}{x_{2}}\right)Y_{W}(Y(u,x_{0})v,x_{2}),
\end{eqnarray}
which gives 
\begin{eqnarray}\label{ecross-product-module-proof}
& &(x_{1}-x_{2})^{n}Y_{W}(u,x_{1})Y_{W}(v,x_{2})
-(-x_{2}+x_{1})^{n}\sum_{i=1}^{r}f_{i}(x_{2}-x_{1})Y_{W}(v^{(i)},x_{2})Y_{W}(u^{(i)},x_{1})
\nonumber\\
&=&\sum_{j=0}^{s}Y_{W}(u_{j+n}v,x_{2})\frac{1}{j!}\left(\frac{\partial}{\partial x_{2}}\right)^{j}
x_{2}^{-1}\delta\left(\frac{x_{1}}{x_{2}}\right)\nonumber\\
&=&\sum_{j=0}^{s}Y_{W}(c^{j},x_{2})\frac{1}{j!}\left(\frac{\partial}{\partial x_{2}}\right)^{j}
x_{2}^{-1}\delta\left(\frac{x_{1}}{x_{2}}\right),
\end{eqnarray}
proving (\ref{e-cross-bracket-module}).

On the other hand, assume that 
$(W,Y_{W})$ is a faithful $V$-module and (\ref{e-cross-bracket-module}) holds.
Just as in the first part, first we  have (\ref{e-Slocality-property-proof-module})
for any nonnegative integer $l\ge n,s+1$, and then we have
(\ref{eSjacobi-identity-module-proof}), which gives
the first equality in (\ref{ecross-product-module-proof}).
Using (\ref{e-cross-bracket-module}) we have the second equality 
in (\ref{ecross-product-module-proof}).
As $(W,Y_{W})$ is faithful, we get
$u_{j+n}v=c^{j}$ for $j=0,\dots,s$ and $u_{n+j}v=0$ for $j>s$.
{}From (\ref{e-Slocality-property-proof-module}), 
by Corollary \ref{cSlocality-algebra-module}
(\ref{e-Slocality-property-proof}) holds for some nonnegative integer $k$,
then (\ref{etemp-jacobi}) holds, and then
(\ref{etemp-cross}) holds. Finally we obtain (\ref{e-cross-bracket-W}).
\end{proof}

In the following we present an example to illustrate how
the general construction theorems can be used.
Let $E=K\oplus \g$ be a semidirect product Lie algebra 
with $K$ an ideal and with $\g$ a Lie subalgebra.  
Suppose that $\<\cdot,\cdot\>$
is a symmetric invariant bilinear form on $K$ such that
$\< [a,u],v\>=-\<u,[a,v]\>$ for $u,v\in K,\; a\in \g$.
Extend $\<\cdot,\cdot\>$ to a bilinear form on $E$ by
$\<u+a,v+b\>=\<u,v\>$ for $u,v\in K,\; a,b\in \g$.  
Then $\<\cdot,\cdot\>$ is symmetric and invariant.
Associated to $(E,\<\cdot,\cdot\>)$ we have the affine Lie
algebra $\hat{E}=E\otimes \C[t,t^{-1}]\oplus C c$.  
Consider the following Lie subalgebra
\begin{eqnarray}
{\mathcal{L}}
=K\otimes \C[t,t^{-1}]\oplus \C c\oplus (\g\otimes t^{-1}\C[t^{-1}]),
\end{eqnarray}
which is the semidirect product Lie algebra of $\hat{K}$ 
with $\g\otimes t^{-1}\C[t^{-1}]$.
Let $\ell$ be any complex number. 
Let $\C_{\ell}$ be the $1$-dimensional module
for $K[t]\oplus \C c$ with $K[t]$ acting as zero 
and with $c$ acting as scalar $\ell$.
Form the induced $\mathcal{L}$-module
\begin{eqnarray}
V_{\mathcal{L}}(\ell,0)
=U({\mathcal{L}})\otimes_{U(\g[t]\oplus \C c)} \C_{\ell}.
\end{eqnarray}
In view of the P-B-W theorem, we have
\begin{eqnarray}
V_{\mathcal{L}}(\ell,0)=U(E\otimes t^{-1}\C[t^{-1}]),
\end{eqnarray}
as a vector space. Set ${\bf 1}=1\otimes 1\in V_{\mathcal{L}}(\ell,0)$.
Embed $E$ into $V_{\mathcal{L}}(\ell,0)$
through the map $v\mapsto v(-1){\bf 1}$, where for $v\in E,\; n\in
\Z$, $v(n)$ denotes the corresponding operator of $v\otimes t^{n}$.

\bp{psemiproduct}
There exists a unique weak quantum vertex algebra structure
on $V_{\mathcal{L}}(\ell,0)$ with ${\bf 1}$ as the vacuum vector
such that for $u\in K,\; a\in \g$,
\begin{eqnarray}
Y(u,x)=u(x)=\sum_{n\in \Z}u(n)x^{-n-1},\;\;\ \ 
Y(a,x)=a(x)^{+}=\sum_{n\ge 1}a(-n)x^{n-1}.
\end{eqnarray}
\ep

\begin{proof} We shall apply Theorem \ref{tqva-construction}.
Note that $1\otimes d/dt$ is a derivation of ${\mathcal{L}}$.
It follows (cf. \cite{ll}, Remark 6.2.9) that 
there exists a (unique) operator 
$L(-1)$ on $V_{\mathcal{L}}(\ell,0)$ such that $L(-1){\bf 1}=0$ and
$$[L(-1),u(x)]=(d/dx)u(x),\;\;\;[L(-1),a(x)^{+}]=(d/dx)a(x)^{+}$$
for $u\in K,\; a\in \g$. 
Set $U=\{1, u(x),a(x)^{+}\;|\; u\in K,\; a\in \g\}.$
By a straightforward calculation (cf. \cite{efk}) we have
\begin{eqnarray*}
& &[u(x_{1}),v(x_{2})]
=[u,v](x_{2})x_{1}^{-1}\delta\left(\frac{x_{2}}{x_{1}}\right)
+\ell \<u,v\> \frac{\partial}{\partial
x_{2}}x_{1}^{-1}\delta\left(\frac{x_{2}}{x_{1}}\right),\\
& &[a(x_{1})^{+},u(x_{2})]=[a,u](x_{2})(x_{2}-x_{1})^{-1},\\
& &[a(x_{1})^{+},b(x_{2})^{+}]=
(x_{2}-x_{1})^{-1}([a,b](x_{2})^{+}-[a,b](x_{1})^{+})
\end{eqnarray*}
for $u,v\in K,\; a,b\in \g$.
It follows that $U$ is $\S$-local. 
By Theorem \ref{tSlocality-key}, $U$ generates 
a weak quantum vertex algebra $\<U\>$. It follows from
Proposition \ref{pmodule-algebra-relations-case} that $\<U\>$ is naturally 
an ${\mathcal{L}}$-module of level $\ell$. Then there exists 
a unique ${\mathcal{L}}$-homomorphism $\psi$ from $V$ to $\<U\>$, sending 
$1_{W}$ to ${\bf 1}$. Now, we have all the conditions 
assumed in Theorem \ref{tqva-construction}.
\end{proof}

\end{document}